\documentclass[hidelinks,onefignum,onetabnum]{siamart250211}

\usepackage{amsfonts}
\usepackage{amsmath}
\usepackage{lmodern}
\usepackage{mathtools}
\usepackage[T1]{fontenc}
\usepackage{graphicx}
\usepackage{epstopdf}
\usepackage{algorithmic}
\usepackage{booktabs}
\usepackage[T1]{fontenc}
\usepackage{tikz,caption}
\usepackage{subcaption}
\usepackage{pgfplots}
\pgfplotsset{compat=newest}
\pgfplotsset{plot coordinates/math parser=false}
\usepackage{url}
\usepackage{color}
\usetikzlibrary{plotmarks} 
\usetikzlibrary{external}
\usepgfplotslibrary{external} 
\tikzset{external/mode=graphics if exists}
\usepackage{cleveref}
\usepackage{multirow}
\usepackage{comment}
\usetikzlibrary{patterns}

\newlength\figureheight
\newlength\figurewidth 

\ifpdf
\DeclareGraphicsExtensions{.eps,.pdf,.png,.jpg}
\else
\DeclareGraphicsExtensions{.eps}
\fi

\usepackage{xcolor}

\newcommand{\R}{\mathbb{R}}

\newcommand{\D}{\,\mathrm{d}}

\def\amen{{{\sc AMEn}}}
\def\ttgmres{{{\sc TT-GMRES}}}

\numberwithin{theorem}{section}
\newsiamremark{remark}{Remark}
\newsiamremark{hypothesis}{Hypothesis}
\crefname{hypothesis}{Hypothesis}{Hypotheses}
\newsiamthm{claim}{Claim}
\newsiamremark{fact}{Fact}
\crefname{fact}{Fact}{Facts}

\headers{A Low-Rank tensor framework for THB-Splines}{T.-C.~Riemer and M.~Stoll}

\title{A Low-Rank tensor framework for THB-Splines}

\author{Tom-Christian Riemer\thanks{Technische Universit\"at Chemnitz, Department of Mathematics, Chair of Scientific Computing, 09107 Chemnitz, Germany, \email{tom-christian.riemer@mathematik.tu-chemnitz.de}}
	\and
	Martin Stoll\thanks{Technische Universit\"at Chemnitz, Department of Mathematics, Chair of Scientific Computing, 09107 Chemnitz, Germany, \email{martin.stoll@mathematik.tu-chemnitz.de}}}

\usepackage{amsopn}

\begin{document}

\maketitle

\begin{abstract}
	We introduce a low-rank framework for adaptive isogeometric analysis with truncated hierarchical B-splines (THB-splines) that targets the main bottleneck of local refinement: memory- and time-intensive matrix assembly once the global tensor-product structure is lost. The method interpolates geometry-induced weight and source terms in separable spline spaces and computes their tensor-train (TT) representations via the alternating minimal energy (AMEn) solver, enabling level-wise assembly of system operators using univariate quadrature. To recover separability in the adaptive setting, we reduce the active basis to tensor-product domains and partition active/non-active cells into a small number of Cartesian cuboids, so each contributes a Kronecker factor that is accumulated and rounded in TT. We realize the two-scale relation with truncation in low rank and assemble the global hierarchical operators in a block TT format suitable for iterative solvers. A prototype \textsc{MATLAB} implementation built on the \textsc{GeoPDEs} package and the \textsc{TT-Toolbox} demonstrates that, for model problems with moderately complex refinement regions, the approach reduces memory footprint and assembly time while maintaining accuracy; we also discuss limitations when ranks grow with geometric or refinement complexity. This framework advances scalable adaptive IgA with THB-splines, particularly in three dimensions.
\end{abstract}

\begin{keywords}
	isogeometric analysis, tensor-train format, low-rank decompositions, low-rank methods, addaptivity, hierarchical B-splines, truncated hierarchical B-splines
\end{keywords}

\begin{AMS}
	65F10, 65F50, 15A69, 93C20 
\end{AMS}

\section{Motivation}

Isogeometric Analysis (IgA) is a discretization technique used for approximating solutions to partial differential equations (PDEs) defined on a given geometry $\Omega$. It was introduced by Hughes, Cottrell, and Bazilevs in 2005 (see~\cite{CAD}). In Isogeometric Analysis, the computational domain $\Omega$ and the solution space for solving the PDE using a Galerkin approach~\cite{strang} are parameterized by the same spline functions, typically B-splines or NURBS (Non-Uniform Rational B-Splines). These basis functions are globally defined and have overlapping supports depending on their degrees. As such, these discretizations have higher computational complexity, increasing exponentially with respect to the dimension of the problem~\cite{Mantzaflaris_space_time}, but they also allow for the relatively easy approximation of domains that are difficult to treat with traditional finite element methods. In recent years, low-rank tensor methods have emerged as a powerful response: by exposing and exploiting separability and Kronecker product structure in the interpolation of the integrands and operators, they reduce both the cost of matrix assembly and the memory footprint of subsequent linear algebra. This has made low-rank approaches one of the major research interests in IgA, which tries to find strategies to overcome this complexity drawback~\cite{Riemer2025, BuengerDolgovStoll:2020, montardini2023lowsp, montardini2023low, montardini2024lowrank, angelos1, angelos2, Antolin, Hughes, pan2020fast, pan2021efficient}.

At the same time, adaptivity has also become one of the most active lines of research in the IgA community~\cite{Buffa2022AdaptiveIGA, Cesare2019}. Local refinement in IgA can be achieved by using a hierarchical B-spline basis~\cite{VUONG20113554} but the partition of unity property gets lost and the overlapping basis functions lead to very dense system matrices. Truncated hierarchical B-Splines (short THB-splines), presented in~\cite{GiannelliJuettlerSpeleers:2012}, restore the partition of unity property, reduce excessive overlap across levels, and typically yields sparser discrete operators than hierarchical B-splines, which is crucial when adaptivity targets small regions of interest. The literature now offers both the foundational constructions and practical algorithms and data structures that make adaptive IgA with THB-splines accessible~\cite{GIANNELLI2016337, GarauVazquez:2018, 10.1007/978-3-642-54382-1_18, PAN2021113278}.

Nevertheless, even with THB-splines, the resulting linear systems are very memory-intensive, and solving them becomes numerically complex even for a few refinement levels. With this work, we want to present a way to implement local refinement in a memory-efficient and computationally efficient manner using low-rank methods, which until now could only be achieved by exploiting low-rank multi-patch approaches, where different refinements can be used for different patches in the case of not fully matching patches~\cite{Riemer2025}. The proposed low-rank approach is tailored to adaptive IgA with THB-splines for the model problem
\begin{equation}
	\label{eq:forward}
	\begin{aligned}
		-\Delta y &= f \quad \mbox{ in } \Omega, \\
		y & = 0  \quad \mbox{ on } \partial \Omega.
	\end{aligned}
\end{equation}
with the goal of fast and memory-efficient assembly and solution. We achieve this by following the idea of Mantzaflaris et al.~\cite{angelos1,angelos2} of using a low-rank tensor method, which exploits the tensor structure of the basis functions and separates the variables of the integrals. As a result the system matrices of B-splines without adaptivity are then approximated to high accuracy by a sum of Kronecker products of smaller matrices, which are assembled via univariate integration. We here rely on the method of~\cite{BuengerDolgovStoll:2020} where the assembly is carried out using an interpolation step and a low-rank representation of the resulting coefficient tensor. The authors there combine the low-rank method of Mantzaflaris et al.\ with low-rank tensor-train (TT) calculations~\cite{osel-tt-2011, DoOs-dmrg-solve-2011}. Exploiting the tensor product nature of the arising interpolation, we can calculate a low-rank TT approximation without prior assembly of the full coefficient tensor by means of the Alternating Minimal Energy (AMEn) method~\cite{amen}. For the assembly of hierarchical system matrices, we follow the heuristic in~\cite{GarauVazquez:2018}, whereby the key challenge is that local refinement breaks global tensor-product structure. Our contribution is to recover local separability where it exists (e.g., within unions of cells or unions of B-splines) and to propagate it through the THB machinery: level-wise assembly on active regions, efficient realization of the two-scale relation with truncation, and block-structured accumulation across levels, all in TT format. In particular, we show that one can detect cuboid-like subsets of the active mesh where separable quadrature applies, accumulate the resulting rank-one Kronecker factors in TT, and keep ranks under control by systematic rounding. The outcome is a linear system that can be assembled and solved in low rank in certain scenarios, with complexity governed by interpolation ranks and the local refinement pattern rather than by the full number of degrees of freedom.

The paper is structured as follows: In \cref{section:Preliminaries}, we set notation and background, where we introduce in \cref{subsection:Low-rank_tensor_format} the low-rank tensor theory and describe in \cref{subsection:hierarchical_setup} the hierarchical THB-spline setup used throughout. In \cref{section:assembling_linear_system}, we present how to assemble a low-rank linear system with THB-splines: we present low-rank interpolation of weight and source terms to expose separability in \cref{subsection:Low-Rank_interpolation}, we show how a domain-wise integration can be performed efficiently to form levelwise mass matrices in \cref{subsection:domain_wise_integration}. We realize the low-rank assembly of the truncation operation in \cref{subsection:Low-rank_Truncation_Operator} and we present the assembly of the the global hierarchical system tensor in \cref{subsection:Gobal_Hierarchical_Mass_Tensor}. In \cref{section:Hierarchical_Low-Rank_solver} we present a method for solving the resulting linear system and the used preconditioners. In \cref{section:numerics}, we present the results of our numerical experiments. In \cref{section:conclusion}, we conclude and give an outlook to future steps.

\section{Preliminaries} 
\label{section:Preliminaries}
\subsection{Low-rank tensor format}
\label{subsection:Low-rank_tensor_format}
The most well-known technique for low-rank approximations is the singular value decomposition, illustrated for a matrix $W \in \R^{m \times n}$ as
\begin{equation}
	\label{eq:SVD}
	W = U \Sigma V^{\top} \approx \sum_{r=1}^R u_r \sigma_r v_r^{\top} = \sum_{r=1}^R \left( u_r \sqrt{\sigma_r} \right) \otimes \left( v_r \sqrt{\sigma_r} \right).
\end{equation}
with $U\in \R^{m \times m}$, $V \in \R^{n \times n}$ with their columns denoted by $u_r$ and $v_r$, and $\Sigma \in \R^{m \times n}$ is the rectangular matrix holding the ordered singular values $\sigma_i$, $i=1,\ldots,\min \left( \left\{ m, n \right\} \right)$ on its main diagonal. The best low-rank approximation is obtained by the truncated SVD where we truncate all singular values below some given threshold resulting in a rank-$R$ approximation, where $R$ is the number of used singular values and therefore the number of summands in \cref{eq:SVD}.

In the high-dimensional case we need low-rank tensor approximations of a $D$-dimensional tensor. Such approximations are given by, e.g., the higher-order singular value decomposition (HOSVD)~\cite{mlsvd} or a canonical polyadic (CP) decomposition~\cite{CPD}. However, the approximation problem in the CP format is typically ill-posed~\cite{desilva-2008} and might be numerically unstable. The HOSVD (also known as the Tucker format) still contains the curse of dimensionality as it relies on a core tensor of the original tensor's dimension. We switch to the more robust tensor-train (TT) decomposition~\cite{osel-tt-2011} in this paper, given the availability of appropriate methods within a robust software framework.

A tensor $\mathbf{W} \in \mathbb{R}^{\left(n^{(1)}, \ldots, n^{(D)} \right)}$ is given in the TT format if its entries can be written as
\begin{equation}
	\label{eq:TT_format}
	\begin{aligned}
		\mathbf{W}_{\mathbf{i}} & = \mathbf{W}_{\left(i^{(1)},\ldots,i^{(D)} \right)} = \mathbf{W}^{\left(1\right)}_{\left( :, i^{(1)}, : \right)} \cdots \mathbf{W}^{\left( D \right)}_{\left( :, i^{(D)}, : \right)} =  \\
		& = \mathbf{W}^{\left( 1 \right)} \left( i^{(1)} \right) \cdots \mathbf{W}^{\left( D \right)} \left( i^{(D)} \right) \in \mathbb{R}
	\end{aligned}
\end{equation}
where $\mathbf{i} = \left(i^{(1)},\ldots,i^{(D)} \right) \in \bigtimes^D_{d = 1} \left\{ 1, \ldots, n^{(d)} \right\}$ is a multi-index, $\mathbf{W}^{\left(d\right)} \in \mathbb{R}^{\left(R^{(d-1)}, n^{(d)}, R^{(d)} \right)}$ are the so-called TT cores, which can be understood as parameter dependent matrices $\mathbf{W}^{\left(d\right)} \left( i^{(d)} \right) = \mathbf{W}^{\left( d \right)}_{\left( :, i^{(d)}, : \right)}$, $i^{(d)} = 1, \ldots, n^{(d)}$, of size $R^{(d-1)}\times R^{(d)}$ with $R^{(0)} = R^{(D)} = 1$~\cite{osel-tt-2011}. In the following, we use $:$ in the index position of dimension $d$ to denote the set of all indices in this dimension. The TT format can be rewritten into a CP decomposition as
\begin{equation*} 
	\mathbf{W} = \sum_{r^{(1)}=1}^{R^{(1)}}\cdots \sum_{r^{(D-1)}=1}^{R^{(D-1)}} \bigotimes_{d=1}^D \mathbf{W}^{\left(d\right)}_{\left( r^{(d-1)}, :, r^{(d)} \right)} = \sum_{r = 1}^{R} \bigotimes_{d=1}^D w^{\left(d\right)}_{r}
\end{equation*} 
with $R = R^{(1)} \cdot \ldots \cdot R^{(D-1)}$. The mode sizes $R^{(d-1)}$ and $R^{(d)}$ of a TT core $\mathbf{W}^{\left(d\right)} \in \mathbb{R}^{\left(R^{(d-1)}, n^{(d)}, R^{(d)} \right)}$ are the so-called TT ranks. In general, one is always interested in approximating a tensor $\mathbf{W}$ as best as possible with the TT format \cref{eq:TT_format} while keeping the TT ranks of each TT core small. If this succeeds, one speaks of a low-rank setup or method.

\subsection{Hierarchical setup}
\label{subsection:hierarchical_setup}
Let $\Omega_0 = \left[ 0, 1 \right]^3$ be the unit cube and 
\begin{equation}
	\label{eq:nested_spline_spaces}
	V_0 \subset V_1 \subset \ldots \subset V_{L-1}
\end{equation}
a sequence of nested $3$-variate tensor product B-spline spaces of the same degree $p = \left[ p^{(1)}, p^{(3)}, p^{(3)} \right] \in \mathbb{N}^3$ in each dimension defined on $\Omega_0$. In the following $\ell = 0, \ldots, L-1$ denotes the level of the corresponding variable. Each spline space $V_{\ell}$ is spanned by a tensor product B-spline basis $\mathcal{B}_{\ell} = \mathcal{B}^{(1)}_{\ell} \otimes \mathcal{B}^{(2)}_{\ell} \otimes \mathcal{B}^{(3)}_{\ell}$, where the splines of each univariate basis $\mathcal{B}^{(d)}_{\ell} = \left\{ \beta^{(d)}_{\ell, i} \colon \left[ 0, 1 \right] \to \left[ 0, 1 \right] \colon \, i = 1, \ldots, n^{(d)}_{\ell} \right\}$ are uniquely defined by the degree $p^{(d)}$ and an open knot vector $\xi^{(d)}_{\ell} = \left\{ \xi^{(d)}_{\ell, 1}, \ldots, \xi^{(d)}_{\ell, n^{(d)}_{\ell}+p^{(d)}+1} \right\}$ with
\begin{equation*} 
	\xi^{(d)}_{\ell, 1} = \cdots = \xi^{(d)}_{\ell, p^{(d)}+1}  < \xi^{(d)}_{\ell, p^{(d)}+2} \leq \cdots \leq \xi^{(d)}_{\ell, n^{(d)}_{\ell}} < \xi^{(d)}_{\ell, n^{(d)}_{\ell}+1} = \cdots = \xi^{(d)}_{\ell, n^{(d)}_{\ell}+p^{(d)}+1}
\end{equation*}
with $\xi^{(d)}_{\ell, 1} = 0$ and $\xi^{(d)}_{\ell, n^{(d)}_{\ell}+p^{(d)}+1} = 1$ using the Cox-de Boor recursion formula~\cite{DEBOOR197250, cox:b_splines}. The size of $\mathcal{B}_{\ell}$ is $n_{\ell} = \left( n^{(1)}_{\ell}, n^{(2)}_{\ell}, n^{(3)}_{\ell} \right)$ and $V_{\ell}$ is therefore spanned by $n^{(1)}_{\ell} n^{(2)}_{\ell} n^{(3)}_{\ell}$ many splines. We point out that since the B-spline spaces are nested, see \cref{eq:nested_spline_spaces}, the same holds true for the bases $\mathcal{B}_{\ell}$ as for the univariate knot vectors $\xi^{(d)}_{\ell}$, $\ell = 0, \ldots, L-1$, $d = 1, 2, 3$, and we assume that the nestedness of the knot vectors emerges from dyadic refinement from one level to the next. Since each B-spline basis $\mathcal{B}_{\ell}$ is defined by means of the tensor product, its basis functions $\beta_{\ell, \mathbf{i}} \in \mathcal{B}_{\ell}$ are defined by 
\begin{equation}
	\label{eq:multivariate_B_spline}
	\beta_{\ell, \mathbf{i}} \left(x\right) = \beta^{(1)}_{\ell, i^{(1)}} \left(x^{(1)}\right) \beta^{(2)}_{\ell, i^{(2)}} \left(x^{(2)}\right) \beta^{(3)}_{\ell, i^{(3)}} \left(x^{(3)}\right),
\end{equation}
with a multi-index $\mathbf{i} = \left( i^{(1)}, i^{(2)}, i^{(3)}\right) \in \mathcal{I}_{\mathcal{B}_{\ell}} = \bigtimes^3_{d = 1} \left\{ 1, \ldots, n^{(d)}_{\ell} \right\}$ and univariate basis functions $\beta^{(d)}_{\ell, i^{(d)}} \in \mathcal{B}^{(d)}_{\ell}$, $d = 1, 2, 3$. By setting $B^{(d)}_{\ell} = \left[ \beta^{(d)}_{\ell, 1}, \ldots, \beta^{(d)}_{\ell, n^{(d)}_{\ell}} \right]^\top \in \mathbb{R}^{n^{(d)}_{\ell}}$ we are able to represent all splines of level $\ell$ by using the tensor or Kronecker product
\begin{equation}
	\label{eq:splines_tensor_product}
	\mathbf{B}_{\ell} \left( x \right) =  B^{(1)}_{\ell} \left(x^{(1)}\right) \otimes B^{(2)}_{\ell} \left(x^{(2)}\right) \otimes B^{(3)}_{\ell} \left(x^{(3)}\right) \in \mathbb{R}^{\left( n^{(1)}_{\ell}, n^{(2)}_{\ell}, n^{(3)}_{\ell} \right)}.
\end{equation}
We point out that in the following we consider the Kronecker product and tensor product as equivalent operations, and thus \cref{eq:splines_tensor_product} can also be regarded as a vector in the space $\mathbb{R}^{n^{(1)}_{\ell} n^{(2)}_{\ell} n^{(3)}_{\ell}}$. 

\begin{figure}
	\begin{center}
		\begin{tikzpicture}[scale=3, line cap=round, line join=round]
			\fill[gray!35] (0,0) rectangle (0.125,0.125);
			\fill[gray!35] (0.125,0) rectangle (0.25,0.125);
			\fill[gray!35] (0,0.125) rectangle (0.125,0.25);
			\fill[gray!35] (0.125,0.125) rectangle (0.25,0.25);
			
			\draw[very thick] (0,0) rectangle (1,1);
			
			\draw[very thick] (0.5,0) -- (0.5,1);
			\draw[very thick] (0,0.5) -- (1,0.5);
			
			\draw[very thick] (0.25,0) -- (0.25,0.5);
			\draw[very thick] (0,0.25) -- (0.5,0.25);
			
			\draw[very thick] (0.125,0) -- (0.125,0.25);
			\draw[very thick] (0,0.125) -- (0.25,0.125);
			
			\draw[very thin, dotted] (0.25,0.5) -- (0.25,1);
			\draw[very thin, dotted] (0,0.75) -- (0.5,0.75);
			
			\draw[very thin, dotted] (0.75,0.5) -- (0.75,1);
			\draw[very thin, dotted] (0.5,0.75) -- (1,0.75);
			
			\draw[very thin, dotted] (0.75,0) -- (0.75,0.5);
			\draw[very thin, dotted] (0.5,0.25) -- (1,0.25);
		\end{tikzpicture}  
		\hspace{20mm} 
		\begin{tikzpicture}[scale=3, line cap=round, line join=round]
			\fill[gray!35] (0,0.5) rectangle (0.5,1);   
			\fill[gray!35] (0.5,0.5) rectangle (1,1);   
			\fill[gray!35] (0.5,0)   rectangle (1,0.5); 
			
			\fill[gray!35] (0,0)       rectangle (0.125,0.125);
			\fill[gray!35] (0.125,0)   rectangle (0.25,0.125);
			\fill[gray!35] (0,0.125)   rectangle (0.125,0.25);
			\fill[gray!35] (0.125,0.125) rectangle (0.25,0.25);
			
			\draw[very thick] (0,0) rectangle (1,1);
			
			\draw[very thick] (0.5,0) -- (0.5,1);
			\draw[very thick] (0,0.5) -- (1,0.5);
			
			\draw[very thick] (0.25,0) -- (0.25,0.5);
			\draw[very thick] (0,0.25) -- (0.5,0.25);
			
			\draw[very thick] (0.125,0) -- (0.125,0.25);
			\draw[very thick] (0,0.125) -- (0.25,0.125);
			
			\draw[very thin, dotted] (0.25,0.5) -- (0.25,1);
			\draw[very thin, dotted] (0,0.75) -- (0.5,0.75);
			
			\draw[very thin, dotted] (0.75,0.5) -- (0.75,1);
			\draw[very thin, dotted] (0.5,0.75) -- (1,0.75);
			
			\draw[very thin, dotted] (0.75,0) -- (0.75,0.5);
			\draw[very thin, dotted] (0.5,0.25) -- (1,0.25);
		\end{tikzpicture}    
	\end{center}
	\caption{A hierarchical mesh for $L = 3$. On the left the deactivated cells and on the right the non-active cells of level $\ell = 1$ are colored.}
	\label{fig:non_active_cells}
\end{figure}
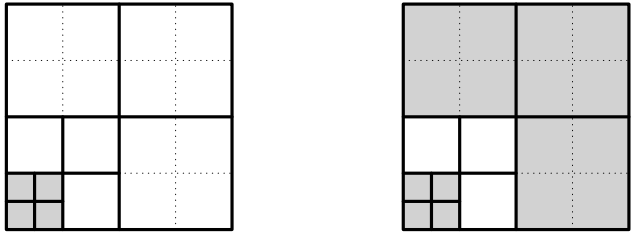

We define the set of breaking points of level $\ell$ for dimension $d$ as the set of non-repeated knot values, i.e.~$\hat{\xi}^{(d)}_{\ell} = \left\{ \hat{\xi}^{(d)}_{\ell,1}, \ldots, \hat{\xi}^{(d)}_{\ell,m^{(d)}_{\ell} + 1} \right\} \subset \xi^{(d)}_{\ell}$, which enables us to define a Cartesian mesh for each level $\ell = 1, \ldots, L-1$
\begin{equation*}
	G_{\ell} = \left\{ q_{\ell, \mathbf{j}} = q^{(1)}_{\ell, j^{(1)}} \times q^{(2)}_{\ell, j^{(2)}} \times q^{(3)}_{\ell, j^{(3)}} \colon \, j^{(d)} = 1, \ldots, m^{(d)}_{\ell}, \, d = 1, 2, 3 \right\}
\end{equation*}
where $q^{(d)}_{\ell, j^{(d)}} = \left[ \hat{\xi}^{(d)}_{\ell,j^{(d)}}, \hat{\xi}^{(d)}_{\ell,j^{(d)}+1} \right]$ and we say that $q \in G_{\ell}$ is a cell of level $\ell$. We point out that $G_{\ell}$ has a tensor product structure with multi-index set $\mathbf{j} = \left( j^{(1)}, j^{(2)}, j^{(3)}\right) \in \mathcal{J}_{G_{\ell}} = \bigtimes^3_{d = 1} \left\{ 1, \ldots, m^{(d)}_{\ell} \right\}$. \\
To introduce adaptivity we consider a sequence of nested domains  
\begin{equation}
	\label{eq:nested_domains}
	\left[ 0, 1 \right]^3 = \Omega_0 \supseteq \Omega_1 \supseteq \ldots \supseteq \Omega_{L-1} \supseteq \Omega_{L} = \emptyset
\end{equation}
where each $\Omega_\ell \subset \left[ 0, 1 \right]^3$ represents the region selected to be refined at level $\ell$ and its boundary $\partial \Omega_\ell$ is aligned with the breaking points in each dimension, i.e.~it is a union of cells from $G_{\ell}$. With that we are able to define a the hierarchical mesh as
\begin{equation}
	\label{eq:hierarchical_mesh}
	\mathcal{G} = \bigcup^{L-1}_{l = 0} \left\{ q \in G_{\ell} \colon \, q \subset \Omega_\ell \wedge q \not\subset \Omega_{\ell+1} \right\}.
\end{equation}
A cell $q$ is called active if $q \in \mathcal{G}$ and it is called an active cell of level $\ell$ if $q \in G_{\ell} \cap \mathcal{G}$. In the following we call a cell $q \in G_{\ell} \setminus \mathcal{G}$ a non-active cell of level $\ell$ and we point out that this should not be confused with deactivated cells of level $\ell$ from~\cite{GarauVazquez:2018}, which are cells $q \in G_{\ell}$ with $q \subset \Omega_{\ell+1}$. Non-active cells of level $\ell$ are either deactivated ($q \subset \Omega_{\ell+1}$) or non-selected ($q \not\subset \Omega_\ell$), as visualized for a two-dimensional example in \cref{fig:non_active_cells}. Only the former coincide with the deactivated cells of~\cite{GarauVazquez:2018}.

We define the hierarchical B-spline basis $\mathcal{H}$ by the following algorithm:
\begin{enumerate}
	\item $\mathcal{H}_0 = \mathcal{B}_0$
	\item for $\ell = 0, \ldots, L-2$
	\begin{equation*}
		\mathcal{H}_{\ell+1} = \left\{ \beta \in \mathcal{H}_\ell : \, \operatorname{supp} \left( \beta \right) \not \subset \Omega_{\ell+1} \right\} \cup \left\{ \beta \in \mathcal{B}_{\ell+1} : \, \operatorname{supp} \left( \beta \right) \subset \Omega_{\ell+1} \right\}
	\end{equation*}
	\item $\mathcal{H} = \mathcal{H}_{L-1}$.
\end{enumerate}
Hierarchical B-splines, while flexible, lack the partition of unity property and can exhibit significant overlap between basis functions from different levels. This overlap not only complicates the function space but also leads to denser system matrices. To address these issues, an alternative basis known as Truncated Hierarchical B-splines (THB-splines) was introduced in~\cite{GiannelliJuettlerSpeleers:2012}. \\
THB-splines preserve the desirable properties of hierarchical splines while restoring the partition of unity and reducing the extent of overlapping supports. These improvements lead to sparser matrices and better numerical behavior, making THB-splines a more suitable choice for adaptive refinement strategies and approximation theory. For that the so-called two-scale relation with non-negative coefficients is utilized, which states that B-splines of level $\ell$ can be written as linear combinations of B-splines of level $\ell + 1$ with non-negative coefficients, i.e. for $\beta^{\ell}_{\mathbf{i}^{\ell}} \in \mathcal{B}_{\ell}$
\begin{equation}
	\label{eq:two_scale_relation}
	\begin{aligned}
		\beta_{\ell, \mathbf{j}} & = \beta^{(1)}_{\ell, i^{(1)}} \beta^{(2)}_{\ell, i^{(2)}} \beta^{(3)}_{\ell, i^{(3)}} \\
		& = \sum_{\mathbf{i} \in \mathcal{I}_{\mathcal{B}_{\ell+1}}} \left( \mathbf{C}^{\left( \ell + 1, \ell \right)} \right)_{\mathbf{i}, \mathbf{j}} \beta_{\ell+1, \mathbf{i}} \\
		& = \left( {\mathbf{C}^{\left( \ell + 1, \ell \right)}}^{\top} \cdot \mathbf{B}_{\ell+1} \right)_{\mathbf{j}}
	\end{aligned}
\end{equation}
where $\cdot$ denotes the contracted product of the respective dimensions and the coefficient tensor factors as the Kronecker product
\begin{equation}
	\label{eq:two_scale_relation_coefficients}
	\mathbf{C}^{\left( \ell + 1, \ell \right)} = C^{\left( \ell + 1, \ell \right), (1)} \otimes C^{\left( \ell + 1, \ell \right), (2)} \otimes C^{\left( \ell + 1, \ell \right), (3)} 
\end{equation}
of size $\left( n^{(1)}_{\ell+1}, n^{(2)}_{\ell+1}, n^{(3)}_{\ell+1} \right) \times \left( n^{(1)}_{\ell}, n^{(2)}_{\ell}, n^{(3)}_{\ell} \right)$. The univariate coefficient matrices $C^{\left( \ell + 1, \ell \right), (d)} \in \mathbb{R}^{n^{(d)}_{\ell+1} \times n^{(d)}_{\ell}}$ with $C^{\left( \ell + 1, \ell \right), (d)}_{i, j} \geq 0$ can be computed using the knot refinement algorithm (see~\cite[5.3]{piegl}). The two-scale relation \cref{eq:two_scale_relation} holds across multiple levels, which means that we can represent a basis function from level $\ell$ as a linear combination of basis functions from level $\ell + m$, i.e.
\begin{align*}
	\beta_{\ell, \mathbf{j}} & = \sum_{\mathbf{i} \in \mathcal{I}_{\mathcal{B}_{\ell+m}}} \left( \mathbf{C}^{\left( \ell + m, \ell \right)} \right)_{\mathbf{i}, \mathbf{j}} \beta_{\ell+m, \mathbf{i}} \\
	& = \left( {\mathbf{C}^{\left( \ell + m, \ell \right)}}^{\top} \cdot \mathbf{B}_{\ell+m} \right)_{\mathbf{j}}
\end{align*}
where the coefficient tensor factors as a Kronecker product
\begin{equation*}
	\mathbf{C}^{\left( \ell + m, \ell \right)} = C^{\left( \ell + m, \ell \right), (1)} \otimes C^{\left( \ell + m, \ell \right), (2)} \otimes C^{\left( \ell + m, \ell \right), (3)},
\end{equation*}
which is of size $\left( n^{(1)}_{\ell+m}, n^{(2)}_{\ell+m}, n^{(3)}_{\ell+m} \right) \times \left( n^{(1)}_{\ell}, n^{(2)}_{\ell}, n^{(3)}_{\ell} \right)$ and the unitary coefficient matrices are just the products of the unitary matrices of the levels in succession, i.e. 
\begin{gather*}
	C^{\left( \ell + m, \ell \right), (d)} = C^{\left( \ell+m, \ell+m-1\right), (d)} C^{\left(\ell+m-1, \ell+m-2\right), (d)} \cdots C^{\left( \ell + 1, \ell \right), (d)} \in \mathbb{R}^{n^{(d)}_{\ell+m} \times n^{(d)}_{\ell}}.
\end{gather*}
The construction of THB-splines is based on a truncation operator, which consists of removing in the two-scale relation \cref{eq:two_scale_relation} the contribution of active and deactivated functions of level $\ell + 1$, that are functions in $\mathcal{H}_{\ell+1} \cap \mathcal{B}_{\ell+1}$. Let \cref{eq:two_scale_relation} be the representation of $\beta_{\ell, \mathbf{i}} \in \mathcal{B}_{\ell}$ with respect to the finer basis of $\mathcal{B}_{\ell+1}$, then the truncation of $\beta_{\ell, \mathbf{i}}$ with respect to $\mathcal{B}_{\ell + 1}$ and $\Omega_{\ell+1}$ is defined as
\begin{align*}
	\operatorname{trunc}^{\ell + 1} \left( \beta_{\ell, \mathbf{j}} \right) & = \sum_{\substack{\mathbf{i} \in \mathcal{I}_{\mathcal{B}_{\ell+1}} \\ \operatorname{supp} \left(\beta_{\ell+1, \mathbf{i}}\right) \not\subseteq \Omega_{\ell + 1}}} \left(\mathbf{C}^{\left( \ell + 1, \ell \right)} \right)_{\mathbf{i}, \mathbf{j}} \beta_{\ell+1, \mathbf{i}} \\
	& = \sum_{\mathbf{i}\in \mathcal{I}_{\mathcal{B}_{\ell+1}}} \left(\mathbf{C}_{t}^{\left( \ell + 1, \ell \right)}\right)_{\mathbf{i}, \mathbf{j}} \beta_{\ell+1, \mathbf{i}} \\
	& = \left( {\mathbf{C}^{\left( \ell + 1, \ell \right)}_{t}}^{\top} \cdot \mathbf{B}_{\ell+1} \right)_{\mathbf{j}}
\end{align*}
with coefficients
\begin{equation}
	\label{eq:truncation_coefficients}
	\left(\mathbf{C}_{t}^{\left( \ell + 1, \ell \right)} \right)_{\mathbf{i}, \mathbf{j}} =
	\begin{cases}
		0 &: \text{ if } \beta^{\ell+1}_{\mathbf{i}} \in \mathcal{H}_{\ell+1} \cap \mathcal{B}_{\ell+1}, \\
		\left( \mathbf{C}^{\left( \ell + 1, \ell \right)} \right)_{\mathbf{i}, \mathbf{j}}  &: \text{ otherwise}.
	\end{cases}
\end{equation}
We point out that the coefficient tensor $\mathbf{C}_{t}^{\left( \ell + 1, \ell \right)} \in \mathbb{R}^{\left( n^{(1)}_{\ell+1}, n^{(2)}_{\ell+1}, n^{(3)}_{\ell+1} \right) \times \left( n^{(1)}_{\ell}, n^{(2)}_{\ell}, n^{(3)}_{\ell} \right)} $ has in general not a rank-one tensor product structure as in \cref{eq:two_scale_relation_coefficients}, which means it can be only represented as a sum of several rank-one tensors. Since we will only be working with THB-splines and not just hierarchical B-splines in the following, we will simplify the notation and use $\mathbf{C}$ to denote the coefficient tensor of the two-scale relation with truncation \cref{eq:truncation_coefficients}.

By applying the truncation operator to hierarchical B-splines of coarse levels, we get the THB-spline basis by the following recursive algorithm
\begin{enumerate}
	\item $\mathcal{T}_0 = \mathcal{B}_0$
	\item for $\ell = 0, \ldots, L-2$
	\begin{align*}
		\mathcal{T}_{\ell+1} = & \left\{ \operatorname{trunc}^{\ell+1} \left(\beta \right) : \; \beta \in \mathcal{T}_\ell \land \operatorname{supp} \left(\beta\right) \not\subset \Omega_{\ell+1} \right\} \\
		& \cup \left\{ \beta \in \mathcal{B}_{\ell+1} : \; \operatorname{supp} \left( \beta \right) \subset \Omega_{\ell+1} \right\}
	\end{align*}
	\item $\mathcal{T} = \mathcal{T}_{L-1}$.
\end{enumerate}
A basis function of $\mathcal{B}_{\ell}$ is after the truncation not a basis function of $\mathcal{B}_{\ell}$ anymore but represented by basis functions of $\mathcal{B}_{\ell'}$, where $\ell' > \ell$ is the level of the finest function that truncates it. In the following we call a basis function $\beta$ active if $\beta \in \mathcal{T}$, we call it an active basis function of level $\ell$ if $\beta \in \mathcal{B}_{\ell} \cap \mathcal{T}$ and we call it deactivated basis function of level $\ell$ if $\beta \in \mathcal{T}_{\ell} \setminus \mathcal{T}_{\ell+1}$, which means $\beta \in \mathcal{B}_{\ell}$ and $\operatorname{supp} \left( \beta \right) \subset \Omega_{\ell+1}$. We denote by $\mathcal{T}_{\ell} \cap \mathcal{B}_{\ell}$ the set of active and deactivated functions of level $\ell$ and by $\mathcal{B}_{\ell} \setminus \mathcal{T}$ the set of non-active basis functions of level $\ell$.

\section{Assembling a low-rank linear system with THB-splines}
\label{section:assembling_linear_system}
In this section, we present our approach to assembling low-rank system matrices under appropriate conditions using tensors in TT format. This approach consists of four main steps: first, interpolation of the weight function and the source function using multivariate B-splines to obtain a separable approximation; second, assembling a low-rank mass or stiffness matrix for each level by integrating all basis functions over all active cells; third, realizing the two-scale relation by performing linear combinations of the system matrices; and last, assembling the hierarchical mass or stiffness matrix. For the last three steps, we follow the assembling heuristic described in~\cite[Section~3.2]{GarauVazquez:2018}, to which we have added our low-rank strategies developed for this setup. Our goal is to obtain a linear system of the form
\begin{equation}
	\label{eq:linear_system}
	\mathbf{K} \mathbf{y} = \mathbf{f},
\end{equation}
to discretize the elliptic problem \cref{eq:forward}, which keeps memory requirements low and enables efficient computation of approximate solutions using methods from the \textsc{TT-Toolbox}~\cite{tt-toolbox}.

\subsection{Low-Rank interpolation}
\label{subsection:Low-Rank_interpolation}
The weight functions $\omega$ and $Q$ as well as the source function $f$ are generally non-separable, which unfortunately results in expensive multivariate quadrature having to be applied when assembling the linear system. To overcome this issue, we interpolate these with separable multivariate B-splines and represent the corresponding coefficient tensor in TT format with low rank. \\

We assume that the computational domain $\Omega \subset \mathbb{R}^3$ is a B-spline or NURBS geometry, which means it can be represented using a geometry map 
\begin{gather*}
	F \colon \left[0, 1\right]^3 = \hat{\Omega} \rightarrow \Omega \\
	F \left( \hat{x} \right)  = \sum_{\mathbf{i} \in \mathcal{I}} \mathbf{A}_\mathbf{i} \beta_{\mathbf{i}} \left( \hat{x} \right)
\end{gather*}
where $\beta_{\mathbf{i}}$ are multivariate B-splines or NURBS, $\mathbf{A}_{\mathbf{i}} \in \R^3$ are so-called control points. These are chosen such that $F \left(\hat{x}\right)$ is at least injective and $\mathcal{I} = \bigtimes^3_{d = 1} \left\{ 1, \ldots, n^{(d)} \right\}$ is a multi-index set. Let $V_h$ be a B-Spline or NURBS space which should now be used for the Galerkin discretization, resulting in the discrete mass and stiffness terms
\begin{align*}
	a_{m} \left(u_h, v_h \right) &= \int_\Omega u_h(x) v_h(x) \D x = \int_{\hat{\Omega}} \sum_{\mathbf{i} \in \mathcal{I}_0} \mathbf{u}_\mathbf{i} \beta_\mathbf{i} \left( \hat{x} \right)  \sum_{\mathbf{j} \in \mathcal{I}_0} \mathbf{v}_\mathbf{j} \beta_\mathbf{j} \left( \hat{x} \right)  \omega \left( \hat{x} \right)  \D \hat{x}, \\
	a_{s} \left(u_h, v_h \right) &= \int_{\Omega} \nabla u_h(x) \cdot \nabla v_h(x) \D x = \int_{\hat{\Omega}} \left(Q \left( \hat{x} \right)  \sum_{\mathbf{i} \in \mathcal{I}_0} \mathbf{u}_\mathbf{i} \nabla \beta_\mathbf{i} \left( \hat{x} \right) \right) \cdot \sum_{\mathbf{j} \in \mathcal{I}_0} \mathbf{v}_\mathbf{j} \nabla \beta_\mathbf{j} \left( \hat{x} \right)  \D \hat{x},
\end{align*}
for $u_h, v_h \in V_h$ and $\mathcal{I}_0$ a multi-index set that excludes basis functions on the boundary to satisfy the homogenouos Dirichlet conditions, with the additional terms stemming from the coordinate transformation,
\begin{equation}
	\label{eq:weight_functions}
	\omega \left( \hat{x} \right)  = \left| \det \left( \nabla F \left( \hat{x} \right) \right) \right| \in \R, \quad Q \left( \hat{x} \right) = \omega \left( \hat{x} \right) \left( \nabla F \left( \hat{x} \right)^{T} \nabla F \left( \hat{x} \right) \right)^{-1} \in \R^{3\times 3},
\end{equation}
which we will refer to as weight functions from now on, as introduced in~\cite{angelos1}. We point out that the spline space for representing the computational domain by the geometry map and the spline space for the Galerkin discretization do not have to be the same and that from now on we assume for the rest of the work that the spline space for representing the computational domain has no hierarchical structure. \\
The weight functions $\omega \left( \hat{x} \right)$ and $Q \left( \hat{x} \right)$, which are determined by the geometry map $F \left( \hat{x} \right)$, are in general not separable into one-dimensional factors. To overcome this issue we are using interpolation
\begin{equation}
	\label{eq:weightfunction_linear_system}
	\omega \left( \hat{x} \right) \approx \mathbf{W} \colon \mathbf{\hat{B}} \left(\hat{x}\right).
\end{equation}
where $\colon$ denotes the Frobenius product and 
\begin{gather*}
	\mathbf{\hat{B}} \left(\hat{x}\right) = \bigotimes^3_{d = 1} \hat{B}^{(d)} \left( \hat{x}^{(d)} \right) \in \mathbb{R}^{\left( \hat{n}^{(1)}, \hat{n}^{(2)}, \hat{n}^{(3)} \right)} \\
	\hat{B}^{(d)} \left( \hat{x}^{(d)} \right) = \left[ \hat{\beta}^{(d)}_1 \left( \hat{x}^{(d)} \right), \ldots, \hat{\beta}^{(d)}_{\hat{n}^{(d)}} \left( \hat{x}^{(d)} \right) \right]^\top \in \mathbb{R}^{\hat{n}^{(d)}}, \quad \hat{x}^{(d)} \in \left[0, 1\right],
\end{gather*}
are multivariate B-splines of higher order and less continuity (i.e.~$\hat{n}^{(d)} > n^{(d)}$, c.f.~\cite{angelos1}) which are defined as the tensor product of univariate B-splines $\hat{\beta}^{(d)}_{i} \colon \left[0, 1\right] \to \left[0, 1\right]$, $i = 1, \ldots, \hat{n}^{(d)}$, and are therefore separable. For determining this coefficient tensor $\mathbf{W} \in \mathbb{R}^{\left( \hat{n}^{(1)}, \hat{n}^{(2)}, \hat{n}^{(3)} \right)}$ we follow the approach introduced in~\cite{BuengerDolgovStoll:2020}, which means that we evaluate the weight function $\omega$ at the so-called Greville points $\mathbf{\hat{X}} = \hat{X}^{(1)} \otimes \hat{X}^{(3)} \otimes \hat{X}^{(3)} \in \mathbb{R}^{\left( \hat{n}^{(1)}, \hat{n}^{(2)}, \hat{n}^{(3)} \right)}$~\cite{greville}, which are separable. We solve the following linear system 
\begin{equation}
	\label{eq:interpolation_system}
	\mathbf{W} \colon \left( \bigotimes^3_{d = 1} \hat{B}^{(d)} \left( \hat{X}^{(d)} \right) \right) = \omega \left( \mathbf{\hat{X}} \right)
\end{equation}
using the alternating minimal energy, \amen, solver introduced in~\cite{amen}, which can work directly with the factor matrices $\hat{B}^{(d)} \left( \hat{X}^{(d)} \right) \in \mathbb{R}^{\hat{n}^{(d)} \times \hat{n}^{(d)}}$. The resulting solution tensor $\mathbf{W}_{R} \in \mathbb{R}^{\left( \hat{n}^{(1)}, \hat{n}^{(2)}, \hat{n}^{(3)} \right)}$ is then in tensor train format, which can also be represented as CP format via 
\begin{align*}
	\mathbf{W}_{R} & = \sum_{r^{(1)}=1}^{R^{(1)}} \sum_{r^{(2)}=1}^{R^{(2)}} \mathbf{W}^{\left(1\right)}_{R} \left( 1,:,r^{(1)} \right) \otimes \mathbf{W}^{\left(2\right)}_{R} \left( r^{(1)},:,r^{(2)} \right) \otimes \mathbf{W}^{\left(3\right)}_{R} \left( r^{(2)},:,1 \right) \\
	& = \sum_{r = 1}^{R} \bigotimes^{3}_{d=1}  w_r^{(d)} \approx \mathbf{W},
\end{align*}
where $\mathbf{W}^{\left( d \right)}_{R} \in \mathbb{R}^{R^{(d-1)}, \hat{n}^{(d)}, R^{(d)}}$ are the TT cores and $w_r^{(d)} \in \mathbb{R}^{\hat{n}^{(d)}}$, $r = 1, \ldots, R$, $R = R^{(1)} R^{(2)}$, represent the TT cores in the CP format. The TT ranks $R^{(d-1)}, R^{(d)} \in \mathbb{N}$ depend on the complexity of the weight function $\omega$ on the one hand and on the accuracy with which the linear system \cref{eq:interpolation_system} is solved on the other (with $R^{(0)} = R^{(3)} = 1$). With this we get a separable low-rank representation of the weight function,
\begin{equation*}
	\omega \left( \hat{x} \right) \approx \mathbf{W}_{R} \colon \mathbf{\hat{B}} \left(\hat{x}\right) = \sum_{r=1}^{R} \prod_{d=1}^3 w_r^{(d)}\cdot \hat{B}^{(d)} \left( \hat{x}^{(d)} \right),
\end{equation*}
where $\hat{B}^{(d)} \left( \hat{x}^{(d)} \right) \in \mathbb{R}^{\hat{n}^{(d)}}$ denotes the vector holding all univariate basis functions evaluated in $\hat{x}^{(d)} \in \left[0, 1\right]$. If a B-spline space $V_h$ without hierarchical structure is used for the Galerkin discretization, with basis functions for each parametric dimension summarized in a vector $B^{(d)} \left( \hat{x}^{(d)} \right)$, we can write the mass tensor as a sum of tensor products of small univariate mass matrices
\begin{equation} 
	\label{eq:massFinal}
	\begin{aligned}
		\mathbf{M} &\approx \sum_{r=1}^R \bigotimes_{d=1}^3 \int_{0}^{1} \left( w_r^{(d)}\cdot \hat{B}^{(d)} \left( \hat{x}^{(d)} \right) \right) \, B^{(d)} \left( \hat{x}^{(d)} \right) \otimes B^{(d)} \left( \hat{x}^{(d)} \right) \D \hat{x}^{(d)} = \\
		&= \sum_{r=1}^R \bigotimes_{d=1}^3 M_r^{\left( d \right)}.
	\end{aligned}
\end{equation}
The same procedure can be applied to each entry of $Q \left( \hat{x} \right) \in \mathbb{R}^{3 \times 3}$ such that we get a low-rank tensor representation of the stiffness tensor by
\begin{equation} 
	\label{eq:stiffnessFinal}
	\mathbf{K} \approx \sum_{k,l=1}^3 \sum_{r=1}^R \bigotimes_{d=1}^3 K_{k,l,r}^{\left( d \right)}.
\end{equation}
Even if $\mathbf{M}$ and $\mathbf{K}$ are represented here in CP format, we will assume in the following that they are in TT format. This is not a restriction, as it is possible to convert between the two formats without much effort. We refer to~\cite{BuengerDolgovStoll:2020} for details and to the codes on our website~\cite{BuengerCode}. We point out that we can influence the accuracy of the interpolation by $h$- and $p$-refinement of the spline space used to represent the computational domain via the geometry map. \\
To obtain a separable representation for the source function $f\colon \mathbb{R}^3 \to \mathbb{R}$ we follow the same methology and interpolate $f$ with multivariate B-splines, i.e.   
\begin{equation*}
	f \left( F  \left( \hat{x} \right) \right) \approx \mathbf{F} \colon \mathbf{\tilde{B}} \left(\hat{x}\right),
\end{equation*}
to get a low-rank representation $\mathbf{F}_{\tilde{R}} \in \mathbb{R}^{\left( \tilde{n}_1, \tilde{n}_2, \tilde{n}_3 \right)}$ of the coefficient tensor $\mathbf{F}$. If $h$- and $p$-refinement is done to the spline space for respresenting the geometry, we use the resulting geometry map $F$, otherwise the initial one. In order to interpolate $f$, it makes sense to choose a spline space that provides the best possible interpolation. The respective spline space should be fine enough and have the appropriate degrees so that the source term $f$, which may have a lot of activity in certain parts of the domain, is interpolated well, but at the same time the number of splines should not be too large, since $f$ is evaluated at the Greville points $\mathbf{\tilde{X}} = \tilde{X}^{(1)} \otimes \tilde{X}^{(3)} \otimes \tilde{X}^{(3)} \in \mathbb{R}^{\left( \tilde{n}_1, \tilde{n}_2, \tilde{n}_3 \right)}$ whose number is as much as the number of splines for the interpolation.

\subsection{Domain-wise Integration}
\label{subsection:domain_wise_integration}
The next step is to assemble a mass or stiffness matrix for each level by integrating all basis functions of this level $\mathcal{B}_{\ell}$ over the active cells of this level $G_{\ell} \cap \mathcal{G}$, as described in~\cite[Section~3.2]{GarauVazquez:2018}. In the following we consider this procedure for assembling a mass matrix, but all operations are transferable to the stiffness matrix and the right-hand side. \\

\noindent The entries of this mass matrix $\mathbf{M}_\ell \in \mathbb{R}^{ \left( n^{(1)}_{\ell}, n^{(2)}_{\ell}, n^{(3)}_{\ell} \right) \times \left( n^{(1)}_{\ell}, n^{(2)}_{\ell}, n^{(3)}_{\ell} \right) }$, $\ell = 0, \ldots, L-1$, are defined by
\begin{equation}
	\label{eq:level_mass}
	\begin{aligned}
		\left( \mathbf{M}_{\ell} \right)_{\mathbf{i}, \mathbf{j}} & = \sum_{q \in G_{\ell} \cap \mathcal{G}} \int_{q} \beta_{\ell, \mathbf{i}} \, \beta_{\ell, \mathbf{j}} \, \omega \, \mathrm{d} x = \sum_{\substack{q \in G_{\ell} \cap \mathcal{G} \\ q = q^{(1)} \times q^{(2)} \times q^{(3)}}} \int_{q^{(1)}} \int_{q^{(2)}} \int_{q^{(3)}} \beta_{\ell, \mathbf{i}} \, \beta_{\ell, \mathbf{j}} \, \omega \, \mathrm{d} x \\
		& = \sum_{\substack{q \in G_{\ell} \cap \mathcal{G} \\ q = q^{(1)} \times q^{(2)} \times q^{(3)}}} \int^{\hat{\xi}^{(1)}_{j_q + 1}}_{\hat{\xi}^{(1)}_{j_q}} \int^{\hat{\xi}^{(2)}_{j_q + 1}}_{\hat{\xi}^{(2)}_{j_q}} \int^{\hat{\xi}^{(3)}_{j_q + 1}}_{\hat{\xi}^{(3)}_{j_q}} \beta_{\ell, \mathbf{i}} \, \beta_{\ell, \mathbf{j}} \, \omega \, \mathrm{d} x,
	\end{aligned}
\end{equation}
where $\beta_{\ell, \mathbf{i}}, \beta_{\ell, \mathbf{j}} \in \mathcal{B}_{\ell}$ and $q^{(d)} = \left[ \hat{\xi}^{(d)}_{j_q}, \hat{\xi}^{(d)}_{j_q + 1} \right]$ denotes the size of the active cell $q \in G_{\ell} \cap \mathcal{G}$ in dimension $d$, which is just an interval bounded by two consecutive breakpoints $\hat{\xi}^{(d)}_{j_q}, \hat{\xi}^{(d)}_{j_q + 1} \in \hat{\xi}^{(d)}$ resulting from the refinement of this level $\ell$ in this dimension $d$.\\
Combining the cell-wise integration \cref{eq:level_mass} with the low-rank formulation \cref{eq:massFinal} results in
\begin{equation}
	\label{eq:level_mass_one_cell}
	\begin{gathered}
		\mathbf{M}_\ell \approx \sum_{\substack{q \in G_{\ell} \cap \mathcal{G} \\ q = q^{(1)} \times q^{(2)} \times q^{(3)}}} \sum_{r=1}^R \bigotimes_{d=1}^3 M_{\ell, r, q^{\left( d \right)}}^{\left( d \right)}, \\
		\begin{aligned}
			M_{\ell, r, q^{\left( d \right)}}^{\left( d \right)} & = \int_{q^{(d)}} \left( w_r^{(d)}\cdot \hat{B}^{(d)} \left( \hat{x}^{(d)} \right) \right) \, B^{(d)}_{\ell} \left( \hat{x}^{(d)} \right) \otimes B^{(d)}_{\ell} \left( \hat{x}^{(d)} \right) \D \hat{x}^{(d)} \\
			& = \int^{\hat{\xi}^{(d)}_{j_q + 1}}_{\hat{\xi}^{(d)}_{j_q}} \left( w_r^{(d)}\cdot \hat{B}^{(d)} \left( \hat{x}^{(d)} \right) \right) \, B^{(d)}_{\ell} \left( \hat{x}^{(d)} \right) \otimes B^{(d)}_{\ell} \left( \hat{x}^{(d)} \right) \D \hat{x}^{(d)},
		\end{aligned}
	\end{gathered}
\end{equation}
where $B^{(d)}_{\ell} \left( \hat{x}^{(d)} \right) = \left[ \beta^{(d)}_{\ell, 1} \left( \hat{x}^{(d)} \right), \ldots, \beta^{(d)}_{\ell, n^{(d)}_{\ell}} \left( \hat{x}^{(d)} \right) \right] \in \mathbb{R}^{n^{(d)}_{\ell}}$ and $n^{(d)}_{\ell}$ is the number of B-Splines in the dimension $d = 1, 2, 3$ of $\mathcal{B}_{\ell}$. We point out that the B-splines $\hat{B}^{(d)}$ for interpolating the weight function $\omega$ are not dependent on the level $\ell$, just like the B-splines for interpolating the source function $f$ for the right hand side.

In \cref{eq:level_mass_one_cell}, two problems arise. First, the number of splines, $n^{(1)}_{\ell} n^{(2)}_{\ell} n^{(3)}_{\ell}$, grows exponentially with the level $\ell$. Consequently, even tensors in  CP or TT format can no longer be stored once they exceed a certain size, despite the fact that many entries of $\mathbf{M}_{\ell}$ are zero because no integration is required for splines whose support lies entirely on non-active cells of this level. Second, the number of active cells, $\lvert G_{\ell} \cap \mathcal{G} \rvert$, may also increase rapidly. In \cref{eq:level_mass_one_cell} we assemble a low-rank matrix or tensor for every active cell $q \in G_{\ell} \cap \mathcal{G}$ and subsequently sum over all active cells; this procedure is computationally not efficient.\\

To overcome the first issue, i.e.~the high number of splines $n^{(1)}_{\ell} n^{(2)}_{\ell} n^{(3)}_{\ell}$, we reduce the basis $\mathcal{B}_{\ell}$ to a subset that still has a tensor product structure and contains all splines relevant for integration. We utilise the tensor product structure of the spline basis $\mathcal{B}_{\ell} = \mathcal{B}^{(1)}_{\ell} \otimes \mathcal{B}^{(2)}_{\ell} \otimes \mathcal{B}^{(3)}_{\ell}$ and the index set $\mathcal{I}_{\mathcal{B}_{\ell}} = \bigtimes^3_{d = 1} \left\{ 1, \ldots, n^{(d)}_{\ell} \right\}$. To compute \cref{eq:level_mass} we only need to consider the splines from $\mathcal{B}_{\ell}$ whose support intersects with or lies in the set of active cells $G_{\ell} \cap \mathcal{G}$. Since we need a span with a tensor product structure, we cannot simply discard the remaining splines from $\mathcal{B}_{\ell}$ and work with the resulting subset. But we can discard slices of splines from $\mathcal{B}_{\ell}$ whose support does not intersect with $G_{\ell} \cap \mathcal{G}$. By slices, we mean that we fix one component of the multi-index $\left( i^{(1)}_{\ell}, i^{(2)}_{\ell}, i^{(3)}_{\ell} \right)$ and examine the resulting slice, for example $\left( :, i^{(2)}_{\ell}, : \right) = \left\{ 1, \ldots, n^{(1)}_{\ell} \right\} \times \left\{ i^{(2)}_{\ell} \right\} \times \left\{ 1, \ldots, n^{(3)}_{\ell} \right\}$. If all splines which are represented by e.g.~the slices $\left( :, i^{(2)}_{\ell}, : \right)$ do not have a support which intersects with $G_{\ell} \cap \mathcal{G}$, we remove these from the basis $\mathcal{B}_{\ell}$, i.e.~$\mathcal{B}_{\ell} \mapsto \mathcal{B}_{\ell} \setminus \left\{ \beta^{(1)}_{\ell, 1}, \ldots, \beta^{(1)}_{\ell, n^{(1)}_{\ell}} \right\} \otimes \left\{ \beta^{(2)}_{\ell, i^{(2)}_{\ell}} \right\} \otimes \left\{ \beta^{(3)}_{\ell, 1}, \ldots, \beta^{(3)}_{\ell, n^{(3)}_{\ell}} \right\}$. This yields a new tensor product basis $\tilde{\mathcal{B}}_{\ell} = \tilde{\mathcal{B}}^{(1)}_{\ell} \otimes \tilde{\mathcal{B}}^{(2)}_{\ell} \otimes \tilde{\mathcal{B}}^{(3)}_{\ell}$, where $\tilde{\mathcal{B}}^{(d)}_{\ell} = \left\{ \tilde{\beta}^{(d)}_{\ell, i} \colon \left[0, 1\right] \to \mathbb{R} \, \colon \, i = 1, \ldots, \tilde{n}^{(d)}_{\ell} \right\}$, with a new numbering $\mathcal{I}_{\tilde{\mathcal{B}}_{\ell}} = \bigtimes^3_{d = 1} \left\{ 1, \ldots, \tilde{n}^{(d)}_{\ell} \right\}$, which contains all splines used in the integration process. By removing the slices of splines, the information of the tensor product structure of $\mathcal{B}_{\ell}$ is preserved and the splines from $\tilde{\mathcal{B}}_{\ell}$ can be clearly assigned to the splines from $\mathcal{B}_{\ell}$. We point out that there may still be splines in the basis $\tilde{\mathcal{B}}_{\ell}$ that do not have support which intersects $G_{\ell} \cap \mathcal{G}$, but these must be included so that $\tilde{\mathcal{B}}_{\ell}$ can have a tensor product structure. It obviously holds $\tilde{\mathcal{B}}_{\ell} \subset \mathcal{B}_{\ell}$, but $\tilde{n}^{(d)}_{\ell} < n^{(d)}_{\ell}$ does not necessarily have to apply, since the set of active cells $G_{\ell} \cap \mathcal{G}$ can be embedded in $G_{\ell}$ in such a way that all indices $i^{(d)}_{\ell} = 1, \ldots, n^{(d)}_{\ell}$ occur in the $d$-th component of the multi-index $\left( i^{(1)}_{\ell}, i^{(2)}_{\ell}, i^{(3)}_{\ell} \right)$ in the set of splines that have support on $G_{\ell} \cap \mathcal{G}$, then no slice of splines with index in this dimension $d$ can be discarded and $\mathcal{B}_{\ell}$ does not become smaller in this dimension. \\

\tikzset{baseline on center/.style={baseline={(current bounding box.center)}}}
\begin{figure}
	\begin{center}
		\begin{tikzpicture}[scale=0.45, baseline on center]
			
			\foreach \x/\y in {%
				1/0,1.5/0,1/0.5,1.5/0.5,%
				2/0,2.5/0,2/0.5,2.5/0.5,%
				4/0,4.5/0,4/0.5,4.5/0.5,%
				5/0,5.5/0,5/0.5,5.5/0.5,%
				1/1,1.5/1,1/1.5,1.5/1.5,%
				5/1,5.5/1,5/1.5,5.5/1.5,%
				1/3,1.5/3,1/3.5,1.5/3.5,%
				4/3,4.5/3,4/3.5,4.5/3.5,%
				5/3,5.5/3,5/3.5,5.5/3.5,%
				1/4,1.5/4,1/4.5,1.5/4.5,%
				2/4,2.5/4,2/4.5,2.5/4.5,%
				4/4,4.5/4,4/4.5,4.5/4.5,%
				5/4,5.5/4,5/4.5,5.5/4.5}
			{%
				\fill[gray!50] (\x,\y) rectangle ++(0.5,0.5);
			}
			
			\foreach \x in {0,...,12}   \draw (\x/2,0) -- (\x/2,6);
			\foreach \y in {0,...,12}   \draw (0,\y/2) -- (6,\y/2);
			
			\node at (3,-0.7) {$G_{\ell}$};
		\end{tikzpicture}
		$\;\xrightarrow{\substack{\text{Discarding slices of}\\\text{non-active cells and}\\\text{and introducing}\\\text{new numbering}}}\;$
		\begin{tikzpicture}[scale=0.45, baseline on center]
			
			\foreach \x/\y in {%
				1/1,2/1,3/1,4/1,%
				1/2,4/2,%
				1/3,3/3,4/3,%
				1/4,2/4,3/4,4/4}
			{%
				\foreach \dx in {0,0.5}
				\foreach \dy in {0,0.5}
				\fill[gray!50] (\x-1+\dx,\y-1+\dy) rectangle ++(0.5,0.5);
			}
			
			\foreach \x in {0,...,8} \draw (\x/2,0) -- (\x/2,4);
			\foreach \y in {0,...,8} \draw (0,\y/2) -- (4,\y/2);
			
			\node at (2,-0.7) {$\tilde{G}_{\ell}$};
		\end{tikzpicture}
		$\;\xrightarrow{\substack{\text{\textbf{Greedy method}:}\\\text{Partitioning}\\\text{into cuboids}}}\;$
		\begin{tikzpicture}[scale=0.45,pattern color=black, baseline on center]
			
			\foreach \x/\y in {1/1,2/1,3/1,4/1}
			\foreach \dx in {0,0.5}
			\foreach \dy in {0,0.5}
			\fill[blue!40] (\x-1+\dx,\y-1+\dy) rectangle ++(0.5,0.5);
			
			\foreach \x/\y in {1/2,1/3,1/4}
			\foreach \dx in {0,0.5}
			\foreach \dy in {0,0.5}
			\fill[red!60] (\x-1+\dx,\y-1+\dy) rectangle ++(0.5,0.5);
			
			\foreach \x/\y in {4/2,4/3,4/4}
			\foreach \dx in {0,0.5}
			\foreach \dy in {0,0.5}
			\fill[green!50] (\x-1+\dx,\y-1+\dy) rectangle ++(0.5,0.5);
			
			\foreach \x/\y in {3/3,3/4}
			\foreach \dx in {0,0.5}
			\foreach \dy in {0,0.5}
			\fill[purple!50] (\x-1+\dx,\y-1+\dy) rectangle ++(0.5,0.5);
			
			\foreach \dx in {0,0.5}
			\foreach \dy in {0,0.5}
			\fill[yellow!70] (1+\dx,3+\dy) rectangle ++(0.5,0.5);
			
			\foreach \x/\y in {2/2,3/2}
			\foreach \dx in {0,0.5}
			\foreach \dy in {0,0.5}
			\fill[pattern=dots] (\x-1+\dx,\y-1+\dy) rectangle ++(0.5,0.5);
			
			\foreach \dx in {0,0.5}
			\foreach \dy in {0,0.5}
			\fill[pattern=north east lines] (1+\dx,2+\dy) rectangle ++(0.5,0.5);
			
			\foreach \x in {0,...,8} \draw (\x/2,0) -- (\x/2,4);
			\foreach \y in {0,...,8} \draw (0,\y/2) -- (4,\y/2);
			
			\node at (2,-0.7) {\textcolor{white}{$\tilde{G}_{\ell}$}};
		\end{tikzpicture}
	\end{center}
	\caption{Visualisation of the heuristic for determining a partition of $\tilde{G}_{\ell}\cap\mathcal{G}$ and $\tilde{G}_{\ell}\setminus\mathcal{G}$ into cuboids. The grey squares represent the active cells $G_{\ell} \cap \mathcal{G}$ or $\tilde{G}_{\ell}\cap\mathcal{G}$. The coloured squares represent partition $\tilde{G}_{\ell}\cap\mathcal{G} = \biguplus^{n_{\tilde{G}_{\ell}\cap\mathcal{G}}}_{j = 1} \mathbf{q}^{\tilde{G}_{\ell}\cap\mathcal{G}}_j$ with $n_{\tilde{G}_{\ell}\cap\mathcal{G}} = 5$, and the white squares with patterns represent partition $\tilde{G}_{\ell}\setminus\mathcal{G} = \biguplus^{n_{\tilde{G}_{\ell}\setminus\mathcal{G}}}_{j = 1} \mathbf{q}^{\tilde{G}_{\ell}\setminus\mathcal{G}}_j$ with $n_{\tilde{G}_{\ell}\setminus\mathcal{G}} = 2$.}
	\label{fig:integration}
\end{figure}
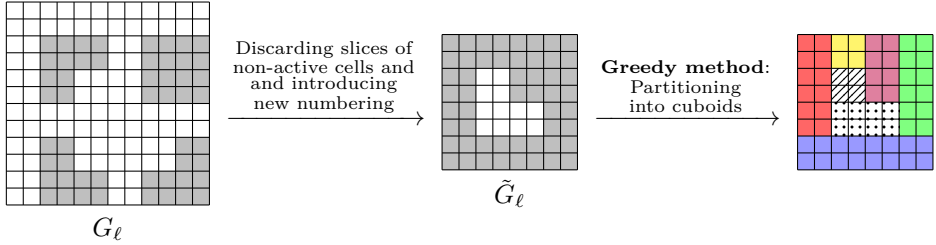

To overcome the second issue, we seek to reduce the number of summands in \cref{eq:level_mass_one_cell}. In a non-hierarchical, non-adaptive setting this is not an issue, because we integrate only once over the entire mesh (all cells), represented by the domain $\hat{\Omega} = \left[ 0, 1 \right]^3$, which has a tensor-product structure by definition and therefore univariate integration can be performed over $\left[ 0, 1 \right]$. In the hierarchical case, however, $G_{\ell}\cap\mathcal{G}$ generally lacks such a structure, even though the mesh of each level $G_{\ell}$ posseses one. Our strategy is therefore to identify a minimal collection of tensor-product structures - each called a cell cuboid - within the active cells $G_{\ell}\cap\mathcal{G}$ and the non-active cells $G_{\ell}\setminus\mathcal{G}$. Integrating over these larger unions of cells reduces the number of summands in \cref{eq:level_mass_one_cell}. This cuboid detection heuristic is illustrated in \cref{fig:integration} for an example in two dimensions for simplicity. 

Before these cell cuboids are detected, we discard slices of cells that are non-active cells inside $G_{\ell}$, as before for the basis functions. This yields a reduced mesh $\tilde{G}_{\ell}$ with a reduced index set $\mathcal{J}_{\tilde{G}_{\ell}} = \bigtimes^3_{d = 1} \left\{ 1, \ldots, \tilde{m}^{(d)}_{\ell} \right\}$, which still has a tensor product structure. This offers two benefits: connected domains of active cells merge, enabling larger active-cell cuboids, and the set of non-active cells $\tilde{G}_{\ell}\setminus\mathcal{G}$ shrinks so that fewer cells must be processed during cuboid detection. 

In the following a cell cuboid of level $\ell = 0, \ldots, L-1$ is defined as a disjoint union of cells such that the area they cover can be represented as a Cartesian product, i.e.
\begin{equation}
	\label{eq:cell_cuboid}
	\begin{gathered}
		\mathbf{q} = \biguplus_{\mathbf{j} \in \mathcal{J}_{\mathbf{q}}} q_{\mathbf{j}} = \mathbf{q}^{(1)} \times \mathbf{q}^{(2)} \times \mathbf{q}^{(3)} \subset \tilde{G}_{\ell}, \\
		\mathcal{J}_{\mathbf{q}} = \bigtimes^{3}_{d = 1} \left\{ j^{(d)}, j^{(d)} + 1, \ldots, j^{(d)} + n^{(d)}_{\mathbf{q}} - 1 \right\} \subset \mathcal{J}_{\tilde{G}_{\ell}}, \\
		\mathbf{q}^{(d)} = \biguplus^{n^{(d)}_{\mathbf{q}} - 2}_{k = 0} \left[ \hat{\xi}^{(d)}_{\ell, j^{(d)} + k}, \hat{\xi}^{(d)}_{\ell, j^{(d)} + k +1} \right],
	\end{gathered}
\end{equation}
where we point out that $\mathbf{q}^{(d)} \subset \left[ 0, 1 \right]$ but $\mathbf{q}^{(d)}$ does not need to be one closed interval. It can also be a set of intervals, since it is the disjoint union of the knot spans of cells which may not be neighbors in the original numbering. We point out that the multi-index set of the cell cuboid $\mathcal{J}_{\mathbf{q}}$ also has a tensor product structure.\\
It is immediately apparent that $\tilde{G}_{\ell}\cap\mathcal{G}$ and $\tilde{G}_{\ell}\setminus\mathcal{G}$ can be represented as disjoint union of such cuboids, i.e
\begin{equation}
	\label{eq:cuboid_partition}
	\tilde{G}_{\ell}\cap\mathcal{G} = \biguplus^{n_{\tilde{G}_{\ell}\cap\mathcal{G}}}_{j = 1} \mathbf{q}^{\tilde{G}_{\ell}\cap\mathcal{G}}_j, \quad \tilde{G}_{\ell}\setminus\mathcal{G} = \biguplus^{n_{\tilde{G}_{\ell}\setminus\mathcal{G}}}_{j = 1} \mathbf{q}^{\tilde{G}_{\ell}\setminus\mathcal{G}}_j,
\end{equation}
where $n_{\tilde{G}_{\ell}\cap\mathcal{G}}$ and $n_{\tilde{G}_{\ell}\setminus\mathcal{G}}$ are the number of cuboids. 

Since each cuboid $\mathbf{q}^{\tilde{G}_{\ell}\cap\mathcal{G}}_j \subset \tilde{G}_{\ell}\cap\mathcal{G}$ and $\mathbf{q}^{\tilde{G}_{\ell}\setminus\mathcal{G}}_j \subset \tilde{G}_{\ell}\setminus\mathcal{G}$ can contain several cells, we have $n_{\tilde{G}_{\ell}\cap\mathcal{G}} \leq \left| \tilde{G}_{\ell}\cap\mathcal{G} \right|$ and $n_{\tilde{G}_{\ell}\setminus \mathcal{G}} \leq \left| \tilde{G}_{\ell} \setminus \mathcal{G} \right|$, which means the number of cell cuboids can be at most as high as the number of active (or non-active) cells, but in many scenarios the number of cuboids is many times smaller, if the set of active cells already has a cuboid-like structure. We refer to the representation in \cref{eq:cuboid_partition} as a partition of $\tilde{G}_{\ell}\cap\mathcal{G}$ and $\tilde{G}_{\ell}\setminus\mathcal{G}$ into cuboids.

\begin{algorithm}
	\caption{Greedy method to detect cuboids in an index set}
	\label{alg:Cuboid_detection}
	\begin{algorithmic}[1]
		\STATE \textbf{Input:} A index set $\mathcal{I} \subset \bigtimes^3_{d = 1} \left\{ 1, \ldots, m^{(d)} \right\}$ which represents active (or non-active) elements with little-endian order 
		\STATE \textbf{Output:} A list of index cuboids $\left\{ \mathcal{I}_1, \ldots, \mathcal{I}_n \right\}$ \\[3mm]
		\STATE $k \gets 1$
		\WHILE{$\mathcal{I} \neq \emptyset$}
		\STATE Choose first index as starting point $\left( i^{(1)}_{s}, i^{(2)}_{s}, i^{(3)}_{s} \right) \in \mathcal{I}$
		\STATE $\left( i^{(1)}_{e}, i^{(2)}_{e}, i^{(3)}_{e} \right) \gets \left( i^{(1)}_{s}, i^{(2)}_{s}, i^{(3)}_{s} \right)$
		\WHILE{$i^{(1)}_{e} \leq m^{(1)}$ \AND $\left( i^{(1)}_{s} \colon i^{(1)}_{e}, i^{(2)}_{s}, i^{(3)}_{s} \right) \subset \mathcal{I}$}
		\STATE $i^{(1)}_{e} \gets i^{(1)}_{e}+1$
		\ENDWHILE 
		\WHILE{$i^{(2)}_{e} \leq m^{(2)}$ \AND $\left( i^{(1)}_{s} \colon i^{(1)}_{e}, i^{(2)}_{s} \colon i^{(2)}_{e}, i^{(3)}_{s} \right) \subset \mathcal{I}$}
		\STATE $i^{(2)}_{e} \gets i^{(2)}_{e}+1$
		\ENDWHILE 
		\WHILE{$i^{(3)}_{e} \leq m^{(3)}$ \AND $\left( i^{(1)}_{s} \colon i^{(1)}_{e}, i^{(2)}_{s} \colon i^{(2)}_{e}, i^{(3)}_{s} \colon i^{(3)}_{e} \right) \subset \mathcal{I}$}
		\STATE $i^{(3)}_{e} \gets i^{(3)}_{e}+1$
		\ENDWHILE 
		\STATE $\mathcal{I}_k \gets \left( i^{(1)}_{s} \colon i^{(1)}_{e}, i^{(2)}_{s} \colon i^{(2)}_{e}, i^{(3)}_{s} \colon i^{(3)}_{e} \right)$
		\STATE $\left\{ \mathcal{I}_1, \ldots, \mathcal{I}_{k-1} \right\} \gets \left\{ \mathcal{I}_1, \ldots, \mathcal{I}_{k-1}, \mathcal{I}_k \right\}$
		\STATE $\mathcal{I} \gets \mathcal{I} \setminus \mathcal{I}_k$
		\STATE $k \gets k + 1$ 
		\ENDWHILE
	\end{algorithmic}
\end{algorithm}

The aim is to find a partition \cref{eq:cuboid_partition} of $\tilde{G}_{\ell}\cap\mathcal{G}$ and $\tilde{G}_{\ell}\setminus\mathcal{G}$ into as few cuboids as possible, in order to minimize the number of summands for constructing $\mathbf{M}_\ell$, which means we want to have $n_{\tilde{G}_{\ell}\cap\mathcal{G}}$ or $n_{\tilde{G}_{\ell}\setminus\mathcal{G}}$ in \cref{eq:cuboid_partition} as small as possible. In our work we use for that a greedy method whose pseudocode is shown in \cref{alg:Cuboid_detection}. Let $\mathcal{J}_{\tilde{G}_{\ell} \cap \mathcal{G}}$ be the set of indices of active cells $\tilde{G}_{\ell}\cap\mathcal{G}$, which we assume are ordered in little-endian (column-major) order and $\mathcal{J}_{\tilde{G}_{\ell} \cap \mathcal{G}}$ will be updated in each iteration of this detection process. The method takes the first index $\left( j^{(1)}_{s}, j^{(2)}_{s}, j^{(3)}_{s} \right) \in \mathcal{J}_{\tilde{G}_{\ell} \cap \mathcal{G}}$ and checks how many cells are active in the fiber of dimension $d = 1$, which means we start with $\left( j^{(1)}_{s}, j^{(2)}_{s}, j^{(3)}_{s} \right)$ and increase the first component $j^{(1)}$ of the multi-index as long as the cell with index $\left( j^{(1)}, j^{(2)}_{s}, j^{(3)}_{s} \right)$ is active, i.e.~$\left( j^{(1)}, j^{(2)}_{s}, j^{(3)}_{s} \right) \in \mathcal{J}_{\tilde{G}_{\ell} \cap \mathcal{G}}$. If $\left( j^{(1)}, j^{(2)}_{s}, j^{(3)}_{s} \right)$ belongs to a non-active cell, i.e.~$\left( j^{(1)}, j^{(2)}_{s}, j^{(3)}_{s} \right) \in \mathcal{J}_{\tilde{G}_{\ell} \setminus \mathcal{G}}$, we stop and we get a fiber of indices $\left( j^{(1)}_{s} \colon j^{(1)}_{e}, j^{(2)}_{s}, j^{(3)}_{s} \right)\subset \mathcal{J}_{\tilde{G}_{\ell} \cap \mathcal{G}}$ of active cells, where from now on $j^{(d)}_{s} \colon j^{(d)}_{e}$ denotes $\left\{ j^{(d)}_{s}, j^{(1)}_{s}+1, \ldots, j^{(d)}_{e} \right\}$ and to be consistent with the notation in \cref{eq:cell_cuboid}, $j^{(d)}_{e} = j^{(d)}_{s} + n^{(d)}_{\mathbf{q}} - 1$ holds. Next, the method checks how many cells in dimension $d = 2$ are active that are already active in dimension $d = 1$, by increasing the second component $j^{(2)}$ of the fiber-index $\left( j^{(1)}_{s} \colon j^{(1)}_{e}, j^{(2)}, j^{(3)}_{s} \right)$, starting from $j^{(2)}_{s}$, and checking whether these fibers also belong to the active cells, i.e.~$\left( j^{(1)}_{s} \colon j^{(1)}_{e}, j^{(2)}, j^{(3)}_{s} \right) \subset \mathcal{J}_{\tilde{G}_{\ell} \cap \mathcal{G}}$. This results in the indices $\left( j^{(1)}_{s} \colon j^{(1)}_{e}, j^{(2)}_{s}\colon j^{(2)}_{e}, j^{(3)}_{s} \right) \subset \mathcal{J}_{\tilde{G}_{\ell} \cap \mathcal{G}}$ of a cell slice. Lastly, the method checks how many cells in dimension $d = 3$ are active that are already active in dimensions $d = 1$ and $d = 2$, by increasing the third component $j^{(3)}$ of the slice-index $\left( j^{(1)}_{s} \colon j^{(1)}_{e}, j^{(2)}_{s}\colon j^{(2)}_{e}, j^{(3)} \right)$, starting from $j^{(3)}_{s}$, and checking whether these slices also belong to the active cells, i.e.~$\left( j^{(1)}_{s} \colon j^{(1)}_{e}, j^{(2)}_{s}\colon j^{(2)}_{e}, j^{(3)} \right) \subset \mathcal{J}_{\tilde{G}_{\ell} \cap \mathcal{G}}$. This results in the indices $\left( j^{(1)}_{s} \colon j^{(1)}_{e}, j^{(2)}_{s}\colon j^{(2)}_{e}, j^{(3)}_{s} : j^{(3)}_{s} \right) \subset \mathcal{J}_{\tilde{G}_{\ell} \cap \mathcal{G}}$ of a cell cuboid, i.e.~$\mathcal{J}_{\mathbf{q}^{\tilde{G}_{\ell}\cap\mathcal{G}}} = \left( j^{(1)}_{s} \colon j^{(1)}_{e}, j^{(2)}_{s} \colon j^{(2)}_{e}, j^{(3)}_{s} \colon j^{(3)}_{e} \right)$. Once the indices $\mathcal{J}_{\mathbf{q}^{\tilde{G}_{\ell}\cap\mathcal{G}}}$ of a cuboid $\mathbf{q} \subset \tilde{G}_{\ell} \cap \mathcal{G}$ have been determined, they are no longer considered, i.e. $\mathcal{J}_{\tilde{G}_{\ell} \cap \mathcal{G}} \mapsto \mathcal{J}_{\tilde{G}_{\ell} \cap \mathcal{G}} \setminus \mathcal{J}_{\mathbf{q}^{\tilde{G}_{\ell}\cap\mathcal{G}}}$, and the method starts from the beginning with the updates index set $\mathcal{J}_{\tilde{G}_{\ell} \cap \mathcal{G}}$, repeating this process until $\tilde{G}_{\ell} \cap \mathcal{G}$ can be represented as a partition of cuboids shown in \cref{eq:cuboid_partition}.

Since we are also interested in a partition of the set of non-active cells $\tilde{G}_{\ell} \setminus \mathcal{G}$ into cuboids, we apply the same method described by \cref{alg:Cuboid_detection} to the set of indices of the non-active cells cells $\mathcal{J}_{\tilde{G}_{\ell} \setminus \mathcal{G}}$. 

\cref{alg:Cuboid_detection} shares the core idea of a greedy method in that it generates an optimal solution locally at each step and continues successively, which is why we point out that the partition of cuboids generated in this way may not be optimal, i.e.~it may be possible to find a partition with fewer cuboids.

Once we have determined a partition of the sets $\tilde{G}_{\ell} \cap \mathcal{G}$ and $\tilde{G}_{\ell} \setminus \mathcal{G}$ into cuboids as in \cref{eq:cuboid_partition} using the indices computed by \cref{alg:Cuboid_detection}, we can perform the integration step. If $n_{\tilde{G}_{\ell}\cap\mathcal{G}} \leq n_{\tilde{G}_{\ell}\setminus\mathcal{G}} + 1$ we set
\begin{equation}
	\label{eq:level_mass_active_cuboids}
	\mathbf{M}_\ell = \sum^{n_{\tilde{G}_{\ell}\cap\mathcal{G}}}_{j = 1} \sum_{r=1}^R \mathbf{M}_{\ell, r, \mathbf{q}_j} = \sum^{n_{\tilde{G}_{\ell}\cap\mathcal{G}}}_{j = 1} \sum_{r=1}^R \bigotimes_{d=1}^3 M_{\ell, r, \mathbf{q}^{\left( d \right)}_j}^{\left( d \right)} \in \mathbb{R}^{\left( \tilde{n}^{(1)}_{\ell}, \tilde{n}^{(2)}_{\ell}, \tilde{n}^{(3)}_{\ell} \right) \times \left( \tilde{n}^{(1)}_{\ell}, \tilde{n}^{(2)}_{\ell}, \tilde{n}^{(3)}_{\ell} \right)}
\end{equation}
where the univariate integral of an active cuboid $\mathbf{q}_j = \mathbf{q}^{(1)}_j \times \mathbf{q}^{(2)}_j \times \mathbf{q}^{(3)}_j \subset \tilde{G}_{\ell} \cap \mathcal{G}$, $j = 1, \ldots, n_{\tilde{G}_{\ell}\cap\mathcal{G}}$, in dimension $d = 1, 2, 3$ is defined by 
\begin{equation}
	\label{eq:univariate_cuboid_integral}
	\begin{gathered}
		M_{\ell, r, \mathbf{q}^{\left( d \right)}_j}^{\left( d \right)} = \int_{\mathbf{q}^{(d)}_j} \left( w_r^{(d)}\cdot \hat{B}^{(d)} \left( \hat{x}^{(d)} \right) \right) \, \tilde{B}^{(d)}_{\ell} \left( \hat{x}^{(d)} \right) \otimes \tilde{B}^{(d)}_{\ell} \left( \hat{x}^{(d)} \right) \D \hat{x}^{(d)}, \\
		\tilde{B}^{(d)}_{\ell} \left( \hat{x}^{(d)} \right) = \left[ \tilde{\beta}^{(d)}_{\ell, 1} \left( \hat{x}^{(d)} \right), \ldots, \tilde{\beta}^{(d)}_{\ell, \tilde{n}^{(d)}_{\ell}} \left( \hat{x}^{(d)} \right) \right]^\top \in \mathbb{R}^{\tilde{n}^{(d)}_{\ell}}.
	\end{gathered}
\end{equation}
We point out that the univariate quadrature over non-closed intervals, which $\mathbf{q}^{\left( d \right)}_j$ may be, does not cause any problems, since the quadrature is performed successively over the single knot spans.

If $n_{\tilde{G}_{\ell}\cap\mathcal{G}} > n_{\tilde{G}_{\ell}\setminus\mathcal{G}} + 1$ then we set
\begin{equation}
	\label{eq:level_mass_non_active_cuboids}
	\mathbf{M}_\ell = \mathbf{M}_{\ell, \tilde{G}_{\ell}} - \sum^{n_{\tilde{G}_{\ell} \setminus \mathcal{G}}}_{j = 1} \sum_{r=1}^R \bigotimes_{d=1}^3 M_{\ell, r, \mathbf{q}^{\left( d \right)}_j}^{\ell, \left( d \right)}
\end{equation}
where $\mathbf{M}_{\ell, \tilde{G}_{\ell}}$ is the mass tensor of the splines $\tilde{\mathcal{B}}_{\ell}$, where the integration is done over all cells with new numbering $\tilde{G}_{\ell}$, which has a tensor product structure, so $\mathbf{M}_{\ell, \tilde{G}_{\ell}}$ can be computed as in \cref{eq:massFinal}, but with adjusted intervals, and $M_{\ell, r, \mathbf{q}_j}^{\ell, \left( d \right)}$ is defined by \cref{eq:univariate_cuboid_integral} for a non-active cuboid $\mathbf{q}_j = \mathbf{q}^{(1)}_j \times \mathbf{q}^{(2)}_j \times \mathbf{q}^{(3)}_j \subset \tilde{G}_{\ell} \setminus \mathcal{G}$, $j = 1, \ldots, n_{\tilde{G}_{\ell}\setminus\mathcal{G}}$.

The choice between \eqref{eq:level_mass_active_cuboids} and \eqref{eq:level_mass_non_active_cuboids} is driven by the cost of the univariate quadratures in the low-rank assembly. Each cuboid contributes one rank-$R$ tensor $\sum_{r=1}^R \bigotimes_{d=1}^3 M^{(d)}_{\ell,r,\mathbf{q}^{(d)}}$, and each Kronecker product $\bigotimes_{d=1}^3 M^{(d)}_{\ell,r,\mathbf{q}^{(d)}}$ requires three one–dimensional integrals as in \eqref{eq:univariate_cuboid_integral}. Hence, the total number of rank-$R$ summands is $n_{\tilde{G}_{\ell}\cap\mathcal{G}}$ for \eqref{eq:level_mass_active_cuboids} and $n_{\tilde{G}_{\ell}\setminus\mathcal{G}} + 1$ for \eqref{eq:level_mass_non_active_cuboids}, where the additional “$1$” accounts for the full-domain tensor $\mathbf{M}_{\ell,\tilde{G}_{\ell}}$. Since the dominant work scales linearly with the number of such summands, we pick the variant with fewer terms. We note that even though the tensors in \cref{eq:level_mass_active_cuboids} and \cref{eq:level_mass_non_active_cuboids} are written in CP form, in our implementation we convert every Kronecker product (a rank-one tensor) to TT format, accumulate these rank-one tensors in TT, and apply TT rounding~\cite{Oseledets:2011} to keep the ranks low. This heuristic therefore minimizes both quadrature work and transient rank growth during TT accumulation. The two formulas are equivalent up to quadrature error, and the heuristic is a purely complexity-aware selection performed independently on each level~$\ell$.

\subsection{Low-rank Truncation Operator}
\label{subsection:Low-rank_Truncation_Operator}
After constructing a mass matrix for each level $\mathbf{M}_\ell$ by integrating all basis functions of this level over the active cells \cref{eq:level_mass}, the global hierarchical mass matrix $\mathbf{M}$ in~\cite{GarauVazquez:2018} is constructed by
\begin{equation}
	\label{eq:hierarchical_mass_matrix}
	\mathbf{M} = \sum_{\ell = 0}^{L-1} 
	\begin{bmatrix}
		\mathbf{C}_{\ell} \\
		\mathbf{0}
	\end{bmatrix} 
	\mathbf{M}_\ell 
	\begin{bmatrix}
		\mathbf{C}_{\ell} & \mathbf{0}
	\end{bmatrix} 
\end{equation}
where $\mathbf{0}$ denotes a zero matrix of suitable size and $\mathbf{C}_{\ell}$ is used to perform both truncation and the inclusion of basis functions of level $\ell$ into the basis $\mathcal{T}$, which is defined by
\begin{enumerate}
	\item[1)] $\mathbf{C}_{0} = \mathbf{J}_{0}$
	\item[2)] $\mathbf{C}_{\ell+1} = 
	\begin{bmatrix}
		\mathbf{C}^{\left( \ell + 1, \ell \right)} \mathbf{C}_{\ell} & \mathbf{J}_{\ell+1}
	\end{bmatrix}$ for $\ell = 0, \ldots, L-2$
\end{enumerate}
where $\mathbf{C}^{\left( \ell + 1, \ell \right)}$ is the coefficient matrix of the two-scale relation with truncation defined by \cref{eq:truncation_coefficients} and $\mathbf{J}_{\ell}$ denotes a Boolean matrix, which realizes the inclusion of active functions of level $\ell$, i.e.~from $\mathcal{B}_{\ell} \cap \mathcal{T}$. \\
The matrix obtained by \cref{eq:hierarchical_mass_matrix} has a $L \times L$ block structure, i.e.
\begin{equation}
	\label{eq:hierarchical_mass_matrix_block_structure}
	\mathbf{M} = 
	\begin{bmatrix}
		\mathbf{M}_{0,0} & \cdots & \mathbf{M}_{0,L-1} \\
		\vdots & & \vdots \\
		\mathbf{M}_{L-1,0} & \cdots & \mathbf{M}_{L-1,L-1}
	\end{bmatrix},
\end{equation}
where the individual blocks for $\ell, \ell' = 0, \ldots, L-1$ are given by 
\begin{equation}
	\label{eq:hierarchical_mass_matrix_blocks}
	\mathbf{M}_{\ell, \ell'} = \mathbf{J}^\top_{\ell} \left( \sum_{k = \max \left( \left\{ \ell, \ell' \right\} \right)}^{L-1} \left( \prod_{m = \ell}^{k-1} {\mathbf{C}^{\left( m + 1, m \right)}}^\top \right) \mathbf{M}_k \left( \prod_{m = \ell'}^{k-1} {\mathbf{C}^{\left( m + 1, m \right)}}^\top \right)^\top \right) \mathbf{J}_{\ell'}.
\end{equation}

Our aim is to convert the global hierarchical mass matrix $\mathbf{M}$ into a low-rank (TT) format usable by \textsc{MATLAB}’s \textsc{TT-Toolbox}~\cite{tt-toolbox}. Since the toolbox provides block solvers, such as TT-GMRES, we adopt the block structure \cref{eq:hierarchical_mass_matrix_block_structure} and therefore represent each block $\mathbf{M}_{\ell,\ell'}$ defined by \eqref{eq:hierarchical_mass_matrix_blocks} in low-rank form. In the following, we assume that $\mathbf{M}_\ell$ has been computed via \cref{eq:level_mass_active_cuboids,eq:level_mass_non_active_cuboids}. To enable efficient assembly of the individual blocks $\mathbf{M}_{\ell,\ell'}$, we present efficient low-rank representations of the operators $\mathbf{C}^{\left( \ell + 1, \ell \right)}$ and $\mathbf{J}_{\ell}$.

In our implementation, we compute the blocks in two loops. First, we add up the mass matrices multiplied by the coefficients of the two-scale relation with truncation, i.e.
\begin{equation}
	\label{eq:hierarchical_mass_matrix_blocks_first_loop}
	\tilde{\mathbf{M}}_{\ell, \ell'} = \sum_{k = \max \left( \left\{ \ell, \ell' \right\} \right)}^{L-1} \left( \prod_{m = \ell}^{k-1} {\mathbf{C}^{\left( m + 1, m \right)}}^\top \right) \mathbf{M}_k \left( \prod_{m = \ell'}^{k-1} {\mathbf{C}^{\left( m + 1, m \right)}}^\top \right)^\top
\end{equation}
and in the second loop, we select the active basis functions of the respective levels, i.e.
\begin{equation}
	\label{eq:hierarchical_mass_matrix_blocks_second_loop}
	\mathbf{M}_{\ell, \ell'} = \mathbf{J}^\top_{\ell} \tilde{\mathbf{M}}_{\ell, \ell'} \mathbf{J}_{\ell'}
\end{equation}
for $\ell, \ell' = 0, \ldots, L-1$.

Since $\mathbf{M}_\ell$ in \cref{eq:hierarchical_mass_matrix_blocks_first_loop} is of size $\left( \tilde{n}^{(1)}_{\ell}, \tilde{n}^{(2)}_{\ell}, \tilde{n}^{(3)}_{\ell} \right) \times \left( \tilde{n}^{(1)}_{\ell}, \tilde{n}^{(2)}_{\ell}, \tilde{n}^{(3)}_{\ell} \right)$, the tensors $\mathbf{C}^{\left( \ell + 1, \ell \right)} \in \mathbb{R}^{\left( \tilde{n}^{(1)}_{\ell+1}, \tilde{n}^{(2)}_{\ell+1}, \tilde{n}^{(3)}_{\ell+1} \right) \times \left( \tilde{n}^{(1)}_{\ell}, \tilde{n}^{(2)}_{\ell}, \tilde{n}^{(3)}_{\ell} \right)}$ take over the numbering of basis functions introduced in \cref{subsection:domain_wise_integration}, which means in $\mathbf{C}^{\left( \ell + 1, \ell \right)}$ the coefficients for the two-scale relation with truncation restricted to the bases $\tilde{\mathcal{B}}_{\ell}$ and $\tilde{\mathcal{B}}_{\ell+1}$ are stored. Respectively $\tilde{\mathbf{M}}_{\ell, \ell'}$ is of size $\left( \tilde{n}^{(1)}_{\ell}, \tilde{n}^{(2)}_{\ell}, \tilde{n}^{(3)}_{\ell} \right) \times \left( \tilde{n}^{(1)}_{\ell'}, \tilde{n}^{(2)}_{\ell'}, \tilde{n}^{(3)}_{\ell'} \right)$ for $\ell, \ell' = 0, \ldots, L-1$. We point out that in $\tilde{\mathcal{B}}_{\ell}$ the active functions of level $\ell$ are included, i.e.~$\tilde{\mathcal{B}}_{\ell} \cap \mathcal{T} = \mathcal{B}_{\ell} \cap \mathcal{T}$, since the support of active basis functions also intersects with or lies in the set of active cells $G_{\ell} \cap \mathcal{G}$. \\

We want to efficiently assemble the coefficient tensor of the two-scale relation with truncation $\mathbf{C}^{\left( \ell + 1, \ell \right)} \in \mathbb{R}^{\left( \tilde{n}^{(1)}_{\ell+1}, \tilde{n}^{(2)}_{\ell+1}, \tilde{n}^{(3)}_{\ell+1} \right) \times \left( \tilde{n}^{(1)}_{\ell}, \tilde{n}^{(2)}_{\ell}, \tilde{n}^{(3)}_{\ell} \right)}$ in low-rank form. We exploit the fact that the coefficient of the two-scale relation \cref{eq:two_scale_relation_coefficients} can be represent as a rank-one tensor, since it inherits the tensor product structure of the bases $\tilde{\mathcal{B}}_{\ell + 1}$ and $\tilde{\mathcal{B}}_{\ell}$. Unfortunately, the coefficient tensor $\mathbf{C}^{\left( \ell + 1, \ell \right)}$ defined by \cref{eq:truncation_coefficients} in general can not be represented as a rank-one tensor, since its entries regarding the basis functions contained in $\tilde{\mathcal{B}}_{\ell+1} \cap \mathcal{T}_{\ell+1}$ are set to zero and the remaining entries can in general not be represented by a tensor product.

To overcome this issue we apply a similar heuristic as described in \cref{subsection:domain_wise_integration} for assembling the mass matrices for each level. We want to find a partition of $\tilde{\mathcal{B}}_{\ell+1} \cap \mathcal{T}_{\ell+1}$ or $\tilde{\mathcal{B}}_{\ell+1} \setminus \mathcal{T}_{\ell+1}$ (for the splines of this subspace the entries in $\mathbf{C}^{\left( \ell + 1, \ell \right)}$ are not zero) into a set spline cuboids in order to successively construct the tensor $\mathbf{C}^{\left( \ell + 1, \ell \right)}$. In the following a spline cuboid of level $\ell = 0, \ldots, L-1$ is defined as a disjoint union of splines of $\tilde{\mathcal{B}}_{\ell}$ which has the following tensor product structure
\begin{equation}
	\label{eq:spline_cuboid}
	\begin{gathered}
		\boldsymbol{\beta} = \biguplus_{\mathbf{i} \in \mathcal{I}_{\boldsymbol{\beta}}} \tilde{\beta}_{\ell, \mathbf{i}} = \biguplus_{\substack{\mathbf{i} \in \mathcal{I}_{\boldsymbol{\beta}} \\ \mathbf{i} = \left( i^{(1)}, i^{(2)}, i^{(3)} \right)}} \tilde{\beta}^{(1)}_{\ell, i^{(1)}} \tilde{\beta}^{(2)}_{\ell, i^{(2)}} \tilde{\beta}^{(3)}_{\ell, i^{(3)}} = \boldsymbol{\beta}^{(1)} \otimes \boldsymbol{\beta}^{(2)} \otimes \boldsymbol{\beta}^{(3)} \subset \tilde{\mathcal{B}}_{\ell}, \\
		\mathcal{I}_{\boldsymbol{\beta}} = \bigtimes^{3}_{d = 1} \left\{ i^{(d)}, i^{(d)} + 1, \ldots, i^{(d)} + n^{(d)}_{\boldsymbol{\beta}} - 1 \right\} \subset \mathcal{I}_{\tilde{\mathcal{B}}_{\ell}}, \\
		\boldsymbol{\beta}^{(d)} = \left\{ \tilde{\beta}^{(d)}_{\ell, i^{(d)}}, \tilde{\beta}^{(d)}_{\ell, i^{(d)} + 1}, \ldots, \tilde{\beta}^{(d)}_{\ell, i^{(d)} + n^{(d)}_{\boldsymbol{\beta}} - 1} \right\} \subset \tilde{\mathcal{B}}^{(d)}_{\ell}.
	\end{gathered}
\end{equation}
As for the active and non-active cells of level $\ell$ in \cref{eq:cuboid_partition}, $\tilde{\mathcal{B}}_{\ell+1} \cap \mathcal{T}_{\ell+1}$ and $\tilde{\mathcal{B}}_{\ell+1} \setminus \mathcal{T}_{\ell+1}$ can be represented as a disjoint union of these spline cuboids via 
\begin{equation}
	\label{eq:spline_partition_truncation}
	\tilde{\mathcal{B}}_{\ell+1} \cap \mathcal{T}_{\ell+1} = \biguplus^{n_{\tilde{\mathcal{B}}_{\ell+1} \cap \mathcal{T}_{\ell+1}}}_{i = 1} \boldsymbol{\beta}^{\tilde{\mathcal{B}}_{\ell+1} \cap \mathcal{T}_{\ell+1}}_i, \quad \tilde{\mathcal{B}}_{\ell+1} \setminus \mathcal{T}_{\ell+1} = \biguplus^{n_{\tilde{\mathcal{B}}_{\ell+1} \setminus \mathcal{T}_{\ell+1}}}_{i = 1} \boldsymbol{\beta}^{\tilde{\mathcal{B}}_{\ell+1} \setminus \mathcal{T}_{\ell+1}}_i,
\end{equation}
where $n_{\tilde{\mathcal{B}}_{\ell+1} \cap \mathcal{T}_{\ell+1}}$ and $n_{\tilde{\mathcal{B}}_{\ell+1} \setminus \mathcal{T}_{\ell+1}}$ are the number of cuboids and we refer to \cref{eq:spline_partition_truncation} as a partition of $\tilde{\mathcal{B}}_{\ell+1} \cap \mathcal{T}_{\ell+1}$ and $\tilde{\mathcal{B}}_{\ell+1} \setminus \mathcal{T}_{\ell+1}$ into cuboids. 

As in the case of active and inactive cells of level $\ell$ in \cref{subsection:domain_wise_integration}, the aim is to find a partition of $\tilde{\mathcal{B}}_{\ell+1} \cap \mathcal{T}_{\ell+1}$ and $\tilde{\mathcal{B}}_{\ell+1} \setminus \mathcal{T}_{\ell+1}$ into as few cuboids as possible. To achieve this, we use a variant of the greedy method represented by \cref{alg:Cuboid_detection}, which returns the indices of the respective spline cuboids. As for the cell cuboids the method starts with an index $\left( i^{(1)}_s, i^{(2)}_s, i^{(3)}_s \right) \in \mathcal{I}_{\tilde{\mathcal{B}}_{\ell+1} \cap \mathcal{T}_{\ell+1}}$ and checks how far it can extend the indices in all directions while staying inside $\tilde{\mathcal{B}}_{\ell+1} \cap \mathcal{T}_{\ell+1}$ to get an index cuboid $\mathcal{I}_{\boldsymbol{\beta}^{\tilde{\mathcal{B}}_{\ell+1} \cap \mathcal{T}_{\ell+1}}} = \left( i^{(1)}_s : i^{(1)}_e, i^{(2)}_s : i^{(2)}_e, i^{(3)}_s : i^{(3)}_e \right)$. These indices are subtracted from the set of indices $\mathcal{I}_{\tilde{\mathcal{B}}_{\ell+1} \cap \mathcal{T}_{\ell+1}} \mapsto \mathcal{I}_{\tilde{\mathcal{B}}_{\ell+1} \cap \mathcal{T}_{\ell+1}} \setminus\mathcal{I}_{\boldsymbol{\beta}^{\tilde{\mathcal{B}}_{\ell+1} \cap \mathcal{T}_{\ell+1}}}$, and the method continues with the next index from $\mathcal{I}_{\tilde{\mathcal{B}}_{\ell+1} \cap \mathcal{T}_{\ell+1}}$. The same is done with $\mathcal{I}_{\tilde{\mathcal{B}}_{\ell+1} \setminus \mathcal{T}_{\ell+1}}$ to determine a partition of $\tilde{\mathcal{B}}_{\ell+1} \setminus \mathcal{T}_{\ell+1}$ into cuboids.

When we have a partition of $\tilde{\mathcal{B}}_{\ell+1} \cap \mathcal{T}_{\ell+1}$ and $\tilde{\mathcal{B}}_{\ell+1} \setminus \mathcal{T}_{\ell+1}$ as in \cref{eq:spline_partition_truncation} and if $n_{\tilde{\mathcal{B}}_{\ell+1} \setminus \mathcal{T}_{\ell+1}} \leq n_{\tilde{\mathcal{B}}_{\ell+1} \cap \mathcal{T}_{\ell+1}} + 1$ we set
\begin{equation}
	\label{eq:truncation_plus}
	\mathbf{C}^{\left( \ell + 1, \ell \right)} = \sum_{i = 1}^{n_{\tilde{\mathcal{B}}_{\ell+1} \setminus \mathcal{T}_{\ell+1}}} \bigotimes^{3}_{d = 1} C_{\boldsymbol{\beta}^{(d)}_i}^{\left( \ell + 1, \ell \right), (d)}
\end{equation}
where the univariate coefficient matrix $C_{\boldsymbol{\beta}^{(d)}}^{\left( \ell + 1, \ell \right), (d)} \in \mathbb{R}^{\tilde{n}^{(d)}_{\ell+1} \times \tilde{n}^{(d)}_{\ell}}$ of dimension $d = 1, 2, 3$ for a spline cuboid $\boldsymbol{\beta} = \boldsymbol{\beta}^{(1)} \otimes \boldsymbol{\beta}^{(2)} \otimes \boldsymbol{\beta}^{(3)} \subset \tilde{\mathcal{B}}_{\ell+1}$ \cref{eq:spline_cuboid} is defined by 
\begin{equation}
	\label{eq:coefficients_cuboid}
	\left( C_{\boldsymbol{\beta}^{(d)}}^{\left( \ell + 1, \ell \right), (d)} \right)_{i,j} = 
	\begin{cases}
		\left( C^{\left( \ell + 1, \ell \right), (d)} \right)_{i, j} &: \text{ if } \tilde{\beta}^{(d)}_{\ell+1, i} \in \boldsymbol{\beta}^{(d)} \\
		0 &: \text{ otherwise}
	\end{cases}
\end{equation}
for $i = 1, \ldots, \tilde{n}^{(d)}_{\ell+1}$ and $j = 1, \ldots, \tilde{n}^{(d)}_{\ell}$, where $C^{\left( \ell + 1, \ell \right), (d)} \in \mathbb{R}^{\tilde{n}^{(d)}_{\ell+1} \times \tilde{n}^{(d)}_{\ell}}$ represents the univariate coefficient matrix in dimension $d$ of the two-scale relation \cref{eq:two_scale_relation_coefficients} restricted to the numbering of the bases $\tilde{\mathcal{B}}^{(d)}_{\ell}$ and $\tilde{\mathcal{B}}^{(d)}_{\ell+1}$.

If $n_{\tilde{\mathcal{B}}_{\ell+1} \setminus \mathcal{T}_{\ell+1}} > n_{\tilde{\mathcal{B}}_{\ell+1} \cap \mathcal{T}_{\ell+1}} + 1$ then we set
\begin{equation}
	\label{eq:truncation_minus}
	\mathbf{C}^{\left( \ell + 1, \ell \right)} = \mathbf{C}^{\left( \ell + 1, \ell \right)} - \sum^{n_{\tilde{\mathcal{B}}_{\ell+1} \cap \mathcal{T}_{\ell+1}}}_{i = 1} \bigotimes_{d=1}^3 C_{\boldsymbol{\beta}^{(d)}_i}^{\left( \ell + 1, \ell \right), (d)}
\end{equation}
where $\mathbf{C}^{\left( \ell + 1, \ell \right)} \in \mathbb{R}^{\left( \tilde{n}^{(1)}_{\ell+1}, \tilde{n}^{(2)}_{\ell+1}, \tilde{n}^{(3)}_{\ell+1} \right) \times \left( \tilde{n}^{(1)}_{\ell}, \tilde{n}^{(2)}_{\ell}, \tilde{n}^{(3)}_{\ell} \right)}$ is the coefficient tensor of the two-scale relation \cref{eq:two_scale_relation_coefficients} restricted to the numbering of the bases $\tilde{\mathcal{B}}_{\ell}$ and $\tilde{\mathcal{B}}_{\ell+1}$ and $C_{\boldsymbol{\beta}^{(d)}_i}^{\left( \ell + 1, \ell \right), (d)}$ is defined by \cref{eq:coefficients_cuboid} for $\boldsymbol{\beta}_i \subset \tilde{\mathcal{B}}_{\ell+1} \cap \mathcal{T}_{\ell+1}$ for $i = 1, \ldots, n_{\tilde{\mathcal{B}}_{\ell+1} \cap \mathcal{T}_{\ell+1}}$.

Using \cref{eq:truncation_plus,eq:truncation_minus}, we efficiently assemble the coefficient tensor $\mathbf{C}^{\left( \ell + 1, \ell \right)} \in \mathbb{R}^{\left(\tilde{n}^{(1)}_{\ell+1},\tilde{n}^{(2)}_{\ell+1},\tilde{n}^{(3)}_{\ell+1} \right)\times \left( \tilde{n}^{(1)}_{\ell},\tilde{n}^{(2)}_{\ell},\tilde{n}^{(3)}_{\ell} \right)}$, which in turn enables the evaluation of $\tilde{\mathbf{M}}_{\ell,\ell'}$ via \cref{eq:hierarchical_mass_matrix_blocks_first_loop} for all $\ell,\ell'=0,\ldots,L-1$.

As in \cref{subsection:domain_wise_integration}, the two different constructions \eqref{eq:truncation_plus} and \eqref{eq:truncation_minus} for $\mathbf{C}^{\left( \ell + 1, \ell \right)}$ are chosen purely for efficiency: we aim to minimize the number of Kronecker summands, since more summands both increase susceptibility to rounding errors and prolong the TT-based assembly (multiplications/additions and subsequent TT rounding). The number of spline cuboids in \eqref{eq:spline_partition_truncation} is not a reliable predictor of the true numerical rank of $\mathbf{C}^{\left( \ell + 1, \ell \right)}$. Nonetheless, each cuboid contributes at most one rank-one Kronecker term, so a natural upper bound for the effective rank is the number of cuboids (plus one in the subtractive case \eqref{eq:truncation_minus}). We therefore prefer the formulation with fewer summands to reduce assembly time and numerical error.

We note that in the case of hierarchical B-splines, where no truncation occurs, the assembly of \cref{eq:hierarchical_mass_matrix_blocks_first_loop} with the coefficient tensor $\mathbf{C}^{\left( \ell + 1, \ell \right)} \in \mathbb{R}^{\left( \tilde{n}^{(1)}_{\ell+1}, \tilde{n}^{(2)}_{\ell+1}, \tilde{n}^{(3)}_{\ell+1} \right) \times \left( \tilde{n}^{(1)}_{\ell}, \tilde{n}^{(2)}_{\ell}, \tilde{n}^{(3)}_{\ell} \right)}$ from \cref{eq:two_scale_relation_coefficients} is done with the new numbering of the basis, which is a rank-1 tensor. This means that the resulting blocks can have lower ranks than for THB-splines and can therefore be more memory-efficient under certain conditions.

\subsection{Gobal Hierarchical Mass Tensor}
\label{subsection:Gobal_Hierarchical_Mass_Tensor}
The final step is to assemble the global hierarchical mass tensor $\mathbf{M}$. We exploit the block structure in \cref{eq:hierarchical_mass_matrix_block_structure}: after computing for each pair $(\ell,\ell')$ the truncated level-coupled contribution $\tilde{\mathbf{M}}_{\ell,\ell'}$ via \cref{eq:hierarchical_mass_matrix_blocks_first_loop}, we have to restrict it to the active bases by applying $\mathbf{J}_\ell$ and $\mathbf{J}_{\ell'}$, cf.~\cref{eq:hierarchical_mass_matrix_blocks_second_loop}. Since the active sets $\mathcal{B}_\ell \cap \mathcal{T}$ generally lack a tensor-product structure, we reuse the cuboid partitioning strategies presented above to obtain a usable low-rank representation. Our objective is a TT-compatible global tensor $\mathbf{M}$ that encodes all active functions $\mathcal{T}$ while keeping ranks and assembly costs under control.

First, we introduce again a new numbering for the basis $\mathcal{B}_{\ell}$ using the slice-discarding heuristic from \cref{subsection:domain_wise_integration}. Since we are interested only in the active functions $\mathcal{B}_\ell \cap \mathcal{T}$, we discard all slices of B-splines which are not active, i.e.~are elemnts from $\mathcal{B}_{\ell} \setminus \mathcal{T}$, from $\mathcal{B}_{\ell}$ to get a new basis $\hat{\mathcal{B}}_{\ell}$ with a new numbering $\mathcal{I}_{\hat{\mathcal{B}}_{\ell}} = \bigtimes^3_{d = 1} \left\{ 1, \ldots, \hat{n}^{(d)}_{\ell} \right\}$, which has a tensor product structure $\hat{\mathcal{B}}_{\ell} = \hat{\mathcal{B}}^{(1)}_{\ell} \otimes \hat{\mathcal{B}}^{(2)}_{\ell} \otimes \hat{\mathcal{B}}^{(3)}_{\ell}$ where $\hat{\mathcal{B}}^{(d)}_{\ell} = \left\{ \hat{\beta}^{(d)}_{\ell, i} \colon \left[0, 1\right] \to \mathbb{R} \, \colon \, i = 1, \ldots, \hat{n}^{(d)}_{\ell} \right\}$. We have $\hat{\mathcal{B}}_\ell \subset \tilde{\mathcal{B}}_\ell \subset \mathcal{B}_\ell$, but $\hat n^{(d)}_\ell < n^{(d)}_\ell$ does not necessarily hold: if the active splines occupy all indices in mode $d$, no slice can be discarded in that direction. Moreover, to preserve the tensor product structure, $\hat{\mathcal{B}}_\ell$ may still contain some inactive functions.\\

We present two different approaches to assemble the individual blocks $\mathbf{M}_{\ell, \ell'}$. For that we need to determine a partition of $\hat{\mathcal{B}}_{\ell} \cap \mathcal{T}$ and $\hat{\mathcal{B}}_{\ell} \setminus \mathcal{T}$ into spline cuboids, which means we apply a variant of \cref{alg:Cuboid_detection} on $\mathcal{I}_{\hat{\mathcal{B}}_{\ell} \cap \mathcal{T}}$ and $\mathcal{I}_{\hat{\mathcal{B}}_{\ell} \setminus \mathcal{T}}$ to get
\begin{equation}
	\label{eq:spline_partition}
	\hat{\mathcal{B}}_{\ell} \cap \mathcal{T} = \biguplus^{n_{\hat{\mathcal{B}}_{\ell} \cap \mathcal{T}}}_{i = 1} \boldsymbol{\beta}^{\hat{\mathcal{B}}_{\ell} \cap \mathcal{T}}_i, \quad \hat{\mathcal{B}}_{\ell} \setminus \mathcal{T} = \biguplus^{n_{\hat{\mathcal{B}}_{\ell} \setminus \mathcal{T}}}_{i = 1} \boldsymbol{\beta}^{\hat{\mathcal{B}}_{\ell} \setminus \mathcal{T}}_i,
\end{equation}
where $\boldsymbol{\beta}^{\hat{\mathcal{B}}_{\ell} \cap \mathcal{T}}_i$, $i = 1, \ldots, n_{\hat{\mathcal{B}}_{\ell} \cap \mathcal{T}}$, and $\boldsymbol{\beta}^{\hat{\mathcal{B}}_{\ell} \setminus \mathcal{T}}_i$, $i = 1, \ldots, n_{\hat{\mathcal{B}}_{\ell} \setminus \mathcal{T}}$, are spline cuboids as defined in \cref{eq:spline_cuboid}, but for $\hat{\mathcal{B}}_{\ell}$. \\

For the first of our two approaches, we only need a partition of $\hat{\mathcal{B}}_{\ell} \cap \mathcal{T}$ into cuboids for each $\ell = 0, \ldots, L-1$, and each block $\mathbf{M}_{\ell, \ell'}$ of $\mathbf{M}$ in \cref{eq:hierarchical_mass_matrix_block_structure} again has a block structure of $n_{\hat{\mathcal{B}}_{\ell} \cap \mathcal{T}} \times n_{\hat{\mathcal{B}}_{\ell'} \cap \mathcal{T}}$ many blocks. To simplify notation, let $n_\ell = n_{\hat{\mathcal{B}}_{\ell} \cap \mathcal{T}}$ and $n_{\ell'} = n_{\hat{\mathcal{B}}_{\ell'} \cap \mathcal{T}}$, then we have
\begin{equation}
	\label{eq:approach_1_blockstructure}
	\mathbf{M}_{\ell, \ell'} = 
	\begin{bmatrix}
		\mathbf{M}_{\ell, \ell',1,1} & \cdots & \mathbf{M}_{\ell, \ell',1,n_{\ell'}} \\
		\vdots & & \vdots \\
		\mathbf{M}_{\ell, \ell',n_{\ell},1} & \cdots & \mathbf{M}_{\ell, \ell',n_{\ell},n_{\ell'}}
	\end{bmatrix}, 
\end{equation}
where each block $\mathbf{M}_{\ell, \ell', i, j}$ is of size $\left( n^{(1)}_{i}, n^{(2)}_{i}, n^{(3)}_{i} \right) \times \left( n^{(1)}_{j}, n^{(2)}_{j}, n^{(3)}_{j} \right)$. Here, $\left( n^{(1)}_{i}, n^{(2)}_{i}, n^{(3)}_{i} \right)$ denotes the size of the cuboid $\boldsymbol{\beta}^{\hat{\mathcal{B}}_{\ell} \cap \mathcal{T}}_i \subset \hat{\mathcal{B}}_{\ell} \cap \mathcal{T}$ and $\left( n^{(1)}_{j}, n^{(2)}_{j}, n^{(3)}_{j} \right)$ the size of the cuboid $\boldsymbol{\beta}^{\hat{\mathcal{B}}_{\ell'} \cap \mathcal{T}}_j \subset \hat{\mathcal{B}}_{\ell'} \cap \mathcal{T}$ for $\ell, \ell' = 0, \ldots, L-1$. The resulting global hierarchical mass tensor $\mathbf{M}$ \cref{eq:hierarchical_mass_matrix_block_structure} has then $\sum_{\ell = 0}^{L-1} n_{\hat{\mathcal{B}}_{\ell} \cap \mathcal{T}} \times \sum_{\ell = 0}^{L-1} n_{\hat{\mathcal{B}}_{\ell} \cap \mathcal{T}}$ many blocks. 

To extract the subblocks $\mathbf{M}_{\ell, \ell', i, j}$ in \cref{eq:approach_1_blockstructure}, we use Boolean selection tensors that pick the corresponding entries of $\tilde{\mathbf{M}}_{\ell, \ell'}$, i.e.
\begin{equation}
	\label{eq:approach_1_blocks}
	\mathbf{M}_{\ell, \ell', i, j} = \mathbf{J}^\top_{\boldsymbol{\beta}^{\hat{\mathcal{B}}_{\ell}\cap\mathcal{T}}_{i}} \tilde{\mathbf{M}}_{\ell, \ell'} \mathbf{J}_{\boldsymbol{\beta}^{\hat{\mathcal{B}}_{\ell'}\cap\mathcal{T}}_{j}} \in \mathbb{R}^{\left( n^{(1)}_{i}, n^{(2)}_{i}, n^{(3)}_{i} \right) \times \left( n^{(1)}_{j}, n^{(2)}_{j}, n^{(3)}_{j} \right)},
\end{equation}
where $\mathbf{J}^\top_{\boldsymbol{\beta}^{\hat{\mathcal{B}}_{\ell}\cap\mathcal{T}}_{i}}$ selects the DoFs of $\tilde{\mathcal{B}}_{\ell}$ associated with the cuboid $\boldsymbol{\beta}^{\hat{\mathcal{B}}_{\ell}\cap\mathcal{T}}_{i} \subset\hat{\mathcal{B}}_{\ell}\cap\mathcal{T}$ and $\mathbf{J}^\top_{\boldsymbol{\beta}^{\hat{\mathcal{B}}_{\ell'}\cap\mathcal{T}}_{j}}$ does the same for $\boldsymbol{\beta}^{\hat{\mathcal{B}}_{\ell'}\cap\mathcal{T}}_{j}\subset\hat{\mathcal{B}}_{\ell'}\cap\mathcal{T}$. This means the tensor $\mathbf{J}^\top_{\boldsymbol{\beta}^{\hat{\mathcal{B}}_{\ell}\cap\mathcal{T}}_{i}}$ assigns the DoFs of the reduced basis $\tilde{\mathcal{B}}_{\ell}$ to the DoFs of the splines contained in the single cuboid $\boldsymbol{\beta}^{\hat{\mathcal{B}}_{\ell} \cap \mathcal{T}}_i$ of the reduced basis $\hat{\mathcal{B}}_{\ell}$.

For a spline cuboid 
\begin{equation}
	\label{eq:spline_cuboid_global}
	\begin{aligned}
		\boldsymbol{\beta} & = \left\{ \hat{\beta}^{(1)}_{\ell, i^{(1)}_1}, \ldots, \hat{\beta}^{(1)}_{\ell, i^{(1)}_{n^{(1)}_{\boldsymbol{\beta}}}} \right\} \otimes \left\{ \hat{\beta}^{(2)}_{\ell, i^{(2)}_1}, \ldots, \hat{\beta}^{(2)}_{\ell, i^{(2)}_{n^{(2)}_{\boldsymbol{\beta}}}} \right\} \otimes \left\{ \hat{\beta}^{(3)}_{\ell, i^{(3)}_1}, \ldots, \hat{\beta}^{(3)}_{\ell, i^{(3)}_{n^{(3)}_{\boldsymbol{\beta}}}} \right\} \\
		& = \boldsymbol{\beta}^{(1)} \otimes \boldsymbol{\beta}^{(2)} \otimes \boldsymbol{\beta}^{(3)} \subset \hat{\mathcal{B}}_{\ell}
	\end{aligned}
\end{equation}
such a selecting Boolean tensor in \cref{eq:approach_1_blocks} factors as a Kronecker product of one-dimensional selectors:
\begin{equation}
	\label{eq:approach_1_Boolean}
	\begin{gathered}
		\mathbf{J}_{\boldsymbol{\beta}} = J^{(1)}_{\boldsymbol{\beta}^{(1)}} \otimes J^{(2)}_{\boldsymbol{\beta}^{(2)}} \otimes J^{(2)}_{\boldsymbol{\beta}^{(1)}} \in \{0,1\}^{\left( \tilde{n}^{(1)}_{\ell}, \tilde{n}^{(2)}_{\ell}, \tilde{n}^{(3)}_{\ell} \right) \times \left( n^{(1)}_{\boldsymbol{\beta}}, n^{(2)}_{\boldsymbol{\beta}}, n^{(3)}_{\boldsymbol{\beta}} \right)}, \\
		\left( J^{(d)}_{\boldsymbol{\beta}^{(d)}} \right)_{i,j} = 
		\begin{cases}
			1 &: \text{ if } \tilde{\beta}^{(d)}_{\ell, i} = \hat{\beta}^{(d)}_{\ell, i^{(d)}_j} \in \boldsymbol{\beta}^{(d)},  \\
			0 &: \text{ otherwise.}
		\end{cases}
	\end{gathered}
\end{equation}
This realizes a restriction from the reduced basis $\tilde{\mathcal{B}}_{\ell}$ to the local DoFs of the cuboid $\boldsymbol{\beta}$ in $\hat{\mathcal{B}}_{\ell}$.\\

For the second approach, we need a partition of $\hat{\mathcal{B}}_{\ell} \cap \mathcal{T}$ as well as of $\hat{\mathcal{B}}_{\ell} \setminus \mathcal{T}$ into spline cuboids. Here the individual block $\mathbf{M}_{\ell, \ell'}$ in \cref{eq:hierarchical_mass_matrix_block_structure} is a single tensor which represents all inner products of the basis functions between the bases $\hat{\mathcal{B}}_{\ell}$ and $\hat{\mathcal{B}}_{\ell'}$, i.e.
\begin{equation*}
	\mathbf{M}_{\ell, \ell'} = \mathbf{J}^\top_{\ell} \tilde{\mathbf{M}}_{\ell, \ell'} \mathbf{J}_{\ell'} \in \mathbb{R}^{\left( \hat{n}^{(1)}_{\ell}, \hat{n}^{(2)}_{\ell}, \hat{n}^{(3)}_{\ell} \right) \times \left( \hat{n}^{(1)}_{\ell'}, \hat{n}^{(2)}_{\ell'}, \hat{n}^{(3)}_{\ell'} \right)}
\end{equation*}
for $\ell, \ell' = 0, \ldots, L-1$. This means that, unlike in the first approach, $\mathbf{M}$ retains the block structure shown in \cref{eq:hierarchical_mass_matrix_block_structure} with $L \times L$ blocks. In this approach the non-diagonal blocks in \cref{eq:hierarchical_mass_matrix_block_structure} are defined by
\begin{equation}
	\label{eq:approach_2_non_diagonal_blocks}
	\mathbf{M}_{\ell, \ell'} = \sum_{i = 1}^{n_{\hat{\mathcal{B}}_{\ell} \cap \mathcal{T}}} \sum_{j = 1}^{n_{\hat{\mathcal{B}}_{\ell'} \cap \mathcal{T}}} \mathbf{J}^\top_{\boldsymbol{\beta}^{\hat{\mathcal{B}}_{\ell}\cap\mathcal{T}}_{i}} \tilde{\mathbf{M}}_{\ell, \ell'} \mathbf{J}_{\boldsymbol{\beta}^{\hat{\mathcal{B}}_{\ell'}\cap\mathcal{T}}_{j}}, \quad \ell \neq \ell', 
\end{equation}
where $\mathbf{J}^\top_{\boldsymbol{\beta}^{\hat{\mathcal{B}}_{\ell}\cap\mathcal{T}}_{i}}$ selects the DoFs of $\tilde{\mathcal{B}}_{\ell}$ associated with the cuboid $\boldsymbol{\beta}^{\hat{\mathcal{B}}_{\ell}\cap\mathcal{T}}_{i} \subset\hat{\mathcal{B}}_{\ell}\cap\mathcal{T}$ and $\mathbf{J}^\top_{\boldsymbol{\beta}^{\hat{\mathcal{B}}_{\ell'}\cap\mathcal{T}}_{j}}$ does the same for $\boldsymbol{\beta}^{\hat{\mathcal{B}}_{\ell'}\cap\mathcal{T}}_{j}\subset\hat{\mathcal{B}}_{\ell'}\cap\mathcal{T}$. In contrast to \cref{eq:approach_1_blocks}, the selector $\mathbf{J}^{\top}_{\boldsymbol{\beta}^{\hat{\mathcal{B}}_{\ell}\cap\mathcal{T}}_{i}}$ here embeds the cuboid $\boldsymbol{\beta}^{\hat{\mathcal{B}}_{\ell}\cap\mathcal{T}}_{i}$ into the \emph{global} reduced basis $\hat{\mathcal{B}}_{\ell}$, i.e., it maps into $\hat{\mathcal{B}}_{\ell}$. For a spline cuboid $\boldsymbol{\beta}=\boldsymbol{\beta}^{(1)}\otimes\boldsymbol{\beta}^{(2)}\otimes\boldsymbol{\beta}^{(3)}\subset\hat{\mathcal{B}}_{\ell}$ as in \cref{eq:spline_cuboid_global}, the Boolean selector in \cref{eq:approach_2_non_diagonal_blocks} factors as
\begin{gather*}
	\mathbf{J}_{\boldsymbol{\beta}} = J^{(1)}_{\boldsymbol{\beta}^{(1)}} \otimes J^{(2)}_{\boldsymbol{\beta}^{(2)}} \otimes J^{(2)}_{\boldsymbol{\beta}^{(1)}} \in \{0,1\}^{\left( \tilde{n}^{(1)}_{\ell}, \tilde{n}^{(2)}_{\ell}, \tilde{n}^{(3)}_{\ell} \right) \times \left( \hat{n}^{(1)}_{\ell}, \hat{n}^{(2)}_{\ell}, \hat{n}^{(3)}_{\ell} \right)}, \\
	\left( J^{(d)}_{\boldsymbol{\beta}^{(d)}} \right)_{i,j} = 
	\begin{cases}
		1 &: \text{ if } \tilde{\beta}^{(d)}_{\ell, i} = \hat{\beta}^{(d)}_{\ell, j} \in \boldsymbol{\beta}^{(d)},  \\
		0 &: \text{ otherwise,}
	\end{cases}
\end{gather*}
The difference to \cref{eq:approach_1_Boolean} is thus the \emph{codomain}: here the selector spans the full basis $\hat{\mathcal{B}}_{\ell}$ (size $\hat n^{(1)}_{\ell}\times \hat n^{(2)}_{\ell}\times \hat n^{(3)}_{\ell}$), rather than the local cuboid of size $n_{\boldsymbol{\beta}}^{(1)}\times n_{\boldsymbol{\beta}}^{(2)}\times n_{\boldsymbol{\beta}}^{(3)}$. 

Since $\hat{\mathcal{B}}_{\ell}$ generally also contains non-active basis functions, using \cref{eq:approach_2_non_diagonal_blocks} with $\ell=\ell'$ would produce \emph{singular} diagonal blocks in \cref{eq:hierarchical_mass_matrix_block_structure}: all entries of $\mathbf{M}_{\ell,\ell}$ involving basis functions from $\hat{\mathcal{B}}_{\ell}\setminus\mathcal{T}$ are zero. Consequently $\mathbf{M}$ would be singular. To avoid this, we set the diagonal entries corresponding to non-active basis functions to one. Using the partition of $\hat{\mathcal{B}}_{\ell}\setminus\mathcal{T}$ into cuboids from \cref{eq:spline_partition}, we add Boolean diagonal masks on those indices
\begin{equation}
	\label{eq:approach_2_diagonal_blocks}
	\mathbf{M}_{\ell, \ell} = \sum_{i, j = 1}^{n_{\hat{\mathcal{B}}_{\ell} \cap \mathcal{T}}} \mathbf{J}^\top_{\boldsymbol{\beta}^{\hat{\mathcal{B}}_{\ell}\cap\mathcal{T}}_{i}} \tilde{\mathbf{M}}_{\ell, \ell} \mathbf{J}_{\boldsymbol{\beta}^{\hat{\mathcal{B}}_{\ell}\cap\mathcal{T}}_{j}} + \sum_{k = 1}^{n_{\hat{\mathcal{B}}_{\ell} \setminus \mathcal{T}}} \mathbf{I}_{\boldsymbol{\beta}^{\hat{\mathcal{B}}_{\ell} \setminus \mathcal{T}}_k}.
\end{equation}
For a spline cuboid $\boldsymbol{\beta}=\boldsymbol{\beta}^{(1)}\otimes\boldsymbol{\beta}^{(2)}\otimes\boldsymbol{\beta}^{(3)}\subset\hat{\mathcal{B}}_{\ell}$ (cf.~\cref{eq:spline_cuboid_global}), the mask factors as
\begin{gather*}
	\mathbf{I}_{\boldsymbol{\beta}} = I^{(1)}_{\boldsymbol{\beta}^{(1)}} \otimes I^{(2)}_{\boldsymbol{\beta}^{(2)}} \otimes I^{(2)}_{\boldsymbol{\beta}^{(1)}} \in \{0,1\}^{\left( \hat{n}^{(1)}_{\ell}, \hat{n}^{(2)}_{\ell}, \hat{n}^{(3)}_{\ell} \right) \times \left( \hat{n}^{(1)}_{\ell}, \hat{n}^{(2)}_{\ell}, \hat{n}^{(3)}_{\ell} \right)}, \\
	\left( I^{(d)}_{\boldsymbol{\beta}^{(d)}} \right)_{i,i} = 
	\begin{cases}
		1 &: \text{ if } \hat{\beta}^{(d)}_{\ell, i} \in \boldsymbol{\beta}^{(d)},  \\
		0 &: \text{ otherwise.}
	\end{cases}
\end{gather*}
This leaves the operator on the active subspace unchanged while preventing singular diagonal blocks.

\section{Hierarchical Low-Rank solver} 
\label{section:Hierarchical_Low-Rank_solver}
Since the linear system \cref{eq:linear_system}, constructed using the heuristics described in \cref{section:assembling_linear_system}, has a block structure as in \cref{eq:hierarchical_mass_matrix_block_structure} for both approaches \cref{eq:approach_1_blockstructure} and \cref{eq:approach_2_non_diagonal_blocks} \& \cref{eq:approach_2_diagonal_blocks}, it is reasonable to use the block version of \ttgmres~(implemented as \texttt{tt\_gmres\_block.m} in the \textsc{TT-Toolbox}~\cite{tt-toolbox}, cf.~\cite{Dolgov:2013}). This method is tensor-matrix free, i.e., we can define the linear operator $\mathbf{K}$ as a function handle. However, choosing the efficient block version of \amen~(implemented as \texttt{amen\_block\_solve.m} in the \textsc{TT-Toolbox}~\cite{tt-toolbox}) introduces a complication: all blocks in $\mathbf{K}$ must have the same size for this solver. One option is to pad the smaller blocks with zeros in the corresponding entries, but this significantly impairs numerical stability, and the performance of the method is not competitive. Furthermore, although the system is symmetric and positive definite, a variant of the (preconditioned) conjugate gradient method -- typically the standard choice for such systems -- is not available in~\cite{tt-toolbox}.
The linear system \cref{eq:linear_system} constructed by both approaches \cref{eq:approach_1_blockstructure} and \cref{eq:approach_2_non_diagonal_blocks} \& \cref{eq:approach_2_diagonal_blocks} is extremely structured, which is why computing the solution without a preconditioner would be extremely expensive and time-consuming. \\

For the approach \cref{eq:approach_2_non_diagonal_blocks} \& \cref{eq:approach_2_diagonal_blocks}, $\mathbf{K}$ has an $L \times L$ block structure. We therefore use a block-diagonal left preconditioner $\mathbf{P}$ of the same size, whose diagonal blocks are given by
\begin{equation}
	\label{eq:block_preconditioner_2}
	\mathbf{P}_{\ell, \ell} = \mathbf{K}_{\ell, \ell}, \quad \ell = 0, \ldots, L-1.
\end{equation}
When applying ${\mathbf{P}_{\ell, \ell}}^{-1}$ within the function handle, the level-wise linear systems are solved using the standard the standard \amen~method (cf.~\texttt{amen\_solve2.m} in~\cite{tt-toolbox}). The method iteratively updates the TT cores in \cref{eq:TT_format} by minimizing an energy function associated with the linear system that represents the residual error of the current approximation, while maintaining its low-rank structure. The algorithm expands the search space using an inexact gradient direction at each step and can achieve geometric convergence. The computational complexity is linear with respect to the dimension and the size of the tensor (cf.~\cite{amen} for a more detailed derivation). 

For the approach in \cref{eq:approach_1_blockstructure}, \cref{eq:block_preconditioner_2} is applicable; however, applying ${\mathbf{P}_{\ell,\ell}}^{-1}$ would require a solver for block-structured systems, which is numerically costly in many cases, since each block $\mathbf{K}_{\ell,\ell}$ again has an $n_{\hat{\mathcal{B}}_{\ell} \cap \mathcal{T}} \times n_{\hat{\mathcal{B}}_{\ell'} \cap \mathcal{T}}$ block structure. To use the efficient \amen~method, we choose as preconditioner the diagonal blocks of the diagonal blocks of $\mathbf{K}$; i.e., $\mathbf{P}$ has the same $\sum_{\ell = 0}^{L-1} n_{\hat{\mathcal{B}}_{\ell} \cap \mathcal{T}} \times \sum_{\ell = 0}^{L-1} n_{\hat{\mathcal{B}}_{\ell} \cap \mathcal{T}}$ block structure as $\mathbf{K}$, and its diagonal blocks are given by
\begin{equation}
	\label{eq:block_preconditioner_1}
	\mathbf{P}_{\ell,\ell,i,i} = \mathbf{K}_{\ell,\ell,i,i}, \quad i = 1, \ldots, n_{\hat{\mathcal{B}}_{\ell} \cap \mathcal{T}}, \; \ell = 0, \ldots, L-1.
\end{equation}

Additionally, we use a standard Jacobi preconditioner for both approaches, \cref{eq:approach_1_blockstructure} and \cref{eq:approach_2_non_diagonal_blocks} \& \cref{eq:approach_2_diagonal_blocks}. For \cref{eq:approach_2_non_diagonal_blocks} \& \cref{eq:approach_2_diagonal_blocks}, this preconditioner has an $L \times L$ block structure, and its diagonal blocks correspond to the diagonal entries of \cref{eq:block_preconditioner_2}, which can be recovered from $\mathbf{K}_{\ell, \ell}$ using the diagonal entries of the TT cores $\mathbf{K}^{(d)}_{\ell, \ell, \left(r_{d-1}, :, :, r_d  \right)} \in \mathbf{R}^{\hat{n}^{(d)}_{\ell} \times \hat{n}^{(d)}_{\ell}}$. Similarly, for \cref{eq:approach_1_blockstructure}, the $\sum_{\ell = 0}^{L-1} n_{\hat{\mathcal{B}}_{\ell} \cap \mathcal{T}}$ diagonal blocks of the standard Jacobi preconditioner can be constructed from the diagonal elements of the TT cores of $\mathbf{K}_{\ell,\ell,i,i}$,  $i = 1, \ldots, n_{\hat{\mathcal{B}}_{\ell} \cap \mathcal{T}}$, $\ell = 0, \ldots, L-1$. For both approaches, the Jacobi preconditioner can be applied block-wise using the efficient \amen~solver. Its straightforward implementation encourages its use.

Alternatively, one could use a BPX preconditioner, as presented in~\cite{BRACCO2021113742} and considered state of the art~\cite{Buffa2022AdaptiveIGA}, but its implementation in low-rank form is not yet clear.

\section{Numerical experiments} 
\label{section:numerics}
In this section, we present the results of our numerical experiments. The experiments were implemented in \textsc{MATLAB}~\cite{MATLAB_R2024b} using the \textsc{NURBS} package~\cite{nurbs_toolbox} version 1.4.3, the \textsc{GeoPDEs} package~\cite{geopdes3.0, geopdes_new} and hierarchical package version 3.2.2 and the \textsc{TT-Toolbox}~\cite{tt-toolbox} version 2.2.2. We investigate the extent to which the method presented in this work is numerically efficient and memory intensive in comparison with methods from the \textsc{GeoPDEs} package, which fully assemble the system. The experiments were performed in \textsc{MATLAB} R2024b on a Lenovo ThinkPad L14 Gen 4 with a 13th Gen Intel Core i7-1355U x 12, 32 GiB RAM and Ubuntu 24.04.2LTS operating system.\footnote{The code is public at \url{https://github.com/TomRiem/Low_Rank_THB_splines.git}.}

Our experiments consist of approximating the solution of the Poisson equation with homogenouos Dirichlet boundary conditions \cref{eq:forward} on the unit cube $\Omega = \left[0 , 1\right]^3$. We use different source functions $f$ and apply hierarchical refinement in subdomains of $\Omega$ where we expect significant activity of the respective solution $y$. We construct an approximate linear system \cref{eq:linear_system} using the heuristics presented in \cref{section:assembling_linear_system} to assemble the hierarchical stiffness tensor $\mathbf{K}$ and the hierarchical right-hand side $\mathbf{f}$ in TT format, and solve it using the block version of \ttgmres~(implemented as \texttt{tt\_gmres\_block.m} in~\cite{tt-toolbox}, cf.~\cite{Dolgov:2013}) to obtain a coefficient tensor $\mathbf{y}$ that represents the degrees of freedom (DoFs) of the B-spline basis.

\begin{figure}[htbp]
\centering
\begin{subfigure}[t]{0.48\textwidth}
  \centering
  \includegraphics[width=\linewidth]{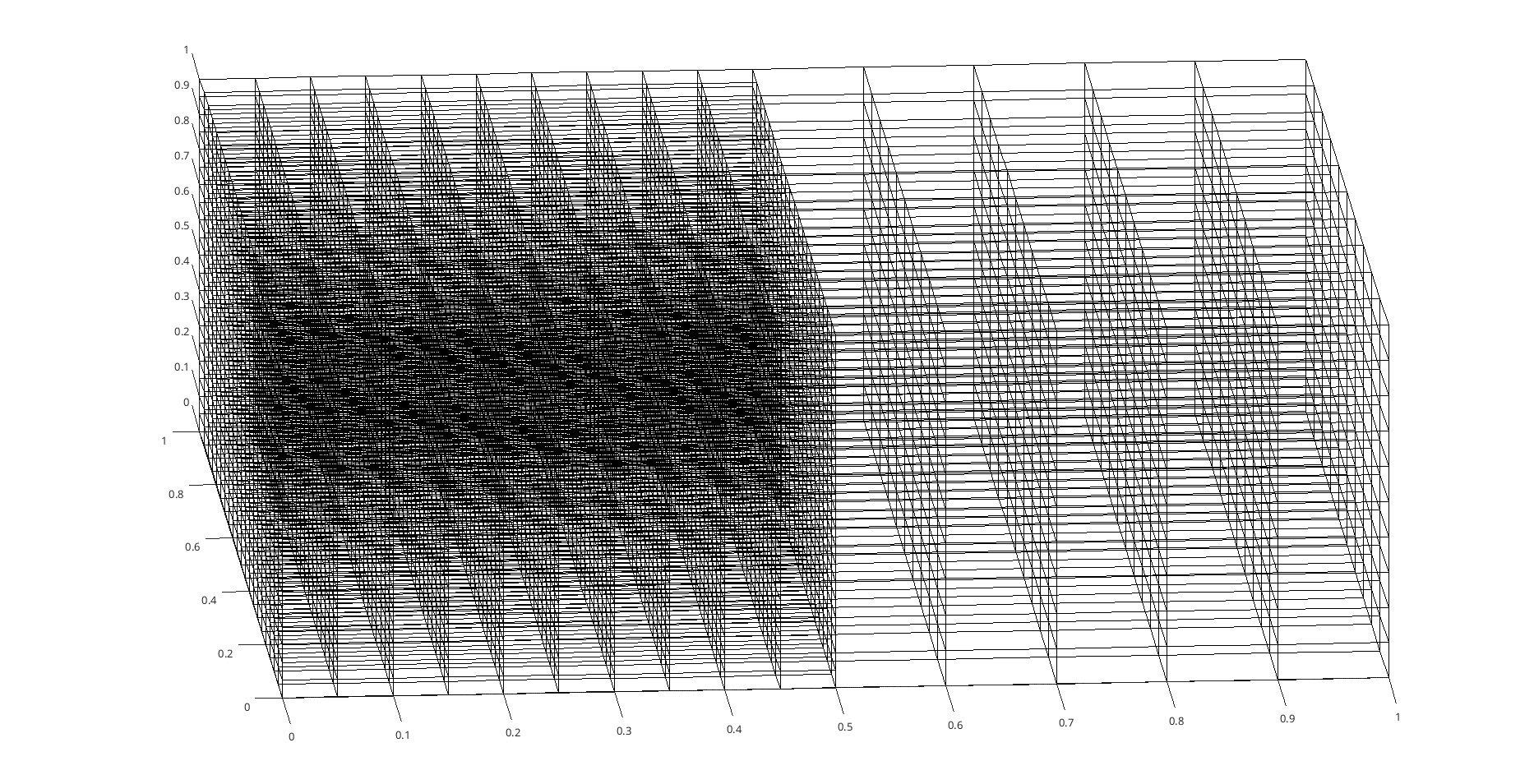}
  \caption{\cref{eq:ref_1} for $k = 2$.}
  \label{fig:panel-a}
\end{subfigure}\hfill
\begin{subfigure}[t]{0.48\textwidth}
  \centering
  \includegraphics[width=\linewidth]{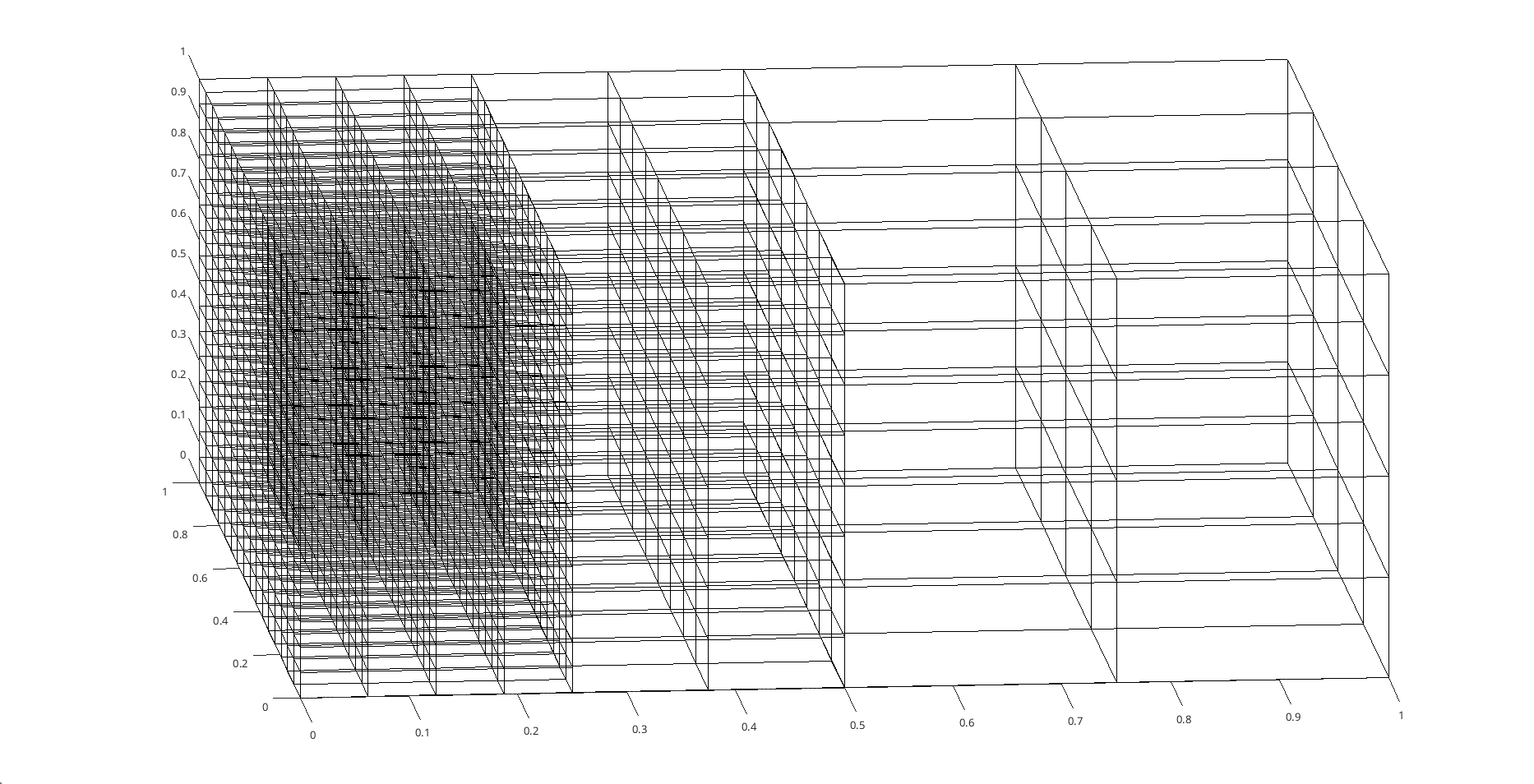}
  \caption{\cref{eq:ref_2} for $L = 3$.}
  \label{fig:panel-b}
\end{subfigure}
\medskip
\begin{subfigure}[t]{0.48\textwidth}
  \centering
  \includegraphics[width=\linewidth]{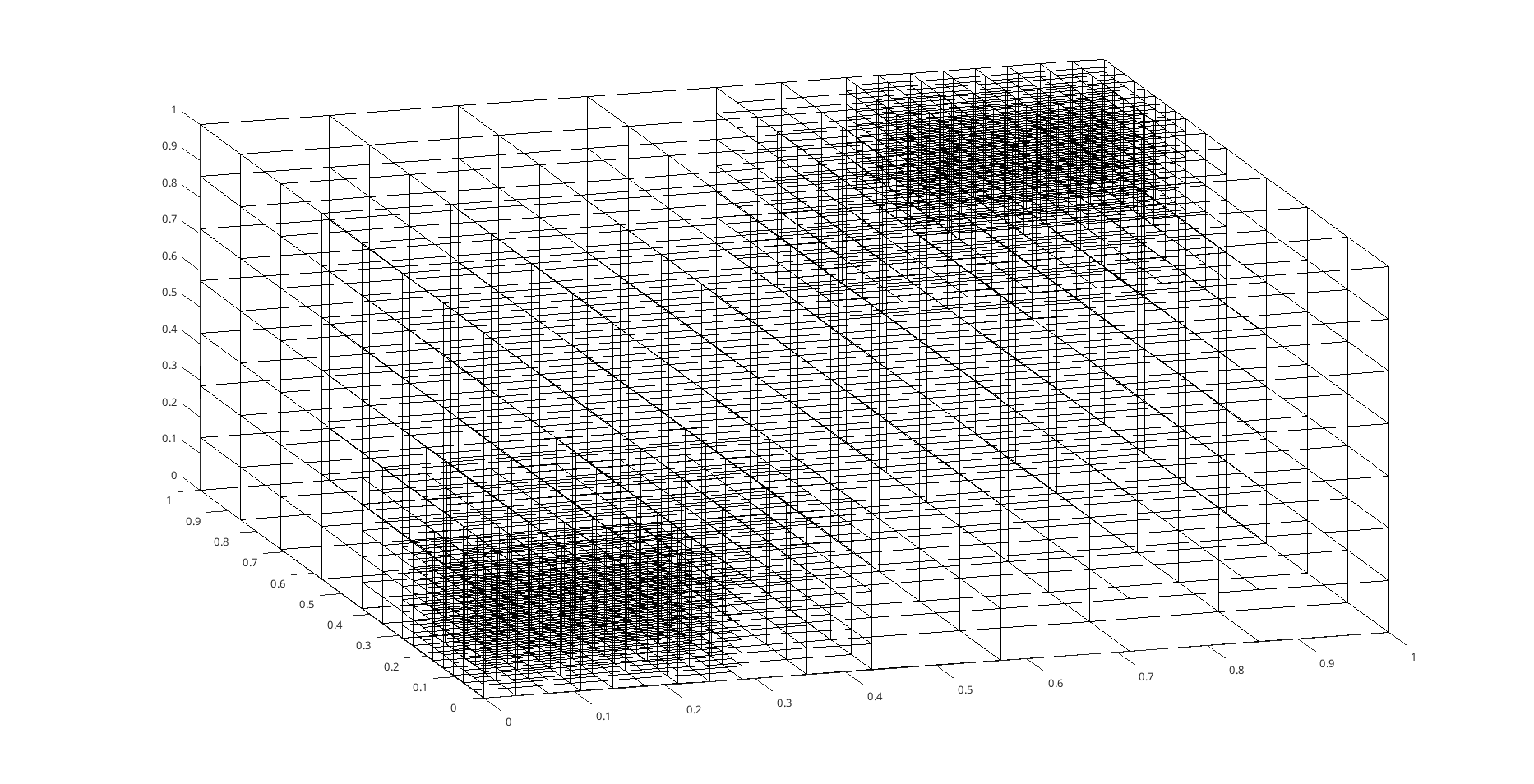}
  \caption{\cref{eq:ref_3} for $L = 3$.}
  \label{fig:panel-c}
\end{subfigure}\hfill
\begin{subfigure}[t]{0.48\textwidth}
  \centering
  \includegraphics[width=\linewidth]{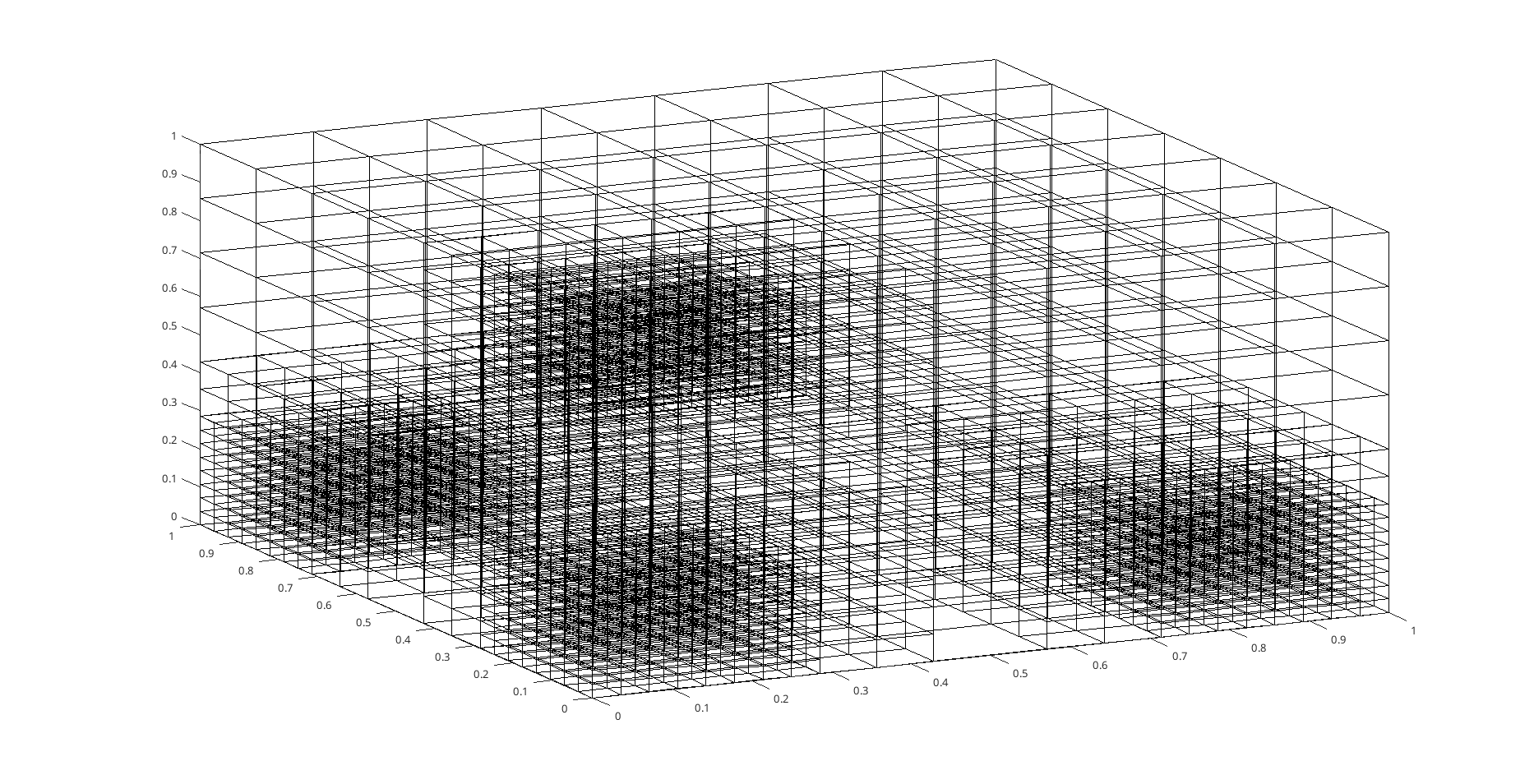}
  \caption{\cref{eq:ref_4} for $L = 3$.}
  \label{fig:panel-d}
\end{subfigure}
\caption{Visualization in \textsc{MATLAB} of the cells of the used hierarchical meshes using \texttt{hmsh\_plot\_cells.m} for the case $p = 3$.}
\label{fig:four-up}
\end{figure}

For representing the computational domain $\Omega = \left[ 0, 1 \right]^3$ we use a B-spline basis with $p^{(d)} = 1$ and $n^{(d)} = 2$, $d = 1,2,3$. The resulting B-spline basis for interpolating the weight function $Q$ \cref{eq:weight_functions} then has degrees $p^{(d)} = 2$ and $n^{(d)} = 4$, $d = 1,2,3$. In the following, we approximate the solution of \cref{eq:forward} using B-spline bases with $p^{(d)} = 3$ and $p^{(d)} = 5$; accordingly, we use a B-spline basis with the same degrees in each case to interpolate the source function $f$. In the case of $p^{(d)} = 3$, we interpolate $f$ with $n^{(d)} = 159$ splines, and in the case of $p^{(d)} = 5$, we interpolate with $n^{(d)} = 171$ splines, $d = 1,2,3$.

For each numerical experiment, we measure the computational time in seconds $s$ to compute an approximate solution $\mathbf{y}$, the $L^2$ error to the analytical solution $y \colon \Omega \to \mathbb{R}$ of the respective elliptic problem \cref{eq:forward}, the estimated memory footprint in number of bytes $n_{\text{bytes}}$ of the linear system $\mathbf{K}$ and the solution $\mathbf{y}$ and the number of iterations $n_{\text{it}}$ the block version of \ttgmres~needs fo convergence. To generate the data structures for hierarchical IgA, we rely on functions from the \textsc{GeoPDEs} package. After generating these structures, we start measuring the computational time of our method, that is, from the interpolation of the weight function until the linear system \cref{eq:linear_system} has been solved, i.e., once \texttt{tt\_gmres\_block.m} has converged. We compare the computational time of our method with the time required using \textsc{GeoPDEs}, starting from the call of \texttt{op\_gradu\_gradv\_hier.m} to set up a hierarchical stiffness tensor and ending after solving the resulting linear system with the usual backslash operator. We estimate memory usage with a generic recursive size estimator that walks through each data structure (e.g., matrices and tensors in TT format), summing the bytes of all components it can reach. To avoid double counting, shared references are counted once. The estimate is approximate because low-level runtime headers and copy-on-write effects are not visible to user-level code. We compare the memory usage to that of the hierarchical stiffness tensor and the full solution vector from \textsc{GeoPDEs}. For the error calculation, we transform our low-rank solution $\mathbf{y}$ into a full coefficient vector, i.e., $\operatorname{vec}(\mathbf{y}) \in \mathbb{R}^{N}$, where $N$ is the total number of degrees of freedom, and then use the function \texttt{sp\_l2\_error.m} from the \textsc{GeoPDEs} package to compute the $L^2$ error with respect to the respective analytical solution $y \colon \Omega \to \mathbb{R}$. We compare the error to that of the solution resulting from the full solution vector from \textsc{GeoPDEs}. 

In the following, $\varepsilon$ denotes the tolerance for solving the linear system \cref{eq:linear_system} with \texttt{tt\_gmres\_block.m}, with a maximum of $900$ iterations and restarts every $30$ iterations. We interpolate the weight function $Q$ and the source function $f$ to tolerance $\varepsilon \cdot 10^{-2}$ using \texttt{amen\_block\_solve.m}\footnote{with parameters \texttt{kickrank} $= 2$, \texttt{resid\_damp} $= 10$, \texttt{nswp} $= 10$, and \texttt{exitdir} $= -1$}. To apply the preconditioner, we use \texttt{amen\_solve2.m}\footnote{with parameters \texttt{nswp} $= 20$ and \texttt{kickrank} $= 2$} with tolerance $\varepsilon \cdot 10^{-2}$.

Since the solution of \cref{eq:forward} must satisfy the homogeneous Dirichlet conditions, the function  
\begin{equation}
	\label{eq:sol_0}
	y_0 \left( x^{(1)}, x^{(2)}, x^{(3)} \right) = x^{(1)} \left( x^{(1)} - 1 \right) x^{(2)} \left( x^{(2)} - 1 \right) x^{(3)} \left( x^{(3)} - 1 \right).
\end{equation}
will be a factor in every analytical solution used here. We multiply \cref{eq:sol_0} by a term that increases the activity in a subdomain of $\Omega$ and then apply hierarchical refinement in this subdomain.

The analytical solution of our first experiment is
\begin{equation}
	\label{eq:sol_1}
	y_1 \left( x^{(1)}, x^{(2)}, x^{(3)} \right) = y_0 \left( x^{(1)}, x^{(2)}, x^{(3)} \right) \operatorname{exp}\left( -{x^{(1)}}^2 \right).
\end{equation}
We approximate \eqref{eq:sol_1} in a two-level THB-spline space $\mathcal{T}$, i.e., $L = 2$. For each degree $p$, the initial mesh has $m^{(d)}_0 = 6 + 2k$ cells per parametric direction $d=1,2,3$, where $k$ is increasing. The domain $\Omega_1$ can be represented in $G_0$ by the cells
\begin{equation}
	\label{eq:ref_1}
	\mathcal{M} = \left\{ q_{0,\left( j^{(1)},j^{(2)},j^{(3)} \right)} \in G_0:\; j^{(1)} = 1, \ldots, \frac{6 + 2k}{2} \right\},
\end{equation}
which means we apply single dyadic subdivision (factor $2$ in each direction) of the cells contained in $\mathcal{M}$. The refined subdomain is then the left half $\Omega_1 = \left[0, \frac{1}{2} \right] \times \left[ 0, 1 \right] \times \left[ 0, 1 \right]$. We employ the associated two-level truncated hierarchical basis. We repeat this for $k=0,\ldots,4$ when $p=3$ and for $k=0,1,2$ when $p=5$. A visualization of the hierarchical mesh for the case $p = 3$ and $k = 2$ is shown in \cref{fig:panel-a}. 

We note that for such hierarchical refinement, the two approaches \cref{eq:approach_1_blockstructure} and \cref{eq:approach_2_non_diagonal_blocks} \& \cref{eq:approach_2_diagonal_blocks} are equivalent, since there is only one spline cuboid for each level $\ell = 0, 1$, and thus $\mathbf{M}_{\ell, \ell'}$ consists of just one block, that has the same size $\left( \hat{n}^{(1)}_{\ell}, \hat{n}^{(2)}_{\ell}, \hat{n}^{(3)}_{\ell} \right) \times \left( \hat{n}^{(1)}_{\ell'}, \hat{n}^{(2)}_{\ell'}, \hat{n}^{(3)}_{\ell'} \right)$.

\begin{figure}
\definecolor{mycolor1}{rgb}{0.00000,0.44700,0.74100}%
\definecolor{mycolor2}{rgb}{0.85000,0.32500,0.09800}%
\definecolor{mycolor3}{rgb}{0.92900,0.69400,0.12500}%
\definecolor{mycolor4}{rgb}{0.49400,0.18400,0.55600}%
\definecolor{mycolor5}{rgb}{0.46600,0.67400,0.18800}%
\begin{subfigure}[t]{0.49\textwidth}
\centering
\begin{tikzpicture}[font=\footnotesize]
\begin{axis}[%
width=1.0\textwidth,
height=0.6\textwidth,
xmin=0,
xmax=6,
xlabel={$k$},
xtick = {1, 2, 3, 4, 5},
xticklabels={0, 1, 2, 3, 4},
ymode=log,
ymin=0,
ymax=500,
ylabel={$s$},
yminorticks=true,
ymajorgrids=true,
ytick = {1, 10, 100}
]
\addplot [color=mycolor5, mark=star, mark size=1.0pt, mark options={solid, mycolor5}]
table[row sep=crcr]{%
1	6.437089\\
2	15.748442\\
3	30.567591\\
4	53.393953\\
5	85.825123\\
};
\addplot [color=mycolor1, mark=square, mark size=1.0pt, mark options={solid, mycolor1}]
table[row sep=crcr]{%
1	0.546032\\
2	0.237786\\
3	0.313628\\
4	0.38465\\
5	0.444305\\
};
\addplot [color=mycolor2, mark=square, mark size=1.0pt, mark options={solid, mycolor2}]
table[row sep=crcr]{%
1	0.232127\\
2	0.330818\\
3	0.358429\\
4	0.510788\\
5	0.633511\\
};
\addplot [color=mycolor3, mark=square, mark size=1.0pt, mark options={solid, mycolor3}]
table[row sep=crcr]{%
1	0.285345\\
2	0.378133\\
3	0.511956\\
4	0.724546\\
5	0.988991\\
};
\addplot [color=mycolor1, mark=triangle, mark size=1.0pt, mark options={solid, mycolor1}]
table[row sep=crcr]{%
1	0.20548\\
2	0.298463\\
3	0.398014\\
4	0.584984\\
5	0.796022\\
};
\addplot [color=mycolor2, mark=triangle, mark size=1.0pt, mark options={solid, mycolor2}]
table[row sep=crcr]{%
1	1.18794\\
2	0.571289\\
3	0.671851\\
4	0.946882\\
5	1.266138\\
};
\addplot [color=mycolor3, mark=triangle, mark size=1.0pt, mark options={solid, mycolor3}]
table[row sep=crcr]{%
1	0.986322\\
2	1.230567\\
3	1.460663\\
4	1.972483\\
5	2.440537\\
};
\end{axis}
\end{tikzpicture}
\caption{Comp. time for $p = 3$.} \label{figure:1_ex_3_1}
\end{subfigure}
\hfill
\begin{subfigure}[t]{0.49\textwidth}
\centering
\begin{tikzpicture}[font=\footnotesize]
\begin{axis}[%
width=1.0\textwidth,
height=0.6\textwidth,
xmin=0,
xmax=6,
xlabel={$k$},
xtick = {1, 2, 3, 4, 5},
xticklabels={0, 1, 2, 3, 4},
ymode=log,
ymin=1e-9,
ymax=1e-5,
ylabel={$\left\Vert \cdot \right\Vert_{L^2}$},
yminorticks=true,
ymajorgrids=true,
ytick = {1e-8, 1e-7, 1e-6}
]
\addplot [color=mycolor5, mark=star, mark size=1.0pt, mark options={solid, mycolor5}]
table[row sep=crcr]{%
1	2.29835480913388e-07\\
2	6.19409792414484e-08\\
3	2.49799312076262e-08\\
4	1.15391399534563e-08\\
5	6.1398898769322e-09\\
};
\addplot [color=mycolor1, mark=square, mark size=1.0pt, mark options={solid, mycolor1}]
table[row sep=crcr]{%
1	3.79629734424419e-07\\
2	5.04668021183064e-07\\
3	5.33515110104063e-07\\
4	9.18588974283876e-07\\
5	2.54996842131607e-06\\
};
\addplot [color=mycolor2, mark=square, mark size=1.0pt, mark options={solid, mycolor2}]
table[row sep=crcr]{%
1	2.29958256760937e-07\\
2	6.21535030430294e-08\\
3	2.55454963129188e-08\\
4	1.25518401679194e-08\\
5	1.23560967980938e-08\\
};
\addplot [color=mycolor3, mark=square, mark size=1.0pt, mark options={solid, mycolor3}]
table[row sep=crcr]{%
1	2.29852343920688e-07\\
2	6.19574619123901e-08\\
3	2.50083294335833e-08\\
4	1.15979280825337e-08\\
5	6.24528447758447e-09\\
};
\addplot [color=mycolor1, mark=triangle, mark size=1.0pt, mark options={solid, mycolor1}]
table[row sep=crcr]{%
1	1.0788783537371e-06\\
2	1.59933161201336e-06\\
3	1.27250868431548e-06\\
4	1.45254674228231e-06\\
5	3.02196527181912e-06\\
};
\addplot [color=mycolor2, mark=triangle, mark size=1.0pt, mark options={solid, mycolor2}]
table[row sep=crcr]{%
1	2.31186504152141e-07\\
2	6.25986281942675e-08\\
3	2.72146605399977e-08\\
4	1.29170094988667e-08\\
5	1.313469403026e-08\\
};
\addplot [color=mycolor3, mark=triangle, mark size=1.0pt, mark options={solid, mycolor3}]
table[row sep=crcr]{%
1	2.2985511452501e-07\\
2	6.19585912228353e-08\\
3	2.50110475530184e-08\\
4	1.16003072645573e-08\\
5	6.24858697338128e-09\\
};
\end{axis}
\end{tikzpicture}
\caption{$L^2$ error for $p = 3$.} \label{figure:1_ex_3_2}
\end{subfigure} \\
\begin{subfigure}[t]{0.49\textwidth}
\centering
\begin{tikzpicture}[font=\footnotesize]
\begin{axis}[%
width=1.0\textwidth,
height=0.6\textwidth,
xmin=0,
xmax=6,
xlabel={$k$},
xtick = {1, 2, 3, 4, 5},
xticklabels={0, 1, 2, 3, 4},
ymode=log,
ymin=5000,
ymax=5e8,
ylabel={$n_{\text{bytes}}$},
yminorticks=true,
ymajorgrids=true,
ytick = {1e4, 1e5, 1e6, 1e7, 1e8}
]
\addplot [color=mycolor5, mark=star, mark size=1.0pt, mark options={solid, mycolor5}]
table[row sep=crcr]{%
1	3700344\\
2	9725984\\
3	20127240\\
4	36091248\\
5	58805144\\
};
\addplot [color=mycolor1, mark=*, mark size=1.0pt, mark options={solid, mycolor1}]
table[row sep=crcr]{%
1	26212\\
2	40564\\
3	58660\\
4	80500\\
5	106084\\
};
\addplot [color=mycolor2, mark=*, mark size=1.0pt, mark options={solid, mycolor2}]
table[row sep=crcr]{%
1	26212\\
2	40564\\
3	58660\\
4	80500\\
5	106084\\
};
\addplot [color=mycolor3, mark=*, mark size=1.0pt, mark options={solid, mycolor3}]
table[row sep=crcr]{%
1	26212\\
2	40564\\
3	58660\\
4	80500\\
5	106084\\
};
\end{axis}
\end{tikzpicture}
\caption{Memory of $\mathbf{K}$ for $p = 3$.} \label{figure:1_ex_3_3} 
\end{subfigure}%
\hfill
\begin{subfigure}[t]{0.49\textwidth}
\centering
\begin{tikzpicture}[font=\footnotesize]
\begin{axis}[%
width=1.0\textwidth,
height=0.6\textwidth,
xmin=0,
xmax=6,
xlabel={$k$},
xtick = {1, 2, 3, 4, 5},
xticklabels={0, 1, 2, 3, 4},
ymode=log,
ymin=500,
ymax=300000,
ylabel={$n_{\text{bytes}}$},
yminorticks=true,
ymajorgrids=true,
ytick = {1e3, 1e4, 1e5}
]
\addplot [color=mycolor5, mark=star, mark size=1.0pt, mark options={solid, mycolor5}]
table[row sep=crcr]{%
1	8720\\
2	20072\\
3	38528\\
4	65816\\
5	103664\\
};
\addplot [color=mycolor1, mark=square, mark size=1.0pt, mark options={solid, mycolor1}]
table[row sep=crcr]{%
1	2202\\
2	2322\\
3	2442\\
4	2562\\
5	3234\\
};
\addplot [color=mycolor2, mark=square, mark size=1.0pt, mark options={solid, mycolor2}]
table[row sep=crcr]{%
1	2314\\
2	2322\\
3	2850\\
4	3042\\
5	3234\\
};
\addplot [color=mycolor3, mark=square, mark size=1.0pt, mark options={solid, mycolor3}]
table[row sep=crcr]{%
1	2842\\
2	2658\\
3	2850\\
4	3354\\
5	2922\\
};
\addplot [color=mycolor1, mark=triangle, mark size=1.0pt, mark options={solid, mycolor1}]
table[row sep=crcr]{%
1	3298\\
2	3762\\
3	5130\\
4	5146\\
5	9666\\
};
\addplot [color=mycolor2, mark=triangle, mark size=1.0pt, mark options={solid, mycolor2}]
table[row sep=crcr]{%
1	4450\\
2	4354\\
3	5610\\
4	3666\\
5	7386\\
};
\addplot [color=mycolor3, mark=triangle, mark size=1.0pt, mark options={solid, mycolor3}]
table[row sep=crcr]{%
1	5842\\
2	8610\\
3	15042\\
4	19674\\
5	15354\\
};
\end{axis}
\end{tikzpicture}
\caption{Memory of $\mathbf{y}$ for $p = 3$.} \label{figure:1_ex_3_4}
\end{subfigure} \\
\begin{subfigure}[t]{0.49\textwidth}
\centering
\begin{tikzpicture}[font=\footnotesize]
\begin{axis}[%
width=1.0\textwidth,
height=0.6\textwidth,
xmin=0,
xmax=4,
xlabel={$k$},
xtick = {1, 2, 3},
xticklabels={0, 1, 2},
ymode=log,
ymin=-10,
ymax=2000,
yminorticks=true,
ylabel={$s$},
yminorticks=true,
ymajorgrids=true,
ytick = {1, 10, 100, 1000}
]
\addplot [color=mycolor5, mark=star, mark size=1.0pt, mark options={solid, mycolor5}]
table[row sep=crcr]{%
1	80.867162\\
2	191.684059\\
3	373.107632\\
};
\addplot [color=mycolor1, mark=square, mark size=1.0pt, mark options={solid, mycolor1}]
table[row sep=crcr]{%
1	0.283979\\
2	0.394862\\
3	0.564576\\
};
\addplot [color=mycolor2, mark=square, mark size=1.0pt, mark options={solid, mycolor2}]
table[row sep=crcr]{%
1	0.837763\\
2	1.936993\\
3	4.115756\\
};
\addplot [color=mycolor3, mark=square, mark size=1.0pt, mark options={solid, mycolor3}]
table[row sep=crcr]{%
1	1.802386\\
2	3.910988\\
3	8.624838\\
};
\addplot [color=mycolor1, mark=triangle, mark size=1.0pt, mark options={solid, mycolor1}]
table[row sep=crcr]{%
1	0.41804\\
2	0.98172\\
3	1.516326\\
};
\addplot [color=mycolor2, mark=triangle, mark size=1.0pt, mark options={solid, mycolor2}]
table[row sep=crcr]{%
1	0.987641\\
2	2.002253\\
3	2.345259\\
};
\addplot [color=mycolor3, mark=triangle, mark size=1.0pt, mark options={solid, mycolor3}]
table[row sep=crcr]{%
1	9.278867\\
2	15.142992\\
3	21.853264\\
};
\end{axis}
\end{tikzpicture}
\caption{Comp. time for $p = 5$.} \label{figure:1_ex_5_1}
\end{subfigure}
\hfill
\begin{subfigure}[t]{0.49\textwidth}
\centering
\begin{tikzpicture}[font=\footnotesize]
\begin{axis}[%
width=1.0\textwidth,
height=0.6\textwidth,
xmin=0,
xmax=4,
xlabel={$k$},
xtick = {1, 2, 3},
xticklabels={0, 1, 2},
ymode=log,
ymin=1e-11,
ymax=5e-5,
yminorticks=true,
ylabel={$\left\Vert \cdot \right\Vert_{L^2}$},
yminorticks=true,
ymajorgrids=true,
ytick = {1e-10, 1e-9, 1e-8, 1e-7, 1e-6, 1e-5}
]
\addplot [color=mycolor5, mark=star, mark size=1.0pt, mark options={solid, mycolor5}]
table[row sep=crcr]{%
1	3.50820102615985e-09\\
2	5.26188409793613e-10\\
3	1.29618409647307e-10\\
};
\addplot [color=mycolor1, mark=square, mark size=1.0pt, mark options={solid, mycolor1}]
table[row sep=crcr]{%
1	1.47034870124211e-06\\
2	1.50033121246615e-06\\
3	9.40148553245068e-07\\
};
\addplot [color=mycolor2, mark=square, mark size=1.0pt, mark options={solid, mycolor2}]
table[row sep=crcr]{%
1	7.0685333025166e-09\\
2	6.64724009102415e-09\\
3	8.61140050642347e-09\\
};
\addplot [color=mycolor3, mark=square, mark size=1.0pt, mark options={solid, mycolor3}]
table[row sep=crcr]{%
1	3.51243470754565e-09\\
2	8.60650157606781e-09\\
3	2.8216214176269e-09\\
};
\addplot [color=mycolor1, mark=triangle, mark size=1.0pt, mark options={solid, mycolor1}]
table[row sep=crcr]{%
1	1.02439642129711e-05\\
2	2.88107808959581e-06\\
3	2.83802709526462e-06\\
};
\addplot [color=mycolor2, mark=triangle, mark size=1.0pt, mark options={solid, mycolor2}]
table[row sep=crcr]{%
1	8.37553406325938e-08\\
2	6.77978376204489e-08\\
3	3.7012577200797e-08\\
};
\addplot [color=mycolor3, mark=triangle, mark size=1.0pt, mark options={solid, mycolor3}]
table[row sep=crcr]{%
1	5.02480429576024e-09\\
2	1.8795782439686e-09\\
3	1.0049109309339e-09\\
};
\end{axis}
\end{tikzpicture}
\caption{$L^2$ error for $p = 5$.} \label{figure:1_ex_5_2}
\end{subfigure} \\
\begin{subfigure}[t]{0.49\textwidth}
\centering
\begin{tikzpicture}[font=\footnotesize]
\begin{axis}[%
width=1.0\textwidth,
height=0.6\textwidth,
xmin=0,
xmax=4,
xlabel={$k$},
xtick = {1, 2, 3},
xticklabels={0, 1, 2},
ymode=log,
ymin=5000,
ymax=5e8,
yminorticks=true,
ylabel={$n_{\text{bytes}}$},
yminorticks=true,
ymajorgrids=true,
ytick = {1e4, 1e5, 1e6, 1e7, 1e8}
]
\addplot [color=mycolor5, mark=star, mark size=1.0pt, mark options={solid, mycolor5}]
table[row sep=crcr]{%
1	18013064\\
2	44420336\\
3	88698008\\
};
\addplot [color=mycolor1, mark=*, mark size=1.0pt, mark options={solid, mycolor1}]
table[row sep=crcr]{%
1	35364\\
2	52212\\
3	72804\\
};
\addplot [color=mycolor2, mark=*, mark size=1.0pt, mark options={solid, mycolor2}]
table[row sep=crcr]{%
1	41980\\
2	61948\\
3	86268\\
};
\addplot [color=mycolor3, mark=*, mark size=1.0pt, mark options={solid, mycolor3}]
table[row sep=crcr]{%
1	41980\\
2	61948\\
3	86268\\
};
\end{axis}
\end{tikzpicture}
\caption{Memory of $\mathbf{K}$ for $p = 5$.} \label{figure:1_ex_5_3} 
\end{subfigure}%
\hfill
\begin{subfigure}[t]{0.49\textwidth}
\centering
\begin{tikzpicture}[font=\footnotesize]
\begin{axis}[%
width=1.0\textwidth,
height=0.6\textwidth,
xmin=0,
xmax=4,
xlabel={$k$},
xtick = {1, 2, 3},
xticklabels={0, 1, 2},
ymode=log,
ymin=1000,
ymax=1e5,
yminorticks=true,
ylabel={$n_{\text{bytes}}$},
yminorticks=true,
ymajorgrids=true,
ytick = {1e3, 1e4, 1e5}
]
\addplot [color=mycolor5, mark=star, mark size=1.0pt, mark options={solid, mycolor5}]
table[row sep=crcr]{%
1	13536\\
2	27960\\
3	50256\\
};
\addplot [color=mycolor1, mark=square, mark size=1.0pt, mark options={solid, mycolor1}]
table[row sep=crcr]{%
1	3058\\
2	3610\\
3	4258\\
};
\addplot [color=mycolor2, mark=square, mark size=1.0pt, mark options={solid, mycolor2}]
table[row sep=crcr]{%
1	5866\\
2	7298\\
3	8890\\
};
\addplot [color=mycolor3, mark=square, mark size=1.0pt, mark options={solid, mycolor3}]
table[row sep=crcr]{%
1	6346\\
2	15434\\
3	19994\\
};
\addplot [color=mycolor1, mark=triangle, mark size=1.0pt, mark options={solid, mycolor1}]
table[row sep=crcr]{%
1	11034\\
2	17186\\
3	21754\\
};
\addplot [color=mycolor2, mark=triangle, mark size=1.0pt, mark options={solid, mycolor2}]
table[row sep=crcr]{%
1	11034\\
2	14730\\
3	21882\\
};
\addplot [color=mycolor3, mark=triangle, mark size=1.0pt, mark options={solid, mycolor3}]
table[row sep=crcr]{%
1	16938\\
2	27922\\
3	46442\\
};
\end{axis}
\end{tikzpicture}
\caption{Memory of $\mathbf{y}$ for $p = 5$.} \label{figure:1_ex_5_4}
\end{subfigure}\\
\begin{subfigure}[t]{0.49\textwidth}
\centering
\begin{tikzpicture}[font=\footnotesize]
\begin{axis}[%
width=1.0\textwidth,
height=0.6\textwidth,
xmin=0,
xmax=6,
xlabel={$k$},
xtick = {1, 2, 3, 4, 5},
xticklabels={0, 1, 2, 3, 4},
ymin=-10,
ymax=60,
yminorticks=true,
ylabel={$n_{\text{it}}$},
yminorticks=true,
ymajorgrids=true,
ytick = {10, 20, 30, 40, 50, 60}
]
\addplot [color=mycolor1, mark=square, mark size=1.0pt, mark options={solid, mycolor1}]
table[row sep=crcr]{%
1	5\\
2	5\\
3	5\\
4	5\\
5	5\\
};
\addplot [color=mycolor2, mark=square, mark size=1.0pt, mark options={solid, mycolor2}]
table[row sep=crcr]{%
1	7\\
2	7\\
3	7\\
4	8\\
5	8\\
};
\addplot [color=mycolor3, mark=square, mark size=1.0pt, mark options={solid, mycolor3}]
table[row sep=crcr]{%
1	9\\
2	9\\
3	9\\
4	10\\
5	11\\
};
\addplot [color=mycolor1, mark=triangle, mark size=1.0pt, mark options={solid, mycolor1}]
table[row sep=crcr]{%
1	8\\
2	10\\
3	11\\
4	13\\
5	17\\
};
\addplot [color=mycolor2, mark=triangle, mark size=1.0pt, mark options={solid, mycolor2}]
table[row sep=crcr]{%
1	21\\
2	22\\
3	21\\
4	22\\
5	24\\
};
\addplot [color=mycolor3, mark=triangle, mark size=1.0pt, mark options={solid, mycolor3}]
table[row sep=crcr]{%
1	54\\
2	53\\
3	51\\
4	51\\
5	50\\
};
\end{axis}
\end{tikzpicture}
\caption{Iterations for $p = 3$.} \label{figure:1_ex_3_5} 
\end{subfigure}%
\hfill
\begin{subfigure}[t]{0.49\textwidth}
\centering
\begin{tikzpicture}[font=\footnotesize]
\begin{axis}[%
width=1.0\textwidth,
height=0.6\textwidth,
xmin=0,
xmax=4,
xlabel={$k$},
xtick = {1, 2, 3},
xticklabels={0, 1, 2},
ymin=-50,
ymax=400,
yminorticks=true,
ylabel={$n_{\text{it}}$},
yminorticks=true,
ymajorgrids=true,
ytick = {0, 100, 200, 300}
]
\addplot [color=mycolor1, mark=square, mark size=1.0pt, mark options={solid, mycolor1}]
table[row sep=crcr]{%
1	6\\
2	6\\
3	6\\
};
\addplot [color=mycolor2, mark=square, mark size=1.0pt, mark options={solid, mycolor2}]
table[row sep=crcr]{%
1	11\\
2	13\\
3	16\\
};
\addplot [color=mycolor3, mark=square, mark size=1.0pt, mark options={solid, mycolor3}]
table[row sep=crcr]{%
1	17\\
2	21\\
3	26\\
};
\addplot [color=mycolor1, mark=triangle, mark size=1.0pt, mark options={solid, mycolor1}]
table[row sep=crcr]{%
1	18\\
2	39\\
3	41\\
};
\addplot [color=mycolor2, mark=triangle, mark size=1.0pt, mark options={solid, mycolor2}]
table[row sep=crcr]{%
1	45\\
2	53\\
3	50\\
};
\addplot [color=mycolor3, mark=triangle, mark size=1.0pt, mark options={solid, mycolor3}]
table[row sep=crcr]{%
1	335\\
2	344\\
3	328\\
};
\end{axis}
\end{tikzpicture}
\caption{Iterations for $p = 5$.} \label{figure:1_ex_5_5}
\end{subfigure} \\
\par\smallskip
\hspace{3em}%
\begin{subfigure}[t]{0.9\textwidth}
\centering
\begin{tikzpicture}[font=\footnotesize]
\begin{axis}[%
hide axis,
xmin=10,
xmax=50,
ymin=0,
ymax=0.4,
legend columns=4, 
legend style={draw=white!15!black,legend cell align=left, /tikz/every even column/.append style={column sep=0.1cm}}
]
\addlegendimage{mycolor1}
\addlegendentry{$\varepsilon = 10^{-3}$};
\addlegendimage{mycolor2}
\addlegendentry{$\varepsilon = 10^{-5}$};
\addlegendimage{mycolor3}
\addlegendentry{$\varepsilon = 10^{-7}$};
\addlegendimage{color=mycolor5, mark=star, mark size=1.0pt, mark options={solid, mycolor5}}
\addlegendentry{\textsc{GeoPDEs}};
\addlegendimage{color = black, mark=square, mark size=1.0pt, mark options={solid, black}}
\addlegendentry{\cref{eq:block_preconditioner_1}};
\addlegendimage{color = white, mark=square, mark size=1.0pt, mark options={solid, white}}
\addlegendentry{};
\addlegendimage{color = black, mark=triangle, mark size=1.0pt, mark options={solid, black}}
\addlegendentry{Jacobi};
\end{axis}
\end{tikzpicture}
\end{subfigure}
\caption{Results for approximating \cref{eq:sol_1} using the refinement scheme \cref{eq:ref_1}.} 
\label{figure:1_ex} 
\end{figure}
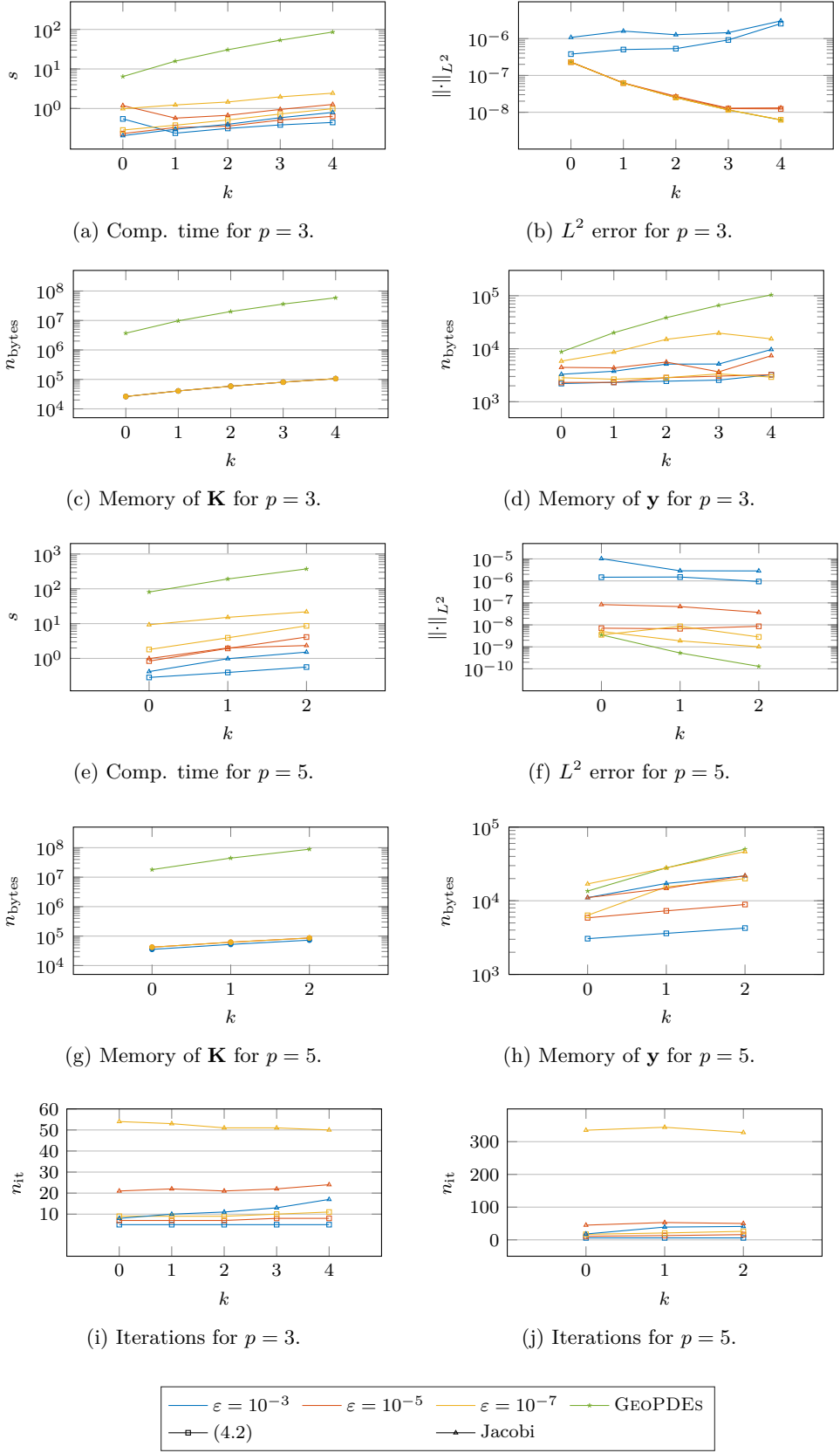

The results of the first experiment are shown in \cref{figure:1_ex}. Since the two approaches \cref{eq:approach_1_blockstructure} and \cref{eq:approach_2_non_diagonal_blocks} \& \cref{eq:approach_2_non_diagonal_blocks} are equivalent, we report only the results for \cref{eq:approach_1_blockstructure}. For $p=3$, the computational times for the block preconditioner remain modest over $k=0,\ldots,4$ (e.g., $\varepsilon=10^{-3}$: $0.55 s \to 0.44 s$, $\varepsilon=10^{-5}$: $0.23 s \to 0.63 s$, $\varepsilon=10^{-7}$: $0.29 s \to 0.99 s$), while \textsc{GeoPDEs} grows from $6.44s$ to $85.83 s$; the $L^2$-errors for $\varepsilon \le 10^{-5}$ essentially coincide with \textsc{GeoPDEs}, staying around $10^{-8}$--$10^{-7}$, whereas $\varepsilon=10^{-3}$ is noticeably less accurate. Since the memory footprint of $\mathbf{K}$ is identical for both preconditioners, we show it in the following only once: It turns out that the memory requirements for this refinement are identical for all tolerances, for it scales from $2.62 \cdot 10^{4}$ to $1.06 \cdot 10^{5}$ bytes for low-rank versus $3.70 \cdot 10^{6}$ to $5.88 \cdot 10^{7}$ bytes for \textsc{GeoPDEs}. For $\mathbf{y}$, the block preconditioner stays in the few-kilobyte range (e.g., $2.20 \cdot 10^{3}$--$3.23 \cdot 10^{3}$ bytes), while Jacobi requires more memory (up to $1.54 \cdot 10^{4}$ bytes). Iteration counts confirm that the block preconditioner \cref{eq:block_preconditioner_1} is preferable (about $5$--$11$ iterations depending on $\varepsilon$) over Jacobi (about $8$--$54$). For $p=5$, low-rank approximation becomes more demanding. With $\varepsilon=10^{-7}$, the block preconditioner attains $L^2$-errors of order $10^{-9}$ (e.g., for $\varepsilon=10^{-7}$: $3.5 \cdot 10^{-9}$, $8.6 \cdot 10^{-9}$, $2.8 \cdot 10^{-9}$), which are close to but generally above the \textsc{GeoPDEs} reference ($3.5 \cdot 10^{-9}$, $5.3 \cdot 10^{-10}$, $1.3 \cdot 10^{-10}$); for $\varepsilon=10^{-3}$ the errors deteriorate ($10^{-6}$--$10^{-5}$). Still, times and storage remain smaller than \textsc{GeoPDEs}; for $\varepsilon=10^{-7}$ and $\varepsilon=10^{-5}$, the number of bytes required to store $\mathbf{K}$ matches $4.2 \cdot 10^{4}$, $6.19 \cdot 10^{4}$ and $8.63 \cdot 10^{4}$. For $\varepsilon=10^{-3}$, we have $3.54 \cdot 10^{4}$, $5.22 \cdot 10^{4}$ and $7.28 \cdot 10^{4}$, which shows that the loss of accuracy conserves memory requirements. All tolerances beat \textsc{GeoPDEs} $1.80 \cdot 10^{7}$--$8.87 \cdot 10^{7}$ bytes. The memory of the solution $\mathbf{y}$ under the block preconditioner stays below $2.0 \cdot 10^{4}$ bytes, whereas Jacobi can approach \textsc{GeoPDEs} (e.g., $4.64 \cdot 10^{4}$ vs $5.03 \cdot 10^{4}$ bytes at $k=2$). Iteration counts again favor the block preconditioner (about $6$--$26$) over Jacobi (about $18$--$335$). Overall, the results are good and they definitely back the theory for easy hierarchical refinement scenarios with a single interface.

%
%

The analytical solution of our second experiment is again \cref{eq:sol_1}. For each degree $p \in \left\{ 3, 5 \right\}$ we start from a uniform coarse mesh $G_0$ with $m^{(d)}_0 = p+1$ cells per parametric direction $d=1,2,3$. We then build a hierarchical THB-spline space with levels $L=1,\ldots,k$ by repeating the following operation: On level $\ell$ the subdomain $\Omega_{\ell + 1}$ is represneted by the cells
\begin{equation}
	\label{eq:ref_2}
	\mathcal{M}_\ell = \left\{ q_{\ell,\left( j^{(1)},j^{(2)},j^{(3)} \right)} \in G_{\ell}: \; j^{(1)} = 1, \ldots, \bigl\lfloor\tfrac{m^{(1)}_\ell}{2^{\ell}}\bigr\rfloor \right\},
\end{equation}
which means we apply single dyadic subdivision (factor $2$ in each direction) of the cells contained in $\mathcal{M}_\ell$ and then do the same for the next level. Equivalently, the refined subdomains form a nested sequence of slices
\begin{equation*}
	\Omega_\ell = [0,2^{-\ell}]\times[0,1]\times[0,1],\qquad \ell=0,\ldots,L-1.
\end{equation*}
We employ the associated truncated hierarchical basis. We consider $k=1,\ldots,5$ for $p=3$ and $k=1,2,3$ for $p=5$. A visualization of the hierarchical mesh for the case $p = 3$ and $L = 3$ is shown in \cref{fig:panel-b}. 

The results of the second experiment are shown in \cref{figure:2_ex}. We note that also here for the refinement heuristic \cref{eq:ref_2} the two approaches \cref{eq:approach_1_blockstructure} and \cref{eq:approach_2_non_diagonal_blocks} \& \cref{eq:approach_2_non_diagonal_blocks} are equivalent, therefore we report only the results for \cref{eq:approach_1_blockstructure}. For the case $p=3$, computational time grows modestly for the block preconditioner across $L=1\ldots5$ (e.g., $\varepsilon=10^{-3}$: $0.14 s \to 0.96 s$; $\varepsilon=10^{-7}$: $0.11 s \to 1.98 s$), while Jacobi is consistently slower (e.g., $\varepsilon=10^{-7}$: $0.30 s \to 10.67 s$), and \textsc{GeoPDEs} is orders of magnitude higher ($0.64 s \to 130.87 s$); the $L^2$ errors essentially coincide with \textsc{GeoPDEs} for $\varepsilon \le 10^{-5}$ and remain in the $10^{-6}$--$10^{-7}$ range, whereas a loose tolerance $\varepsilon=10^{-3}$ under Jacobi yields noticeably larger errors. The memory footprint of $\mathbf{K}$ is for all tolerances the same, it scales from $2.95 \cdot 10^{3}$ to $7.66 \cdot 10^{5}$ bytes for low-rank versus $1.96 \cdot 10^{5}$ to $7.42 \cdot 10^{7}$ bytes for \textsc{GeoPDEs}. For $\mathbf{y}$, the block preconditioner stays compact (e.g., at $L=5$: $6.9 \cdot 10^{3}$ bytes for $\varepsilon=10^{-7}$) while Jacobi is larger ($5.29 \cdot 10^{4}$ bytes). The iteration counts confirm that the block preconditioner \cref{eq:block_preconditioner_1} is preferable (up to $16$ iterations for $\varepsilon=10^{-7}$) over the standard Jacobi smoother (up to $87$). For $p=5$ in \cref{figure:2_ex}, the low-rank solver with the block preconditioner delivers reference level accuracy as soon as the truncation is tight. With $\varepsilon=10^{-7}$, the $L^2$-error essentially matches \textsc{GeoPDEs} across $L=1,2,3$ (e.g., $5.85 \cdot 10^{-9}\to3.62 \cdot 10^{-9}$ vs.\ $5.85 \cdot 10^{-9}\to3.50 \cdot 10^{-9}$), while fot $\varepsilon=10^{-3}$ it deteriorates to $10^{-6}$. Timings remain very favorable: for $\varepsilon=10^{-7}$ the block preconditioner grows from $0.14s$ to $11.05s$ as $L$ increases, compared with $18.06 s \to 331.84 s$ for \textsc{GeoPDEs} (Jacobi is slower: $0.90s \to 24.50s$); for $\varepsilon=10^{-5}$ and $10^{-3}$ the computational times are even smaller (e.g., $0.12s \to 3.40s$ and $0.15s \to 1.04s$, respectively). Storage shows the clearest advantage. Here, the number of bytes differs over $L=2,3$ for the three tolerances: for $10^{-3}$ we have $3.54 \cdot 10^{4}$ and $1.41 \cdot 10^{5}$, for $10^{-5}$ we have $4.2 \cdot 10^{4}$ and $1.58 \cdot 10^{5}$ and for $10^{-7}$ we have $4.2 \cdot 10^{4}$ and $1.84 \cdot 10^{5}$. At $L=2$, the stiffness memory $\mathbf{K}$ for \textsc{GeoPDEs} is $7.28 \cdot 10^{7}$ bytes, i.e., a roughly $400$ times magnitude. The solution memory $\mathbf{y}$ with the block preconditioner remains compact -- about $1.00 \cdot 10^{4}$--$1.88 \cdot 10^{4}$ bytes at $L=2$ -- whereas Jacobi can exceed the \textsc{GeoPDEs} value (up to $4.51 \cdot 10^{4}$ vs.\ $3.91 \cdot 10^{4}$ bytes). Iteration counts consistently favor the block preconditioner: at $L=2$ it needs $11$, $17$, $41$ iterations for $\varepsilon=10^{-3},10^{-5},10^{-7}$, respectively, while Jacobi requires $38$, $55$, $330$. Overall, even at the higher degree $p=5$, these results are good and they definitely back the theory for easy hierarchical refinement: by choosing $\varepsilon\le10^{-5}$ and using the block preconditioner, we obtain (near-)reference accuracy at a fraction of the time and with orders-of-magnitude smaller $\mathbf{K}$-storage.

%
%
\begin{figure}
\definecolor{mycolor1}{rgb}{0.00000,0.44700,0.74100}%
\definecolor{mycolor2}{rgb}{0.85000,0.32500,0.09800}%
\definecolor{mycolor3}{rgb}{0.92900,0.69400,0.12500}%
\definecolor{mycolor4}{rgb}{0.49400,0.18400,0.55600}%
\definecolor{mycolor5}{rgb}{0.46600,0.67400,0.18800}%
\begin{subfigure}[t]{0.49\textwidth}
\centering
\begin{tikzpicture}[font=\footnotesize]
\begin{axis}[%
width=1.0\textwidth,
height=0.6\textwidth,
xmin=0,
xmax=6,
xlabel={$L$},
xtick = {1, 2, 3, 4, 5},
xticklabels={1, 2, 3, 4, 5},
ymode=log,
ymin=0,
ymax=500,
ylabel={$s$},
yminorticks=true,
ymajorgrids=true,
ytick = {1, 10, 100}
]
\addplot [color=mycolor5, mark=star, mark size=1.0pt, mark options={solid, mycolor5}]
  table[row sep=crcr]{%
1	0.64155\\
2	1.907567\\
3	7.785761\\
4	32.497421\\
5	130.87182\\
};
\addplot [color=mycolor1, mark=square, mark size=1.0pt, mark options={solid, mycolor1}]
  table[row sep=crcr]{%
1	0.137373\\
2	0.247471\\
3	0.27034\\
4	0.508065\\
5	0.960143\\
};
\addplot [color=mycolor2, mark=square, mark size=1.0pt, mark options={solid, mycolor2}]
  table[row sep=crcr]{%
1	0.105806\\
2	0.257951\\
3	0.336442\\
4	0.627928\\
5	1.372444\\
};
\addplot [color=mycolor3, mark=square, mark size=1.0pt, mark options={solid, mycolor3}]
  table[row sep=crcr]{%
1	0.110468\\
2	0.21489\\
3	0.475897\\
4	0.867776\\
5	1.978062\\
};
\addplot [color=mycolor1, mark=triangle, mark size=1.0pt, mark options={solid, mycolor1}]
  table[row sep=crcr]{%
1	0.104398\\
2	0.189987\\
3	0.347526\\
4	0.590026\\
5	1.164875\\
};
\addplot [color=mycolor2, mark=triangle, mark size=1.0pt, mark options={solid, mycolor2}]
  table[row sep=crcr]{%
1	0.778028\\
2	0.547906\\
3	0.76569\\
4	1.512986\\
5	3.087217\\
};
\addplot [color=mycolor3, mark=triangle, mark size=1.0pt, mark options={solid, mycolor3}]
  table[row sep=crcr]{%
1	0.29519\\
2	0.913434\\
3	1.971516\\
4	4.482433\\
5	10.666334\\
};
\end{axis}
\end{tikzpicture}
\caption{Comp. time for $p = 3$.} \label{figure:2_ex_3_1}
\end{subfigure}
\hfill
\begin{subfigure}[t]{0.49\textwidth}
\centering
\begin{tikzpicture}[font=\footnotesize]
\begin{axis}[%
width=1.0\textwidth,
height=0.6\textwidth,
xmin=0,
xmax=6,
xlabel={$L$},
xtick = {1, 2, 3, 4, 5},
xticklabels={1, 2, 3, 4, 5},
ymode=log,
ymin=1e-7,
ymax=1e-5,
ylabel={$\left\Vert \cdot \right\Vert_{L^2}$},
yminorticks=true,
ymajorgrids=true,
ytick = {1e-7, 1e-6, 1e-5}
]
\addplot [color=mycolor5, mark=star, mark size=1.0pt, mark options={solid, mycolor5}]
  table[row sep=crcr]{%
1	2.49382570225397e-06\\
2	7.7027790959977e-07\\
3	7.73933091710532e-07\\
4	7.7365622875175e-07\\
5	7.73635743589918e-07\\
};
\addplot [color=mycolor1, mark=square, mark size=1.0pt, mark options={solid, mycolor1}]
  table[row sep=crcr]{%
1	2.49382218505743e-06\\
2	1.71701351422294e-06\\
3	1.19967287615358e-06\\
4	1.03294203709325e-06\\
5	9.85973113756366e-07\\
};
\addplot [color=mycolor2, mark=square, mark size=1.0pt, mark options={solid, mycolor2}]
  table[row sep=crcr]{%
1	2.49379648426266e-06\\
2	7.70116686868015e-07\\
3	7.73840682040688e-07\\
4	7.73527410911887e-07\\
5	7.73404561338595e-07\\
};
\addplot [color=mycolor3, mark=square, mark size=1.0pt, mark options={solid, mycolor3}]
  table[row sep=crcr]{%
1	2.49384913012812e-06\\
2	7.70321215308223e-07\\
3	7.73973975864604e-07\\
4	7.73697046649611e-07\\
5	7.73678740141647e-07\\
};
\addplot [color=mycolor1, mark=triangle, mark size=1.0pt, mark options={solid, mycolor1}]
  table[row sep=crcr]{%
1	3.51891192665318e-06\\
2	2.32109647352791e-06\\
3	2.32279853792807e-06\\
4	2.69369409649812e-06\\
5	2.62105715570243e-06\\
};
\addplot [color=mycolor2, mark=triangle, mark size=1.0pt, mark options={solid, mycolor2}]
  table[row sep=crcr]{%
1	2.49407300969905e-06\\
2	7.71819978446768e-07\\
3	7.74954190484389e-07\\
4	7.75427077475292e-07\\
5	7.74059265130066e-07\\
};
\addplot [color=mycolor3, mark=triangle, mark size=1.0pt, mark options={solid, mycolor3}]
  table[row sep=crcr]{%
1	2.49384901299224e-06\\
2	7.70323289829341e-07\\
3	7.73976972413417e-07\\
4	7.73700206110921e-07\\
5	7.73677535433369e-07\\
};
\end{axis}
\end{tikzpicture}
\caption{$L^2$ error for $p = 3$.} \label{figure:2_ex_3_2}
\end{subfigure} \\
\begin{subfigure}[t]{0.49\textwidth}
\centering
\begin{tikzpicture}[font=\footnotesize]
\begin{axis}[%
width=1.0\textwidth,
height=0.6\textwidth,
xmin=0,
xmax=6,
xlabel={$L$},
xtick = {1, 2, 3, 4, 5},
xticklabels={1, 2, 3, 4, 5},
ymode=log,
ymin=1000,
ymax=5e8,
yminorticks=true,
ylabel={$n_{\text{bytes}}$},
ymajorgrids=true,
ytick = {1e3, 1e4, 1e5, 1e6, 1e7, 1e8}
]
\addplot [color=mycolor5, mark=star, mark size=1.0pt, mark options={solid, mycolor5}]
  table[row sep=crcr]{%
1	195680\\
2	863184\\
3	4022944\\
4	17645808\\
5	74175040\\
};
\addplot [color=mycolor1, mark=*, mark size=1.0pt, mark options={solid, mycolor1}]
  table[row sep=crcr]{%
1	2953\\
2	15604\\
3	59601\\
4	209440\\
5	766081\\
};
\addplot [color=mycolor2, mark=*, mark size=1.0pt, mark options={solid, mycolor2}]
  table[row sep=crcr]{%
1	2953\\
2	15604\\
3	59601\\
4	209440\\
5	766081\\
};
\addplot [color=mycolor3, mark=*, mark size=1.0pt, mark options={solid, mycolor3}]
  table[row sep=crcr]{%
1	2953\\
2	15604\\
3	59601\\
4	209440\\
5	766081\\
};
\end{axis}
\end{tikzpicture}
\caption{Memory of $\mathbf{K}$ for $p = 3$.} \label{figure:2_ex_3_3} 
\end{subfigure}%
\hfill
\begin{subfigure}[t]{0.49\textwidth}
\centering
\begin{tikzpicture}[font=\footnotesize]
\begin{axis}[%
width=1.0\textwidth,
height=0.6\textwidth,
xmin=0,
xmax=6,
xlabel={$L$},
xtick = {1, 2, 3, 4, 5},
xticklabels={1, 2, 3, 4, 5},
ymode=log,
ymin=500,
ymax=300000,
yminorticks=true,
ylabel={$n_{\text{bytes}}$},
ymajorgrids=true,
ytick = {1e3, 1e4, 1e5, 1e6}
]
\addplot [color=mycolor5, mark=star, mark size=1.0pt, mark options={solid, mycolor5}]
  table[row sep=crcr]{%
1	1000\\
2	2744\\
3	9032\\
4	32856\\
5	125544\\
};
\addplot [color=mycolor1, mark=square, mark size=1.0pt, mark options={solid, mycolor1}]
  table[row sep=crcr]{%
1	1021\\
2	2082\\
3	3271\\
4	4716\\
5	6865\\
};
\addplot [color=mycolor2, mark=square, mark size=1.0pt, mark options={solid, mycolor2}]
  table[row sep=crcr]{%
1	1021\\
2	2274\\
3	3463\\
4	5212\\
5	6865\\
};
\addplot [color=mycolor3, mark=square, mark size=1.0pt, mark options={solid, mycolor3}]
  table[row sep=crcr]{%
1	1501\\
2	2858\\
3	4607\\
4	6044\\
5	9289\\
};
\addplot [color=mycolor1, mark=triangle, mark size=1.0pt, mark options={solid, mycolor1}]
  table[row sep=crcr]{%
1	1501\\
2	3162\\
3	5511\\
4	7292\\
5	10321\\
};
\addplot [color=mycolor2, mark=triangle, mark size=1.0pt, mark options={solid, mycolor2}]
  table[row sep=crcr]{%
1	1661\\
2	3890\\
3	7335\\
4	13276\\
5	28945\\
};
\addplot [color=mycolor3, mark=triangle, mark size=1.0pt, mark options={solid, mycolor3}]
  table[row sep=crcr]{%
1	1821\\
2	4042\\
3	8423\\
4	19780\\
5	52937\\
};
\end{axis}
\end{tikzpicture}
\caption{Memory of $\mathbf{y}$ for $p = 3$.} \label{figure:2_ex_3_4}
\end{subfigure} \\
\begin{subfigure}[t]{0.49\textwidth}
\centering
\begin{tikzpicture}[font=\footnotesize]
\begin{axis}[%
width=1.0\textwidth,
height=0.6\textwidth,
xmin=0,
xmax=4,
xlabel={$L$},
xtick = {1, 2, 3},
xticklabels={1, 2, 3},
ymode=log,
ymin=-10,
ymax=1000,
yminorticks=true,
ylabel={$s$},
ymajorgrids=true,
ytick = {1, 10, 100, 1000}
]
\addplot [color=mycolor5, mark=star, mark size=1.0pt, mark options={solid, mycolor5}]
  table[row sep=crcr]{%
1	18.057277\\
2	81.089038\\
3	331.842186\\
};
\addplot [color=mycolor1, mark=square, mark size=1.0pt, mark options={solid, mycolor1}]
  table[row sep=crcr]{%
1	0.149126\\
2	0.286463\\
3	1.043212\\
};
\addplot [color=mycolor2, mark=square, mark size=1.0pt, mark options={solid, mycolor2}]
  table[row sep=crcr]{%
1	0.120114\\
2	0.879375\\
3	3.403857\\
};
\addplot [color=mycolor3, mark=square, mark size=1.0pt, mark options={solid, mycolor3}]
  table[row sep=crcr]{%
1	0.139608\\
2	1.889488\\
3	11.051226\\
};
\addplot [color=mycolor1, mark=triangle, mark size=1.0pt, mark options={solid, mycolor1}]
  table[row sep=crcr]{%
1	0.158978\\
2	0.43668\\
3	1.606134\\
};
\addplot [color=mycolor2, mark=triangle, mark size=1.0pt, mark options={solid, mycolor2}]
  table[row sep=crcr]{%
1	0.281293\\
2	1.019064\\
3	2.905469\\
};
\addplot [color=mycolor3, mark=triangle, mark size=1.0pt, mark options={solid, mycolor3}]
  table[row sep=crcr]{%
1	0.898483\\
2	7.365012\\
3	24.498259\\
};
\end{axis}
\end{tikzpicture}
\caption{Comp. time for $p = 5$.} \label{figure:2_ex_5_1}
\end{subfigure}
\hfill
\begin{subfigure}[t]{0.49\textwidth}
\centering
\begin{tikzpicture}[font=\footnotesize]
\begin{axis}[%
width=1.0\textwidth,
height=0.6\textwidth,
xmin=0,
xmax=4,
xlabel={$L$},
xtick = {1, 2, 3},
xticklabels={1, 2, 3},
ymode=log,
ymin=1e-9,
ymax=1e-4,
yminorticks=true,
ylabel={$\left\Vert \cdot \right\Vert_{L^2}$},
yminorticks=true,
ymajorgrids=true,
ytick = {1e-9, 1e-8, 1e-7, 1e-6, 1e-5, 1e-4}
]
\addplot [color=mycolor5, mark=star, mark size=1.0pt, mark options={solid, mycolor5}]
  table[row sep=crcr]{%
1	5.85484067096534e-09\\
2	3.50820102615985e-09\\
3	3.50081412579466e-09\\
};
\addplot [color=mycolor1, mark=square, mark size=1.0pt, mark options={solid, mycolor1}]
  table[row sep=crcr]{%
1	5.8565596969868e-09\\
2	1.4520885991982e-06\\
3	1.59536919923959e-06\\
};
\addplot [color=mycolor2, mark=square, mark size=1.0pt, mark options={solid, mycolor2}]
  table[row sep=crcr]{%
1	5.85437983799648e-09\\
2	6.36044416995841e-09\\
3	9.45149265530046e-09\\
};
\addplot [color=mycolor3, mark=square, mark size=1.0pt, mark options={solid, mycolor3}]
  table[row sep=crcr]{%
1	5.85438171193443e-09\\
2	3.51235415105102e-09\\
3	3.62040881256985e-09\\
};
\addplot [color=mycolor1, mark=triangle, mark size=1.0pt, mark options={solid, mycolor1}]
  table[row sep=crcr]{%
1	1.8681600784613e-06\\
2	1.02360852999637e-05\\
3	4.37881807730667e-06\\
};
\addplot [color=mycolor2, mark=triangle, mark size=1.0pt, mark options={solid, mycolor2}]
  table[row sep=crcr]{%
1	4.08671171660808e-08\\
2	8.66345894897345e-08\\
3	7.6301504229656e-08\\
};
\addplot [color=mycolor3, mark=triangle, mark size=1.0pt, mark options={solid, mycolor3}]
  table[row sep=crcr]{%
1	5.93582884011939e-09\\
2	4.78216296109688e-09\\
3	4.18563165401021e-09\\
};
\end{axis}
\end{tikzpicture}
\caption{$L^2$ error for $p = 5$.} \label{figure:2_ex_5_2}
\end{subfigure} \\
\begin{subfigure}[t]{0.49\textwidth}
\centering
\begin{tikzpicture}[font=\footnotesize]
\begin{axis}[%
width=1.0\textwidth,
height=0.6\textwidth,
xmin=0,
xmax=4,
xlabel={$L$},
xtick = {1, 2, 3},
xticklabels={1, 2, 3},
ymode=log,
ymin=1000,
ymax=5e8,
yminorticks=true,
ylabel={$n_{\text{bytes}}$},
yminorticks=true,
ymajorgrids=true,
ytick = {1e3, 1e4, 1e5, 1e6, 1e7, 1e8}
]
\addplot [color=mycolor5, mark=star, mark size=1.0pt, mark options={solid, mycolor5}]
  table[row sep=crcr]{%
1	5261984\\
2	18013064\\
3	72784176\\
};
\addplot [color=mycolor1, mark=*, mark size=1.0pt, mark options={solid, mycolor1}]
  table[row sep=crcr]{%
1	6537\\
2	35364\\
3	140625\\
};
\addplot [color=mycolor2, mark=*, mark size=1.0pt, mark options={solid, mycolor2}]
  table[row sep=crcr]{%
1	6537\\
2	41980\\
3	158113\\
};
\addplot [color=mycolor3, mark=*, mark size=1.0pt, mark options={solid, mycolor3}]
  table[row sep=crcr]{%
1	6537\\
2	41980\\
3	183697\\
};
\end{axis}
\end{tikzpicture}
\caption{Memory of $\mathbf{K}$ for $p = 5$.} \label{figure:2_ex_5_3} 
\end{subfigure}%
\hfill
\begin{subfigure}[t]{0.49\textwidth}
\centering
\begin{tikzpicture}[font=\footnotesize]
\begin{axis}[%
width=1.0\textwidth,
height=0.6\textwidth,
xmin=0,
xmax=4,
xlabel={$L$},
xtick = {1, 2, 3},
xticklabels={1, 2, 3},
ymode=log,
ymin=500,
ymax=100000,
yminorticks=true,
ylabel={$n_{\text{bytes}}$},
ymajorgrids=true,
ytick = {1e3, 1e4, 1e5}
]
\addplot [color=mycolor5, mark=star, mark size=1.0pt, mark options={solid, mycolor5}]
  table[row sep=crcr]{%
1	5832\\
2	13536\\
3	39096\\
};
\addplot [color=mycolor1, mark=square, mark size=1.0pt, mark options={solid, mycolor1}]
  table[row sep=crcr]{%
1	1117\\
2	3578\\
3	10015\\
};
\addplot [color=mycolor2, mark=square, mark size=1.0pt, mark options={solid, mycolor2}]
  table[row sep=crcr]{%
1	1117\\
2	5866\\
3	13055\\
};
\addplot [color=mycolor3, mark=square, mark size=1.0pt, mark options={solid, mycolor3}]
  table[row sep=crcr]{%
1	1117\\
2	6762\\
3	18823\\
};
\addplot [color=mycolor1, mark=triangle, mark size=1.0pt, mark options={solid, mycolor1}]
  table[row sep=crcr]{%
1	3349\\
2	11034\\
3	23047\\
};
\addplot [color=mycolor2, mark=triangle, mark size=1.0pt, mark options={solid, mycolor2}]
  table[row sep=crcr]{%
1	4717\\
2	11034\\
3	25639\\
};
\addplot [color=mycolor3, mark=triangle, mark size=1.0pt, mark options={solid, mycolor3}]
  table[row sep=crcr]{%
1	4717\\
2	16362\\
3	45079\\
};
\end{axis}
\end{tikzpicture}
\caption{Memory of $\mathbf{y}$ for $p = 5$.} \label{figure:2_ex_5_4}
\end{subfigure}\\
\begin{subfigure}[t]{0.49\textwidth}
\centering
\begin{tikzpicture}[font=\footnotesize]
\begin{axis}[%
width=1.0\textwidth,
height=0.6\textwidth,
xmin=0,
xmax=6,
xlabel={$L$},
xtick = {1, 2, 3, 4, 5},
xticklabels={1, 2, 3, 4, 5},
ymin=-10,
ymax=100,
yminorticks=true,
ylabel={$n_{\text{it}}$},
yminorticks=true,
ymajorgrids=true,
ytick = {10, 20, 30, 40, 50, 60, 70, 80, 90}
]
\addplot [color=mycolor1, mark=square, mark size=1.0pt, mark options={solid, mycolor1}]
  table[row sep=crcr]{%
1	1\\
2	4\\
3	5\\
4	7\\
5	8\\
};
\addplot [color=mycolor2, mark=square, mark size=1.0pt, mark options={solid, mycolor2}]
  table[row sep=crcr]{%
1	1\\
2	6\\
3	8\\
4	10\\
5	13\\
};
\addplot [color=mycolor3, mark=square, mark size=1.0pt, mark options={solid, mycolor3}]
  table[row sep=crcr]{%
1	1\\
2	9\\
3	12\\
4	14\\
5	16\\
};
\addplot [color=mycolor1, mark=triangle, mark size=1.0pt, mark options={solid, mycolor1}]
  table[row sep=crcr]{%
1	5\\
2	8\\
3	10\\
4	11\\
5	12\\
};
\addplot [color=mycolor2, mark=triangle, mark size=1.0pt, mark options={solid, mycolor2}]
  table[row sep=crcr]{%
1	15\\
2	22\\
3	25\\
4	27\\
5	29\\
};
\addplot [color=mycolor3, mark=triangle, mark size=1.0pt, mark options={solid, mycolor3}]
  table[row sep=crcr]{%
1	29\\
2	56\\
3	77\\
4	82\\
5	87\\
};
\end{axis}
\end{tikzpicture}
\caption{Iterations for $p = 3$.} \label{figure:2_ex_3_5} 
\end{subfigure}%
\hfill
\begin{subfigure}[t]{0.49\textwidth}
\centering
\begin{tikzpicture}[font=\footnotesize]
\begin{axis}[%
width=1.0\textwidth,
height=0.6\textwidth,
xmin=0,
xmax=4,
xlabel={$L$},
xtick = {1, 2, 3},
xticklabels={1, 2, 3},
ymin=-50,
ymax=400,
yminorticks=true,
ylabel={$n_{\text{it}}$},
yminorticks=true,
ymajorgrids=true,
ytick = {0, 100, 200, 300}
]
\addplot [color=mycolor1, mark=square, mark size=1.0pt, mark options={solid, mycolor1}]
  table[row sep=crcr]{%
1	1\\
2	6\\
3	11\\
};
\addplot [color=mycolor2, mark=square, mark size=1.0pt, mark options={solid, mycolor2}]
  table[row sep=crcr]{%
1	1\\
2	11\\
3	17\\
};
\addplot [color=mycolor3, mark=square, mark size=1.0pt, mark options={solid, mycolor3}]
  table[row sep=crcr]{%
1	1\\
2	17\\
3	41\\
};
\addplot [color=mycolor1, mark=triangle, mark size=1.0pt, mark options={solid, mycolor1}]
  table[row sep=crcr]{%
1	8\\
2	18\\
3	38\\
};
\addplot [color=mycolor2, mark=triangle, mark size=1.0pt, mark options={solid, mycolor2}]
  table[row sep=crcr]{%
1	23\\
2	45\\
3	55\\
};
\addplot [color=mycolor3, mark=triangle, mark size=1.0pt, mark options={solid, mycolor3}]
  table[row sep=crcr]{%
1	87\\
2	265\\
3	330\\
};
\end{axis}
\end{tikzpicture}
\caption{Iterations for $p = 5$.} \label{figure:2_ex_5_5}
\end{subfigure} \\
\par\smallskip
\hspace{3em}%
\begin{subfigure}[t]{0.9\textwidth}
\centering
\begin{tikzpicture}[font=\footnotesize]
\begin{axis}[%
hide axis,
xmin=10,
xmax=50,
ymin=0,
ymax=0.4,
legend columns=4, 
legend style={draw=white!15!black,legend cell align=left, /tikz/every even column/.append style={column sep=0.1cm}}
]
\addlegendimage{mycolor1}
\addlegendentry{$\varepsilon = 10^{-3}$};
\addlegendimage{mycolor2}
\addlegendentry{$\varepsilon = 10^{-5}$};
\addlegendimage{mycolor3}
\addlegendentry{$\varepsilon = 10^{-7}$};
\addlegendimage{color=mycolor5, mark=star, mark size=1.0pt, mark options={solid, mycolor5}}
\addlegendentry{\textsc{GeoPDEs}};
\addlegendimage{color = black, mark=square, mark size=1.0pt, mark options={solid, black}}
\addlegendentry{\cref{eq:block_preconditioner_1}};
\addlegendimage{color = white, mark=square, mark size=1.0pt, mark options={solid, white}}
\addlegendentry{};
\addlegendimage{color = black, mark=triangle, mark size=1.0pt, mark options={solid, black}}
\addlegendentry{Jacobi};
\end{axis}
\end{tikzpicture}
\end{subfigure}
\caption{Results for approximating \cref{eq:sol_1} using the refinement scheme \cref{eq:ref_2}.} 
\label{figure:2_ex} 
\end{figure}

The analytical solution of our third experiment is given by 
\begin{multline}
	\label{eq:sol_2}
	y_2 \left( x^{(1)}, x^{(2)}, x^{(3)} \right) = y_0 \left( x^{(1)}, x^{(2)}, x^{(3)} \right) \operatorname{exp} \Biggl( -10 \left( {x^{(1)}}^2 + {x^{(2)}}^2 + {x^{(3)}}^2 \right) \\
	\left( \left( {x^{(1)}} - 1 \right)^2 + \left( {x^{(2)}} - 1 \right)^2 + \left( {x^{(3)}} - 1 \right)^2 \right) \Biggr).
\end{multline}
For each degree $p \in \left\{ 3, 5 \right\}$ we start from a uniform coarse mesh with $m^{(d)}_0 = 2p+1$ cells per parametric direction $d=1,2,3$. We then build a hierarchical THB-spline space with levels $\ell=0,\ldots,L-1$ by repeating the following operation: On level $\ell$, the mesh $G_\ell$ has $m^{(d)}_\ell = 2^\ell m^{(d)}_0$ cells per direction. The subdomain $\Omega_{\ell + 1}$ is represneted by the cells $\mathcal{M}_\ell = \mathcal{M}_\ell^{(0)} \cup \mathcal{M}_\ell^{(1)}$, where
{\small
	\begin{equation}
		\label{eq:ref_3}
		\begin{aligned}
			\mathcal{M}_\ell^{(0)} &= \left\{  q_{\ell,\left( j^{(1)},j^{(2)},j^{(3)} \right)}\in G_\ell : \; 1 \le j^{(1)},j^{(2)},j^{(3)} \le p+\ell-1 \right\},\\
			\mathcal{M}_\ell^{(1)} &= \left\{  q_{\ell,\left( j^{(1)},j^{(2)},j^{(3)} \right)}\in G_\ell : \; m^{(d)}_\ell - (p+\ell-1) + 1 \le j^{(d)} \le m^{(d)}_\ell,\; d=1,2,3 \right\},
		\end{aligned}
\end{equation}}
which means we apply single dyadic subdivision (factor $2$ in each direction) of the cells contained in $\mathcal{M}_\ell$ and then do the same for the next level. Equivalently, defining $s_\ell = \frac{p+\ell-1}{m^{(d)}_\ell}$, the subdomains $\Omega_\ell = \Omega_\ell^{(0)} \cup \Omega_\ell^{(1)}$ on level $\ell$ are the two nested corner cubes
\begin{equation*}
	\Omega_\ell^{(0)} = \left[0,s_\ell \right] \times \left[0, s_\ell \right] \times \left[ 0, s_\ell \right], \quad \Omega_\ell^{(1)} = \left[ 1-s_\ell, 1 \right] \times \left[ 1-s_\ell, 1 \right] \times \left[1-s_\ell, 1 \right],
\end{equation*}
We employ the associated truncated hierarchical basis. We consider $k=0,\ldots,6$ for $p=3$ and $k=0,1,2,3$ for $p=5$. A visualization of the hierarchical mesh for the case $p = 3$ and $L = 3$ is shown in \cref{fig:panel-c}. 

The results of the third experiment are shown in \cref{figure:3_ex}. Because the refinement acts in two opposite corners, \cref{eq:approach_1_blockstructure} and \cref{eq:approach_2_non_diagonal_blocks} \& \cref{eq:approach_2_diagonal_blocks} are no longer equivalent; across all panels \cref{eq:approach_1_blockstructure} systematically dominates \cref{eq:approach_2_non_diagonal_blocks} \& \cref{eq:approach_2_diagonal_blocks}. In computational time, the block preconditioner with \cref{eq:approach_1_blockstructure} is only barely competitive with \textsc{GeoPDEs} at the finest levels: for $\varepsilon=10^{-3}$, runtimes grow from $0.31$ s to $46.09$ s (vs.\ $2.47\to130.60$ s); for $\varepsilon=10^{-7}$, from $0.47$ s to $127.49$ s (vs.\ $130.60$ s). Accuracy quickly saturates: the $L^2$-error tracks the reference and plateaus near $1.39 \cdot 10^{-7}$ from $L\ge3$, essentially independent of the preconditioner and only weakly dependent on $\varepsilon$. The clearest win is the stiffness storage $\mathbf{K}$: For \cref{eq:approach_1_blockstructure} it is far smaller than \textsc{GeoPDEs} across levels -- e.g., for $\varepsilon=10^{-3}$ it ranges from $5.45 \cdot 10^{3}$ to $1.25 \cdot 10^{7}$ bytes at $L=1\to7$, while \textsc{GeoPDEs} requires $1.37 \cdot 10^{6}\to7.18 \cdot 10^{7}$ bytes (\cref{figure:3_ex_3_3}). Solution storage $\mathbf{y}$ remains modest and grows gently with $L$: with \cref{eq:approach_1_blockstructure} it stays below $\approx1.1 \cdot 10^{5}$ bytes at $L=7$, whereas the dashed variant can reach $\approx2.5 \cdot 10^{5}$ bytes (\cref{figure:3_ex_3_4}). Iteration counts corroborate this picture: the block preconditioner is consistently better than Jacobi -- at $L=7$ one observes about $14,19,28$ iterations for $\varepsilon=10^{-3},10^{-5},10^{-7}$ with the block preconditioner, versus $20,42,90$ for Jacobi. In short, the approach \cref{eq:approach_1_blockstructure} is the superior choice; for $p=3$ the block preconditioner yields the best iteration behavior and strong $\mathbf{K}$-memory savings, even though wall-clock times at tight tolerances only just compete with \textsc{GeoPDEs}. For $p=5$ the problem is substantially more demanding and the differences between variants of $\mathbf{M}_{\ell, \ell'}$ and preconditioners are apparent. Across all panels, \cref{eq:approach_1_blockstructure} is consistently superior to \cref{eq:approach_2_non_diagonal_blocks} \& \cref{eq:approach_2_diagonal_blocks}: it is faster, uses markedly less $\mathbf{K}$-memory, and requires fewer iterations when it converges. At tight tolerance $\varepsilon=10^{-7}$, the Jacobi preconditioner with \cref{eq:approach_1_blockstructure} reaches the reference accuracy at the finest level we tested: at $L=4$ we obtain $1.254 \cdot 10^{-9}$, essentially identical to \textsc{GeoPDEs} $1.253 \cdot 10^{-9}$, while reducing computational time from $937.9s$ to $226.1s$ ($\approx4.1$ times faster). The same configuration also cuts the stiffness storage from $1.627 \cdot 10^{8}$ to $9.68 \cdot 10^{6}$ bytes ($\approx17$ times smaller); the solution vector stays modest at $1.21 \cdot 10^{5}$ bytes (vs.\ $8.36 \cdot 10^{4}$ for \textsc{GeoPDEs}). Relaxing the tolerance to $\varepsilon=10^{-5}$ preserves near-reference accuracy ($1.59 \cdot 10^{-9}$ at $L=4$) with much lower time ($60.4$ s) and still tiny $\mathbf{K}$-memory ($8.22 \cdot 10^{6}$ bytes). At $\varepsilon=10^{-3}$ the Jacobi path is fastest (e.g., $22.6s$ at $L=4$) but accuracy levels off around $5.84 \cdot 10^{-8}$. Preconditioner behavior is nuanced. For \cref{eq:approach_1_blockstructure}, Jacobi is the safer choice at tight tolerances: it converges at $\varepsilon=10^{-7}$ with $n_{\mathrm{it}}=76,325,328,414$ for $L=1,\ldots,4$, whereas the block preconditioner stalls and hits the cap $n_{\mathrm{it}}=900$ from $L=2$ onward (non-convergence). At looser tolerances the block preconditioner needs far fewer iterations (e.g., at $L=4$: $44$ for $\varepsilon=10^{-3}$ and $54$ for $\varepsilon=10^{-5}$ vs.\ Jacobi’s $41$ and $103$), but it still trails Jacobi in computational time at $p=5$. The approach \cref{eq:approach_2_non_diagonal_blocks} \& \cref{eq:approach_2_diagonal_blocks} is distinctly less attractive: it is much slower (e.g., $\varepsilon=10^{-7}$, Jacobi: $4.90 \cdot 10^{3}s$ at $L=4$), often inflates $\mathbf{K}$-memory by an order of magnitude at higher levels, and more frequently reaches the $900$-iteration cap. In summary, for $p=5$ the best performing and most reliable configuration is Jacobi preconditioner + \cref{eq:approach_1_blockstructure}. It matches \textsc{GeoPDEs} accuracy at $L=4$ with a $\sim4\times$ speedup and a $\sim17\times$ reduction in $\mathbf{K}$-storage, while the block preconditioner -- although efficient in iteration counts at relaxed tolerances -- can fail to converge at tight tolerances (hitting $900$ iterations).

\input{3_ex.tex}

The analytical solution of our fourth experiment is given by 
{\small
	\begin{multline}
		\label{eq:sol_3}
		y_3 \left( x^{(1)}, x^{(2)}, x^{(3)} \right) = y_0 \left( x^{(1)}, x^{(2)}, x^{(3)} \right) \operatorname{exp} \Biggl( - \left( {x^{(1)}}^2 + {x^{(2)}}^2 + {x^{(3)}}^2 \right) \\
		\left( \left( x^{(1)} - 1 \right)^2 + {x^{(2)}}^2 + {x^{(3)}}^2 \right) \left( {x^{(1)}}^2 + \left( {x^{(2)}} - 1 \right)^2 + {x^{(3)}}^2 \right) \left( {x^{(1)}}^2 + {x^{(2)}}^2 + \left( x^{(3)} - 1 \right)^2 \right) \Biggr).
\end{multline}}
For each degree $p \in \{3,5\}$ we start from a uniform coarse mesh with $m^{(d)}_0 = 2p+1$ cells per parametric direction $d=1,2,3$. We then build a hierarchical THB-spline space with levels $\ell=0,\ldots,k$ by repeating the following operation for $\ell=0,\ldots,k-1$. On level $\ell$, let $G_\ell$ denote the reference grid with $m^{(d)}_\ell$ cells per direction and cell indices $\mathbf{j} = \left( j^{(1)},j^{(2)},j^{(3)} \right)$. We mark four corner blocks,
{\small
	\begin{equation}
		\label{eq:ref_4}
		\begin{aligned}
			\mathcal{M}_\ell^{(0)} &= \left\{ q_{\ell, \mathbf{j}}\in G_\ell : \; 1 \le j^{(1)},j^{(2)},j^{(3)} \le p+\ell-1 \right\},\\
			\mathcal{M}_\ell^{(x)} &= \left\{ q_{\ell, \mathbf{j}}\in G_\ell : \; m^{(1)}_\ell-(p+\ell-1)+1 \le j^{(1)} \le m^{(1)}_\ell,\; 1 \le j^{(2)},j^{(3)} \le p+\ell-1 \right\},\\
			\mathcal{M}_\ell^{(y)} &= \left\{ q_{\ell, \mathbf{j}}\in G_\ell : \; 1 \le j^{(1)},j^{(3)} \le p+\ell-1,\; m^{(2)}_\ell-(p+\ell-1)+1 \le j^{(2)} \le m^{(2)}_\ell \right\},\\
			\mathcal{M}_\ell^{(z)} &= \left\{ q_{\ell, \mathbf{j}}\in G_\ell : \; 1 \le j^{(1)},j^{(2)} \le p+\ell-1,\; m^{(3)}_\ell-(p+\ell-1)+1 \le j^{(3)} \le m^{(3)}_\ell \right\},
		\end{aligned}
\end{equation}}
and create $G_{\ell+1}$ by dyadic subdivision (factor $2$ in each direction) of the marked set $\mathcal{M}_\ell = \mathcal{M}_\ell^{(0)} \cup \mathcal{M}_\ell^{(x)} \cup \mathcal{M}_\ell^{(y)} \cup \mathcal{M}_\ell^{(z)}$. Equivalently, defining $s_\ell = \dfrac{p+\ell-1}{m^{(1)}_\ell}$, the refined subdomains at level $\ell$ are the four corner cubes
\begin{align*}
	\Omega_\ell^{(0)} &= \left[0,s_\ell \right]\times \left[0,s_\ell \right]\times \left[0,s_\ell \right],\\
	\Omega_\ell^{(x)} &= \left[1-s_\ell,1 \right]\times \left[0,s_\ell \right]\times \left[0,s_\ell \right],\\
	\Omega_\ell^{(y)} &= \left[0,s_\ell \right]\times \left[1-s_\ell,1 \right]\times \left[0,s_\ell \right],\\
	\Omega_\ell^{(z)} &= \left[0,s_\ell \right]\times \left[0,s_\ell \right]\times \left[1-s_\ell,1 \right],
\end{align*}
We employ the associated truncated hierarchical basis. We consider $k=0,\ldots,6$ for $p=3$ and $k=0,1,2,3$ for $p=5$. A visualization of the hierarchical mesh for the case $p = 3$ and $L = 3$ is shown in \cref{fig:panel-d}. 

The results of the fourth experiment are shown in \cref{figure:4_ex}. Because refinement acts in all four corners, \cref{eq:approach_1_blockstructure} and \cref{eq:approach_2_non_diagonal_blocks} \& \cref{eq:approach_2_diagonal_blocks} are not equivalent. Throughout the panels, \cref{eq:approach_1_blockstructure} deliver, lower times, smaller memory, and fewer iterations than \cref{eq:approach_2_non_diagonal_blocks} \& \cref{eq:approach_2_diagonal_blocks}. For $p=3$, the block preconditioner with \cref{eq:approach_1_blockstructure} is competitive with \textsc{GeoPDEs}: computational time grows from $0.48s$ to $50.29s$ for $\varepsilon=10^{-3}$ and from $0.56s$ to $145.45s$ for $\varepsilon=10^{-7}$, versus $2.23s \to 249.09 s$ for \textsc{GeoPDEs}. The $L^2$-error quickly plateaus around $1.065 \cdot 10^{-6}$ and closely tracks the reference. We see that the storage of $\mathbf{K}$ depends on $\varepsilon$ and competes with \textsc{GeoPDEs} for \cref{eq:approach_1_blockstructure}. The memory requirements of $\mathbf{y}$ are similar with them of \textsc{GeoPDEs}. Iteration counts with the block preconditioner remain moderate ($13$--$26$ at large $L$, depending on $\varepsilon$) and are substantially lower than with Jacobi ($38$--$107$). For $p=5$, the problem is harder but still manageable, with a notable difference between preconditioners. First, we were only able to use \textsc{GeoPDEs} methods up to $L=2$ due to memory: the memory of $\mathbf{K}$ reaches for \textsc{GeoPDEs} $7.13 \cdot 10^{7}$ already by $L=2$. In contrast, the low-rank method remains light: at $L=2$ it uses $(1.21 - 1.84) \cdot 10^{6}$ bytes (across $\varepsilon$), and even at $L=4$ stays around $(9.09-14.05) \cdot 10^{6}$ bytes. In terms of time, the approach \cref{eq:approach_1_blockstructure} is faster than \textsc{GeoPDEs} at low levels (e.g., $L=1$: $0.38s$ vs.\ $110.75$ s; \cref{figure:4_ex_5_1}), but costs climb with $L$ and tighter $\varepsilon$. Most importantly, Jacobi performs better for $p=5$: it is more robust at tight tolerances and delivers the best accuracy -- e.g., with $\varepsilon=10^{-7}$ it matches the reference with $5 \cdot 10^{-9}$, whereas the block preconditioner often fails to converge, hitting the cap $n_{\mathrm{it}}=900$ already from $L=2$ at $\varepsilon=10^{-7}$ (and for the approach \cref{eq:approach_2_non_diagonal_blocks} \& \cref{eq:approach_2_diagonal_blocks} even earlier). Jacobi does converge across these cases, albeit with large iteration counts at tight tolerances (up to $6.3 \cdot 10^{2}$). Overall, even under the demanding four-corner refinement, the low-rank framework preferably with \cref{eq:approach_1_blockstructure} retains a consistent $\mathbf{K}$-memory advantage and delivers acceptable times and errors for $p=3$. For $p=5$, Jacobi is the preferred preconditioner (better robustness and accuracy), while the block preconditioner can stall at small $\varepsilon$. These results underscore the method’s usefulness when memory is the bottleneck or when moderately tight tolerances suffice.

\input{4_ex.tex}

\section{Conclusion, Limitations, and Outlook}
\label{section:conclusion}
We presented a first methodology for combining local refinement with low-rank techniques in adaptive IgA with THB-splines. The key idea is to recover separability where it still exists: interpolate geometry induced weights in separable spline spaces, assemble level-wise operators via univariate quadrature on cell cuboids, realize truncation in low rank, and accumulate everything in TT with systematic rounding. Our method shows that, for model problems in which the refined regions admit a partition into a small number of cuboids, the approach substantially reduces assembly time and memory while maintaining accuracy. In simple settings the two approaches \cref{eq:approach_1_blockstructure} and \cref{eq:approach_2_non_diagonal_blocks} \& \cref{eq:approach_2_diagonal_blocks} are equivalent and very efficient; for more complex refinements, \cref{eq:approach_1_blockstructure} remains clearly superior in memory and often competitive in time. As for solvers, the block-diagonal strategy tends to yield lower iteration counts for moderate degrees, whereas a Jacobi preconditioner was more robust for tighter tolerances and higher degrees.
The presented approach is not universally applicable. Its effectiveness hinges on two assumptions: the refined and non-refined sets (cells as well as THB-splines) can be partitioned into few cuboids and the interpolation ranks of weights and source terms remain moderate. When refinement patterns become geometrically intricate the number of cuboids and/or TT ranks can grow, inflating assembly cost and, at times, hampering solver robustness. Our greedy method to determine partitions into cuboids does not guarantee a minimal number of cuboids, so suboptimal partitions translate directly into more Kronecker summands which attack efficiency. Finally, the preconditioners used here (block-diagonal and Jacobi) were chosen for simplicity rather than optimality; iteration counts can rise with polynomial degree, number of levels, and tighter tolerances, and in isolated high-degree cases the block preconditioner may stall.
A natural next step is an adaptive low-rank preconditioner. In particular, a low-rank formulation of the BPX preconditioner for THB-splines~\cite{BRACCO2021113742} would likely be transformative and is, in our view, a key topic for future work. Beyond preconditioning, several other directions also look promising: improving the partitioning of active/non-active set to reduce the number of cuboids, rank-adaptive interpolation and rounding strategies that balance accuracy and storage automatically and extensions to multi-patch geometries while retaining local separability, to name a few.
We expect the proposed framework to be strong for PDE-constrained optimal control problems ore more generally for space-time formulations. Adding a time dimension typically introduces an additional Kronecker factor and preserves tensor-product structure; likewise, the KKT blocks inherit separable patterns. This extra structured dimension can compensate for complexity introduced by local refinement in space and may yield even larger gains than in the pure forward problem.
Overall, our results indicate that low-rank assembly and solution of hierarchical IgA systems are feasible and advantageous whenever the refinement is simple enough.

\end{document}